\documentclass[10pt, DIV14, twoside]{scrartcl}
\usepackage{latexsym,amsmath,amssymb,graphics,graphicx,amscd,amsfonts}
\usepackage[thmmarks,amsmath]{ntheorem}
\usepackage[arrow, matrix, curve]{xy}
\usepackage{hyperref}
\usepackage{eucal}

\newtheorem{theorem}{Theorem}
\newtheorem{definition}{Definition}
\newtheorem{lemma}{Lemma}
\newtheorem{proposition}{Proposition}
\newtheorem{corollary}{Corollary}
\theorembodyfont{\normalfont}
\theoremsymbol{\ensuremath{\lozenge}}
\newtheorem{example}{Example}
\theoremsymbol{\ensuremath{\triangle}}
\newtheorem{remark}{Remark}
\theoremstyle{nonumberplain}
\theoremsymbol{\ensuremath{\square}}
\theoremheaderfont{\scshape}
\theoremseparator{:}
\theorembodyfont{\normalfont}
\newtheorem{proof}{Proof}
\begin{document}
\title{Dirac structures, nonholonomic systems and reduction}
\author{M. Jotz\footnote{Section de Math{\'e}matiques, Ecole
  Polytechnique F{\'e}d{\'e}rale de Lausanne, 1015 Lausanne,
  Switzerland. madeleine.jotz@a3.epfl.ch
 This work was partially supported by Swiss
    NSF grant  200021-121512.}\,
 and T.S. Ratiu\footnote{Section de Math{\'e}matiques and Bernoulli Center, Ecole
  Polytechnique F{\'e}d{\'e}rale de Lausanne, 1015 Lausanne,
  Switzerland. tudor.ratiu@epfl.ch }
\date{}
}
\maketitle

\begin{abstract}
\centerline{\textbf{Abstract}} 
The reduction of nonholonomic systems is formulated in terms of Dirac reduction. 
 An optimal reduction method for a class of 
nonholonomic systems is formulated. Several examples are studied in detail.
\end{abstract}

\noindent \textbf{AMS Classification:} 53D20, 37J15, 70F25, 70H45, 37J60, 53D17

\noindent \textbf{Keywords:} Dirac structure, symmetry reduction, nonholonomic mechanics, momentum map.

\tableofcontents


\section{Introduction}

The equations of motion of nonholonomic mechanical systems and those in
circuit theory can be geometrically described using a Dirac structure
(introduced by \cite{Courant90a} and \cite{CoWe88}) in taking either a Hamiltonian or Lagrangian point
of view (see, e.g., \cite{Blankenstein00}, \cite{BlRa04},
\cite{BlvdS01},  \cite{MaYo06, MaYo07, MaYo09}). 
A Dirac structure $D$ on a manifold $M$ is a subbundle of the Pontryagin bundle $TM\oplus T^*M$
which is Lagrangian relative to the canonical symmetric pairing on it.
Dirac structures were introduced by \cite{CoWe88} and \cite{Courant90a} to provide a geometric
framework for the study of constrained mechanical systems. An easy example
of a Dirac structure is the graph of a $2$-form $\omega\in\Omega^2(M)$. 
Integrable Dirac structures have an additional integrability condition. 
They have been more intensively studied because they generalize, in a  certain sense, 
Poisson structures. For example, if the Dirac structure is the graph of $\omega\in\Omega^2(M)$, 
then it is integrable if and only if $\mathbf{d}\omega = 0 $. Other examples of integrable Dirac 
structures include various foliated manifolds. 
In general, an integrable Dirac structure determines a singular foliation on $M$ whose leaves carry a 
natural induced presymplectic structure.

Dirac
structures simultaneously generalize symplectic and Poisson structures and
also form the correct setting for the description of implicit Hamiltonian and
Lagrangian systems usually appearing as systems of algebraic-differential
equations. In symplectic and Poisson geometry, as well as geometric mechanics,
a major role is played by the reduction method since it creates, under
suitable hypotheses or in categories weaker than smooth manifolds, new spaces
with the same type of motion equations on them. Briefly put, it is a method
that eliminates variables and hence yields systems on smaller dimensional
manifolds. Due to the spectacular array of applications, reduction has been
extensively studied in various settings, including that of Dirac
manifolds. The present paper connects Dirac
and nonholonomic reduction, introduces an ``optimal'' reduction for special nonholonomic systems with symmetries
and presents several classical examples.

 First we recall the necessary background 
on Dirac geometry in \S\ref{diracquot}, as well as 
descriptions of Dirac reduction by symmetry groups. It is known that under certain assumptions 
beyond the usual ones, the quotient manifold carries a natural Dirac structure. These hypotheses are formulated 
in the literature in two different manners: using sections (see \cite{BlvdS01}) or appealing to the theory of 
fiber bundles (see \cite{BuCaGu07}).
We show in \S\ref{nonhol} that Dirac reduction as presented in \S\ref{red}, coincides 
with the method of reduction for nonholonomic systems
due to \cite{BaSn93}. This is achieved by reformulating their Hamiltonian approach to nonholonomic systems in the 
context of Dirac structures.
We give several standard examples that illustrate this.

We reformulate in \S\ref{sec:nonholonomic_Noether} the nonholonomic Noether Theorem 
(see \cite{BaSn93}, \S6, \cite{BaCuKeSn95}, Theorem 2,
and \cite{Bloch03} \footnote{A somewhat restricted version of the momentum equation
was given in \cite{KoKo78}; see also \cite{Arnold88}})
on the Hamiltonian side and show that, under the \emph{dimension assumption},
it is equivalent to the nonholonomic Noether Theorem in \cite{Bloch03}.
We study in \S\ref{sec:reac_annhili} a distribution where the fundamental vector fields have to lie
to yield constants of the motion. 
This gives an explanation for certain constants of motion that sometimes appear
as a consequence of the nonholonomic Noether theorem (see \cite{FaRaSa07}).

Under certain integrability assumptions imposed
 on a  distribution associated to a Dirac structure
modeling these, it is possible to extend the ideas in Marsden-Weinstein reduction 
to nonholonomic systems. 
This is achieved in \S\ref{nonholoptred}. These integrability conditions are certainly 
strong since they imply that the 
nonholonomic Noether $1$-forms that descend to the
quotient are exact. This is not true in general
but holds in the case of certain systems such as the vertical rolling disk or
the constrained particle. We discuss some of these examples at the end of the paper.

\medskip

\noindent\textbf{Conventions.} Throughout the paper $M$ is a \emph{paracompact} manifold, that is, it is
Hausdorff and every open covering admits a locally finite refinement.  The orientation preserving rotation group $\operatorname{SO}(2)$ of the plane $\mathbb{R}^2$ is also denoted by $\mathbb{S}^1$ and consists of matrices of the form
\[
\begin{bmatrix}
\cos \alpha & - \sin \alpha\\
\sin \alpha & \cos \alpha
\end{bmatrix}, \quad \alpha\in \mathbb{R}.
\]

\medskip

 If $E \rightarrow
M$ is a smooth fiber bundle over a manifold $M$, the spaces of smooth global
and local sections are denoted by $\Gamma_{{\rm global}}(E) $ and $\Gamma(E)$,
respectively. For example, $\mathfrak{X}(M): = \Gamma(T M) $ denotes the Lie algebra of smooth local vector fields endowed with the usual Jacobi-Lie bracket $[X,Y](f) = X[Y[f]] - Y[X[f]]$, where $X, Y \in \mathfrak{X}(M) $, $f$ is a smooth (possibly only locally defined) function on $M$, and $X[f]: = {\boldsymbol{\pounds}}_X f = \mathbf{d}f (X) $ denotes the Lie derivative of $f $ in the direction $X$. If $\wedge^k(M) \rightarrow M $ denotes the vector bundle of exterior $k$-forms on $M$ then $\Omega^k(M) : = \Gamma(\wedge^k(M)) $ is the space of local $k $-forms on the manifold $M$.
 
\medskip

Recall that a subset $ N \subset  M $ is an \emph{initial} submanifold 
of $ M $ if  $ N $ carries a manifold structure such that the inclusion 
$ \iota :N \hookrightarrow M$ is a smooth immersion and satisfies
the following condition: for any smooth manifold $ P $  an arbitrary
map $ g: P \rightarrow  N$ is smooth if and only if $\iota \circ  g: P 
\rightarrow  M$ is smooth. The notion of initial submanifold lies strictly between 
those of injectively immersed and embedded submanifolds.

\section{Background and Dirac structures}\label{diracquot}

This section  summarizes the key facts from the theory of foliations and Dirac
manifolds needed in the rest of the paper. It also establishes notation,
terminology, and conventions, since these are not uniform in the
literature. The proofs of the statements below can be found in 
\cite{Courant90a}, \cite{BlvdS01}, \cite{BlRa04}, \cite{BuCaGu07}. 

\subsection{Distributions and foliations}\label{sec:distributions}

We will need a few standard facts from the theory of generalized distributions on a smooth manifold $ M $ (see \cite{Stefan74a, Stefan74b, Stefan80}, \cite{Sussmann73} for the original articles and \cite{LiMa87}, \cite{Vaisman94}, \cite{Pflaum01}, or \cite{OrRa04},  for a quick review of this theory). 
\medskip

A \emph{generalized distribution} $ \Delta $ on $ M $ is a subset of the
tangent bundle $ TM $ such that $ \Delta (m) : = \Delta \cap T _m M$ is a
vector subspace of $ T _m M $. The number $\dim \Delta (m) $ is called the
\emph{rank}  of $ \Delta $ at $ m  \in  M $. A local \emph{differentiable
  section} of $ \Delta $ is a smooth vector field $ X \in  \mathfrak{X}(M) $
defined on some open subset $ U \subset  M $ such that $ X (u) \in  \Delta (u)
$ for each $ u \in  U $.  A generalized
distribution is said to be \emph{differentiable} or \emph{smooth} if for every
point $ m \in  M $ and every vector $ v \in  \Delta (m) $, there is a
differentiable section $ X \in  \Gamma ( \Delta ) $ defined on an open
neighborhood $ U $ of $ m $ such that $ X (m) = v $.

The term \emph{distribution} is usually synonymous to that of a vector subbundle of $ TM $. Since we shall work mostly with generalized distributions, we shall call below all generalized distributions simply distributions. If the generalized distribution happens to be a vector subbundle we shall always state this fact explicitly.

\medskip

In all that follows, $ \Delta $ is a smooth distribution. An \emph{integral manifold} of $ \Delta $ is an injectively immersed connected manifold $ \iota _L : L \hookrightarrow M $ satisfying the condition $ T _m \iota_L (T _m  L )  \subset  \Delta (m) $  for every $ m \in  L $. The integral manifold $ L $ is of \emph{maximal dimension} at $ m \in  L $ if $ T _m \iota_L (T _m  L )  =  \Delta (m) $. The distribution $ \Delta $ is \emph{completely integrable} if for every $m \in  M $ there is an integral manifold $ L $ of $ \Delta $, $ m \in  L $, everywhere of maximal dimension. The distribution $\Delta$ is \emph{involutive} if it is invariant under the (local) flows associated to differentiable sections of $\Delta$. The distribution $\Delta$ is \emph{algebraically involutive} if for any two smooth vector fields defined on an open set of $M $ which take values in $\Delta$, their bracket also takes values in $\Delta$. Clearly involutive distributions are algebraically involutive and the converse is true if the distribution is a subbundle. The analog of the Frobenius theorem (which deals only with vector subbundles of $ TM $) for distributions is known as the Stefan-Sussmann Theorem. Its statement is the same except that one needs the distribution to be involutive and not just algebraically involutive: \textit{$ \Delta $ is completely integrable if and only if $ \Delta $ is involutive}.

Recall that the Frobenius theorem states that a vector subbundle of $ TM $ is  (algebraically) involutive if and only if it is the tangent bundle of a foliation on $M $. The same is true for distributions: \textit{A smooth distribution is involutive if and only if it coincides with the set of vectors tangent to a generalized foliation.} To give content to this statement and elaborate on it, we need to quickly review the concept and main properties of generalized foliations. 
\medskip

A \emph{generalized foliation} on $ M $ is a partition $\mathfrak{F} : = 
\{ \mathcal{L} _\alpha \}_{ \alpha \in  A}$ of $ M $  into disjoint connected
sets, called \emph{leaves},  such that each point $ m \in  M $ has a
\emph{generalized foliated chart} $ (U, \varphi : U \rightarrow  V \in
\mathbb{R} ^{\dim M}) $, $ m \in  U $. This means that  there is some natural
number $ p_\alpha \leq \dim M $, called the \emph{dimension} of the leaf $ \mathcal{L} _\alpha $, and a subset $S _\alpha \subset  \mathbb{R} ^{\dim M - p_\alpha } $ such that $ \varphi(U \cap \mathcal{L} _\alpha ) = \{(x^1, \ldots ,x^{\dim M} ) \in  V \mid  (x^{p_\alpha+1}, \ldots  ,x^{\dim M} ) \in  S _\alpha \} $ and each $ (x^{p_\alpha+1}_ \circ  , \ldots  ,x^{\dim M}_ \circ  ) \in  S _\alpha$ determines a connected component $ ( U \cap \mathcal{L}_\alpha )_ \circ $ of $U \cap \mathcal{L}_\alpha$, that is, $ \varphi (( U \cap \mathcal{L}_\alpha )_ \circ ) = \{(x^1, \ldots , x^{p_\alpha}, x^{p_\alpha+1}_ \circ  , \ldots  ,x^{\dim M}_ \circ  ) \in V \} $. The key difference with the concept of foliation is that the number $ p_\alpha $ can change from leaf to leaf. The generalized foliated charts induce on each leaf a smooth manifold structure that makes them into initial submanifolds of $ M $. 

A leaf $ \mathcal{L} _\alpha $ is called \emph{regular} if it has an open neighborhood that intersects only leaves whose dimension equals $ \dim \mathcal{L} _\alpha $. If such a neighborhood does not exist, then $ \mathcal{L} _\alpha $ is called a \emph{singular} leaf. A point is called \emph{regular} (\emph{singular}) if it is contained in a regular (singular) leaf. The set of vectors tangent to the leaves of  $ \mathfrak{F} $ is defined by
\[
T(M,  \mathfrak{F}): = \bigcup _{ \alpha \in  A} \bigcup_{m \in  \mathcal{L} _\alpha} T _m \mathcal{L} _\alpha \subset  TM.
\]

Under mild topological conditions on $ M $ a generalized foliation has very useful properties. Assume that $ M $ is second countable. Then for each $ p_\alpha $-dimensional leaf $ \mathcal{L} _\alpha $  and any generalized foliated chart $(U, \varphi : U \rightarrow  V \in  \mathbb{R} ^{\dim M}) $ that intersects it, the corresponding set $ S _\alpha $ is countable. The set of regular points is open and dense in $ M$. Finally, any closed leaf is embedded in $ M $. Note that this last property is specific to (generalized) foliations since an injectively immersed submanifold whose range is closed is not necessarily embedded.

\medskip

Let us return now to the relationship between distributions and generalized foliations. As already mentioned, given an involutive (and hence a completely integrable) distribution $ \Delta $, each point $ m \in  M $ belongs to exactly one connected integral manifold $ \mathcal{L} _m $ that is maximal relative to inclusion. It turns out that $ \mathcal{L} _m $ is an initial submanifold and that it is also the \emph{accessible} set of $ m $, that is, $ \mathcal{L} _m $ equals the subset of points in $ M $ that can be reached by applying to $ m $ a finite number of composition of flows of elements of $\Gamma (\Delta)$.  The collection of all maximal integral submanifolds of $ \Delta $ forms a generalized foliation $ \mathfrak{F} _\Delta $ such that $ \Delta = T(M, \mathfrak{F} _\Delta) $.  Conversely, given a generalized foliation $ \mathfrak{F} $ on $ M $, the subset $ T(M, \mathfrak{F}) \subset  TM $ is a smooth completely integrable (and hence involutive) distribution whose collection of maximal integral submanifolds coincides with $ \mathfrak{F} $. These two statements expand the Stefan-Sussmann Theorem cited above. 
\medskip

In the study of Dirac manifolds we will also need the concept of codistribution. A \emph{generalized codistribution} $\Xi$ on $ M $ is a subset of the cotangent bundle $ T^*M $ such that $\Xi(m) : = \Xi \cap T^* _m M$ is a vector subspace of $ T^* _m M $. The notions of rank, differentiable section, and smooth codistribution are completely analogous to those for distributions. 
\medskip

If $\Delta \subset TM $ is a smooth distribution on $M$, its (\emph{smooth}) \emph{annihilator} $\Delta^ \circ $ is defined by 
\begin{align*}
\Delta^ \circ (m): = \left\{\alpha(m) \left|\begin{array}{c}  
\alpha  \in \Omega^1(M), \quad \left\langle \alpha, X \right\rangle = 0 
\;\text{for all}\; X \in \mathfrak{X}(U),\\m \in U \;\text{open such that}\; X(u) \in \Delta(u) \; \text{for all}\; u \in U\end{array}\right.\right\}.
\end{align*}
We have the, in general strict, inclusion $\Delta \subset \Delta^{ \circ \circ } $.
A similar definition holds for smooth codistributions. Note that the
annihilators are smooth by construction. If a distribution (codistribution) is
a vector subbundle of $TM $ (respectively of $T ^\ast M $), then its
annihilator is also a vector subbundle of $T ^\ast M $ (respectively of $TM
$). If  $\Delta $ is a subbundle then $\Delta = \Delta^{ \circ \circ } $ and similarly for codistributions.

\subsection{Dirac structures} \label{admsble}

For a smooth manifold $M$ denote by $\left\langle \cdot , \cdot  \right\rangle $ the duality pairing between the cotangent bundle $T ^\ast M $ and the tangent bundle $TM $ (or $\Omega^1(M) $ and $\mathfrak{X}(M)$). The \emph{Pontryagin bundle} $TM \oplus T^* M$ is endowed with a nondegenerate symmetric fiberwise bilinear form of signature $( \dim M, \dim M) $ given by
\begin{equation}
\label{pairing}
\left\langle (u_m, \alpha_m), ( v_m, \beta_m ) \right\rangle : = \left\langle\beta_m , u_ m \right\rangle + \left\langle\alpha_m, v _m \right\rangle
\end{equation}
for all $u _m, v _m \in T _mM$ and $\alpha_m, \beta_m \in T^\ast_mM$. A \emph{Dirac structure} (see \cite{Courant90a}) on $M $ is a Lagrangian subbundle $D \subset TM \oplus T^* M $, that is, $ D$ coincides with its orthogonal relative to \eqref{pairing} and so its fibers are necessarily $\dim M $-dimensional.

The space $\Gamma(TM \oplus T ^\ast M) $ of local sections of the Pontryagin bundle is endowed with a $\mathbb{R}$-bilinear skew-symmetric bracket (which does not satisfy the Jacobi identity) given by
\begin{align}
\label{courant_bracket}
[(X, \alpha), (Y, \beta) ] : &= \left( [X, Y],  {\boldsymbol{\pounds}}_X \beta - {\boldsymbol{\pounds}}_Y \alpha + \frac{1}{2} \mathbf{d}\left(\alpha(Y) - \beta(X) \right) \right) \nonumber \\
&= \left([X, Y],  {\boldsymbol{\pounds}}_X \beta - \mathbf{i}_Y \mathbf{d}\alpha - \frac{1}{2} \mathbf{d} \left\langle (X, \alpha), (Y, \beta) \right\rangle
\right)
\end{align}
(see \cite{Courant90a}). The Dirac structure is  \emph{integrable} if $[ \Gamma(D), \Gamma(D) ] \subset \Gamma(D) $. Since 
$\left\langle (X, \alpha), (Y, \beta) \right\rangle = 0 $ if $(X, \alpha), (Y, \beta) \in \Gamma(D)$, integrability of the Dirac structure is often expressed in the literature relative to a non-skew-symmetric bracket that differs from 
\eqref{courant_bracket} by eliminating in the second line the third term of
the second component. This truncated expression which satisfies the Jacobi identity but is no longer skew-symmetric is called the \emph{Courant-Dorfman bracket}
(see  \cite{BuCaGu07}, \cite{BuCrWeZh04}, \cite{BuCr05}, \cite{LiWeXu97},
\cite{SeWe01}).

\medskip

A Dirac structure defines two smooth distributions $\mathsf{G}_0, \mathsf{G}_1 \subset TM $ and two smooth codistributions $\mathsf{P}_0, \mathsf{P}_1 \subset T^*M$:
\begin{align*}
\mathsf{G}_0(m)&:= \{X(m) \in T_mM \mid X \in \mathfrak{X}(M), (X, 0) \in
\Gamma(D) \} \\
\mathsf{G}_1(m)&:= \{X(m) \in T_mM \mid X \in \mathfrak{X}(M), \, \text{there is an}\;
\alpha \in \Omega^1(M), \text{such that}\; 
(X, \alpha) \in \Gamma(D)\}
\end{align*}
and
\begin{align*}
\mathsf{P}_0(m)&:= \{\alpha(m) \in T^*_mM \mid \alpha \in \Omega^1(M), 
(0, \alpha) \in \Gamma(D) \} \\
\mathsf{P}_1(m)&:= \{\alpha(m) \in T^*_mM \mid  \alpha \in \Omega^1(M),\, \text{there is an}\;
X \in\mathfrak{X}(M), \text{such that}\; (X, \alpha) \in \Gamma(D) \}.
\end{align*}
The smoothness of $\mathsf{G}_0, \mathsf{G}_1, \mathsf{P}_0, \mathsf{P}_1$ is obvious since, by definition, they are generated by smooth local sections. In general, these are not vector subbundles of $TM$ and $T ^\ast M $, respectively. It is also clear that $\mathsf{G}_0 \subset \mathsf{G}_1$ and $\mathsf{P}_0 \subset \mathsf{P}_1$.

The \emph{characteristic equations} of a Dirac structure are
\begin{itemize}
\item[{\rm (i)}] $\mathsf{G}_0 = \mathsf{P}_1^\circ, \; \mathsf{P}_0 =
\mathsf{G}_1^\circ$.
\item[{\rm (ii)}] $\mathsf{P}_1 \subset \mathsf{G}_0^\circ, \; \mathsf{G}_1 \subset \mathsf{P}_0^\circ$.
\item[{\rm (iii)}] If $\mathsf{P}_1$ has constant rank, then $\mathsf{P}_1 = 
\mathsf{G}_0^\circ$. If $\mathsf{G}_1$ has constant rank, then $ \mathsf{G}_1 = \mathsf{P}_0^\circ$.
\end{itemize}

If $D$ is a Dirac structure on $M$ having  the property that $\mathsf{G}_1 \subset TM$ is a constant rank distribution on $M$, then there exists a skew-symmetric vector bundle map $\flat: \mathsf{G}_1 \rightarrow \mathsf{G}_1^*$  such that $D$ is given by
\begin{align}\label{gone}
D(m): = \{(X(m),  \alpha(m)) \in T_mM \oplus T_m^*M \mid X\, \text{a smooth local section of}\; \mathsf{P}^\circ, \alpha \in \Omega^1(M), \,
 \alpha|_{\mathsf{P}^\circ} = X^{\flat} \} \end{align}
 with $\mathsf{P} := \mathsf{P}_0 =
\mathsf{G}_1^\circ$. Also, $\ker(\flat:\mathsf{G}_1 \rightarrow 
\mathsf{G}_1^*) = \mathsf{G}_0$.
\medskip

A function $f\in C^\infty(M)$ is called \textit{admissible} if $\mathbf{d}f \in\Gamma(\mathsf{P}_1)$. 
There is an induced  bracket $\{\cdot,\cdot\}_D$ on the admissible functions given by
\begin{equation}\label{poissonbracket}
\{f,g\}_D=X_g[f]=-X_f[g] , 
\end{equation}
where $X_f\in\mathfrak{X}(M)$ is such that $(X_f,\mathbf{d} f)\in\Gamma(D)$.
 If the Dirac structure $D$ on $M$ is integrable, this bracket is a Poisson bracket.
Note that $X_f\in\mathfrak{X}(M)$ is not uniquely determined by this condition.
If the
Dirac structure is not integrable, we get with the same definition an almost
Poisson structure, that is, the Jacobi-identity doesn't necessarily hold.

\paragraph{Integrable Dirac structures as Lie algebroids}\label{algebroiddirac} 
A \emph{Lie algebroid} $E \rightarrow  M$ is a smooth vector bundle over $M$ with a vector bundle
homomorphism $\rho:E\to TM$, called the \emph{anchor},  and a Lie algebra bracket
$[\cdot,\cdot]:\Gamma(E)\times\Gamma(E)\to \Gamma(E)$ satisfying:
\begin{enumerate}
\item $\rho$ is a Lie algebra homomorphism
\item for all $f \in C^\infty(M)$ and $X$, $Y$ $\in \Gamma(E)$:
\[ [X, fY]=f[X,Y]+\rho(X)[f]Y.\]
\end{enumerate}
It is shown in \cite{Courant90a} that for an arbitrary Lie algebroid $ E \rightarrow  M$, the smooth distribution $\rho(E)$ is completely integrable. 
\medskip

Assume that $D $ is an integrable Dirac structure. Then, relative to the Courant bracket 
\eqref{courant_bracket} and the anchor $\pi_1:D\to TM$ given by the projection on the first factor, $ D $ 
becomes a Lie algebroid over $ M $.
The smooth distribution $ \pi _1( D ) \subset  TM $ coincides with $\mathsf{G}_1$. 
Hence, $\mathsf{G_1}$ is completely integrable and 
Theorem 2.3.6 in \cite{Courant90a} states the following result.

\begin{theorem}\label{presympleaves}
An integrable Dirac structure has a generalized foliation by presymplectic leaves.
\end{theorem} 

The presymplectic form $\omega_N$ on a leaf $N$ of the generalized foliation by leaves of $\mathsf{G_1}$ is
given by
\begin{equation}\label{induceddiracN}
\omega_N(\tilde{X},\tilde{Y})(p)=\alpha(Y)(p)=-\beta(X)(p)
\end{equation}
for all $p\in N$ and $\tilde{X},\tilde{Y}\in\mathfrak{X}(N)$, where $i_N:N\hookrightarrow M$ is the inclusion and $X,Y\in
\Gamma(\mathsf{G}_1)$ are $i_N$-related to  $\tilde{X},\tilde{Y}$, respectively; we shall denote $i_N$-relatedness by $\tilde{X}\sim_{i_N}X$ and
$\tilde{Y}\sim_{i_N}Y$. The $1$-forms $\alpha, \beta \in \Omega^1(M)$ are such that 
$(X,\alpha), (Y,\beta)\in \Gamma(D)$. Formula \eqref{induceddiracN} is independent of all the choices involved. 

\paragraph{Implicit Hamiltonian systems} 

Let $D$ be a Dirac structure on  $M $ and $H \in C ^{\infty}(M)$. The
\textit{implicit Hamiltonian system} $(M,D,H)$ is defined as the
set of $C^\infty$ solutions $x(t)$ satisfying the condition
\begin{equation}
\label{hamilton_equations}
(\dot{x},\mathbf{d}H(x(t))) \in D(x(t)), \quad \text{for all} \quad  t.
\end{equation}
In this general situation, \emph{conservation of energy} is still valid: 
$\dot{H}(t) =\langle \mathbf{d}H(x(t)), \dot{x}(t) \rangle = 0$, for all $t$ for which the solution exists. In addition, these equations contain algebraic constraints, namely, $\mathbf{d}H(x(t)) \in \mathsf{P}_1(x(t))$, for all $t$. Note that  $\dot{x}(t) \in \mathsf{G}_1(x(t))$, so the set of \emph{admissible flows}  have velocities in the distribution $\mathsf{G}_1$. Thus, an implicit Hamiltonian system defines a set of \emph{differential and algebraic equations}. 

Note that if $\mathsf{G}_1$ is an involutive subbundle of $TM$, then there are $\dim M - \operatorname{rank}\mathsf{G}_1$ independent
conserved quantities for the Hamiltonian system \eqref{hamilton_equations}. We
want to emphasize that standard existence and uniqueness theorems do not apply
to \eqref{hamilton_equations}, even if all the distributions and
codistributions are subbundles. The only general theorems that ensure the
local existence and  uniqueness of solutions for \eqref{hamilton_equations}
are for the so-called implicit Hamiltonian systems of index one (see \cite{Blankenstein00}, \cite{BlvdS01}).

\paragraph{Restriction of Dirac structures}\label{restriction} 
 Let $D$ be a Dirac structure on $M$ and $N\subset M$ a submanifold of $M$. 
Define the map $\sigma (m): T_mN \times T_m^\ast M \rightarrow T_mN \times T_m^\ast N, \; m \in
N$,  by $\sigma (m) (v_m,\alpha_m) = (v_m, \alpha_m\vert_{T_mN})$.
Assume that the dimension of $\mathsf{G}_1(m) \cap T_mN$ is independent of $m \in N$ and that the rank of
$\mathsf{G}_1$ is constant on $M$. Define the vector
subbundle $D_N \subset TN \oplus T^\ast N$ by
\[
  D_N(m)= \sigma(m) \left(D(m) \cap (T_mN \times T_m^\ast M)
\right), 
  \quad m \in N.
\]
Then $D_N$ is a Lagrangian subbundle in the Pontryagin bundle $TN\oplus T^\ast N $  and it is 
thus a Dirac structure on $N$. 
Let $\iota : N \hookrightarrow M$ denote the inclusion map and define for all $m\in N$ 
\begin{align*} 
E_s(m):=\left\{(X(m),\alpha(m))\in  T_mM\times T_m^*M\left|\begin{array}{c}
\alpha\in\Omega^1(M),\quad  X \in  \mathfrak{X}(M) \text{ such that }\\
  X(n) \in  T _n N \text{ for all } n \in  N \text{ for which } X 
\text{ is defined} \end{array}\right.\right\}
\end{align*}
(where the subscript $s$ stands for submanifold).
This defines a smooth bundle $E_s=\cup_{m\in N}E_s(m)$ on $N$.
\cite{BlvdS01} show that under the assumption that the fibers of $E_s\cap D$
have  constant dimension on $M$, there is another way to give the induced
Dirac structure, namely,  
$\left(\widetilde{X},\widetilde{\alpha}\right)$ is a local
  section of $D_N$ if and only if there exists a local section
  $(X,\alpha)$ of $D$ such that $\widetilde{X} \sim_\iota X$ and
  $\widetilde{\alpha} = \iota^\ast \alpha$.
Otherwise stated,
\begin{align}\label{induced_dirac_eq}
\Gamma(D_N) = \left\{ (\widetilde{X},\widetilde{\alpha}) \in 
   \mathfrak{X}(N) \oplus \Omega^1(N) \left|\begin{array}{c}
 \text{there is}\;
   (X,\alpha) \in \Gamma(D) \text{ such that }\\
    \widetilde{X} \sim_\iota X \text{ and } \widetilde{\alpha} = \iota^\ast \alpha 
\end{array}\right.\right\}.
\end{align}
Furthermore, if $D$ is integrable, then $D_N$ is also integrable. As stated in
\cite{BlvdS01}, if $\mathsf{G}_1$ is constant dimensional, the
assumptions for  both methods of restriction are equivalent.

\medskip

Second, we recall the restriction construction for implicit Hamiltonian systems. Given is the implicit Hamiltonian 
system $(M, D, H)$ and $N \subset M $ an invariant submanifold under the integral curves of $(M, D, H)$ 
(if they exist). Define $H_N: = H|_N = H \circ \iota$. Then every solution $x(t)$ of $(M,D,H)$ which 
leaves $N$ invariant (that is, $x(t)\in N$ for all $t$) is a solution of $(N,D_N,H_N)$. The converse statement 
is not true, in general.
 
\paragraph{Symmetries of Dirac manifolds} Let $G$ be a Lie group and
$\Phi: G\times M \rightarrow M$ a smooth left action. Then $G$ is called a
\emph{symmetry Lie group of} $D$ if for every $g\in G$ the condition
$(X,\alpha) \in \Gamma(D)$ implies that  $\left( \Phi_g^\ast X, \Phi_g^\ast
  \alpha \right) \in \Gamma(D)$. We say then that the Lie group $G$ acts
\textit{canonically} or \textit{by Dirac actions} on $M$.

For any admissible $f\in C^\infty(M)$, i.e., a function such that $(X_f,\mathbf{d}
f)\in \Gamma(D)$ for some $X_f \in  \mathfrak{X}(M)$, this yields
$(\Phi_g^*X_f,\Phi^*_g\mathbf{d}  f)\in
\Gamma(D)$ or $(\Phi_g^*X_f,\mathbf{d}  (\Phi^*_g f))\in \Gamma(D)$. Hence we have
simultaneously the facts that $\Phi^*_g f$ is admissible and that
$\Phi_g^*X_f-X_{\Phi^*_g f}=:Y\in\Gamma(\mathsf{G}_0)$. 
This implies for the almost Poisson bracket on admissible functions (see
\ref{admsble}):
\begin{align*}
\Phi_g^*\{f,h\}_D&=-\Phi_g^*(X_f[h])=-(\Phi_g^* X_f )[\Phi_g^*h]=-(Y+X_{\Phi^*_g f})[\Phi_g^*h]\\ 
&= -\mathbf{d} (\Phi_g^*h)(Y+X_{\Phi^*_g f})
= -\mathbf{d} (\Phi_g^*h)(X_{\Phi^*_g f})=\{\Phi^*_g f,\Phi_g^*h\}_D
\end{align*}
since $\Phi_g ^\ast h $ is an admissible function (and hence $ \mathbf{d} (\Phi_g ^\ast h) \in  
\Gamma(\mathsf{P}_1) \subset  \Gamma\left(\mathsf{G}_0^ \circ\right)$ and $ Y \in  \Gamma (\mathsf{G}_0)$).

The Lie group $G$ is a \emph{symmetry Lie group of the implicit Hamiltonian system}
$(M,D,H)$ if, in addition, $H$ is  $G$-invariant, that is, $H\circ \Phi_g = H$ for
all $g \in G$.

Let $\mathfrak{g}$ be a Lie algebra and $\xi \in \mathfrak{g} \mapsto \xi_M \in
\mathfrak{X}(M)$ be a smooth left Lie algebra action, that is, the map $(x, \xi) \in
M \times \mathfrak{g} \mapsto \xi_M(x) \in TM $ is smooth and $\xi \in\mathfrak{g}
\mapsto \xi_M  \in \mathfrak{X}(M)$ is a Lie algebra anti-homomorphism.  The Lie
algebra $\mathfrak{g}$ is said to be a 
\emph{symmetry Lie algebra of} $D$ if for every $\xi \in \mathfrak{g}$ the condition
$(X,\alpha) \in \Gamma(D)$ implies that  
$\left(\boldsymbol{\pounds}_{\xi_M}X,\boldsymbol{\pounds}_{\xi_M}\alpha \right) \in
\Gamma(D)$. If, in addition, $\boldsymbol{\pounds}_{\xi_M}H=0$ for all $\xi \in
\mathfrak{g}$, then $\mathfrak{g}$ is a \emph{symmetry Lie algebra of the implicit
Hamiltonian system} $(M, D, H) $. Of course, if $\mathfrak{g}$ is the Lie algebra of
$G $ and $\xi\mapsto \xi_M$ is the associated infinitesimal generator, then  $G $
is a symmetry Lie group  of $D$ if and only if  $\mathfrak{g}$ is a symmetry Lie
algebra of $D$.

\subsection{Regular reduction of Dirac structures}\label{red}
In all that follows we shall assume that
$(M, D)$ is a smooth Dirac manifold,
  $G $ is a symmetry Lie group of the Dirac structure $D$ 
on $M $ and that the action is free and proper. Thus, the projection on the quotient $\pi: M \rightarrow M/G:=\bar{M}$ 
defines a left principal $G $-bundle. Note that 
the Dirac structure $D \subset  TM \oplus T^\ast M$ is $G$-invariant as 
a subbundle since for all $g\in G$ and $(X,\alpha)\in \Gamma(D)$ we have 
$(\Phi_g^*X,\Phi_g^*\alpha)\in \Gamma(D)$. 
Set $\mathfrak{g}_M: =  \{\xi_M\mid \xi \in \mathfrak{g}\} \subset  \mathfrak{X} (M)$
and, for $m\in M$, define
the 
vector subspace
$\mathcal{V}(m): = \{\xi_M(m)\mid \xi\in \mathfrak{g}\} \subseteq T_m M $ and the 
distribution $\mathcal{V}:= \cup _{m \in M } \mathcal{V}(m)$. Since the $G$-action is free, $\mathcal{V}$ 
is a $G $-invariant vector subbundle of $TM $ and its annihilator $\mathcal V^\circ$ is a $G$-invariant subbundle of 
$T^*M$. 
 It is worth noting that the space of sections $\Gamma( \mathcal{V}) $ coincides with the 
$C ^{\infty}(M)$-module spanned by $\mathfrak{g}_M$.

\bigskip

For all $m\in M$ the map $T_m\pi:T_mM\to T_{\pi(m)}\bar M$ is 
surjective with kernel $\mathcal{V}(m)$. This yields an isomorphism between $T_mM/\mathcal{V}(m)$ and $T_{\pi(m)}\bar{M}$. 
The Lie 
group $G$ acts smoothly on the quotient vector bundle 
$TM/ \mathcal{V}$ by $ g \cdot  \hat v: = \widehat{T \Phi_g(v)} $, where 
$\hat v \in  TM/\mathcal{V}$.

For $X\in \mathfrak{X}(M)$, we will say that  the section
$\widehat{X}:=X(\operatorname{mod}\mathcal{V})$ of $TM/\mathcal{V}$ is 
\emph{$G$-equivariant}, if there is a representative $X^G$ of $\widehat{X}$ that is $G$-equivariant, i.e., a smooth section $X^G\in\mathfrak{X}(M)^G$ with
$X-X^G\in\Gamma(\mathcal{V})$. This is equivalent to the condition $[X,V]\in  \Gamma(\mathcal{V})$ 
for all representatives $X$ of  $\widehat{X}$ and for all $V\in\Gamma(\mathcal{V})$ (see for instance \cite{JoRaSn11}).
In what follows we shall use these two equivalent definitions interchangeably.
  
The representative  $X^G$ of $\widehat{X}$ uniquely
induces a smooth vector field $\bar{X} $ on $\bar{M} $, where $ \bar{X} $ is
defined by the condition $X^G \sim_ \pi \bar{X}$, that is, $T \pi \circ X^G =
\bar{X} \circ \pi$.
The map  
\begin{eqnarray}\label{iso-of-sections}
\begin{array}{cccc}
\Pi:&\Gamma(TM/\mathcal{V})^G&\to&\mathfrak{X}(\bar M)\\
&X(\operatorname{mod}\mathcal{V}) &\mapsto& \bar{X},
\end{array}
\end{eqnarray}
is a well defined homomorphism of
$C^\infty(\bar M)$-modules (note that $ C^\infty(\bar M)\simeq C^\infty(M)^G$
via $\bar f\mapsto \pi^*\bar f$). This map \eqref{iso-of-sections} is  in fact  an isomorphism
(use, for instance, the results in  \cite{Schwarz80}).

In the same way, for all $\bar{\alpha} \in \Omega^1(\bar M)$, we have
$\pi^*\bar{\alpha} \in \Gamma(\mathcal{V}^\circ)^G$. Note that if $\alpha \in 
\Gamma(\mathcal{V}^\circ)^G$, then the $1$-form $\bar{\alpha} \in \Omega^1(\bar{M}) $
defined by $\left\langle \bar{ \alpha}(\pi(m)), T_m \pi( v _m) \right\rangle :=
\left\langle \alpha(m), v_m \right\rangle$, for all $v _m\in T _mM $, is well defined and satisfies 
$\pi^*\bar{\alpha} = \alpha$. This shows that the map 
$\bar{ \alpha} \in \Omega^1(\bar {M}) \mapsto \pi^\ast \bar{\alpha} \in \Gamma( \mathcal{V}^ \circ)^G $ 
is an isomorphism of $C^\infty(\bar M)$-modules.

\medskip

We close these preliminary remarks by recording that the $G $-action on $(TM/ \mathcal{V}) \oplus \mathcal{V}^ \circ $
\[
g \cdot \left(\hat{v}_m, \alpha_m \right) : = \left(\widehat{T_m \Phi_g(v _m)}, T_{g \cdot m}^\ast  \Phi_{ g ^{-1}} \alpha_m \right)
\]
is free and proper.

\medskip 

Consider the vector subbundle $\mathcal{K}: = 
\mathcal{V} \oplus \{0\}
\subset TM \oplus T^* M $ of
the Pontryagin bundle and its orthogonal  $\mathcal{K}
^\perp = TM \oplus \mathcal{V}^ \circ $. Both vector subbundles are
$G$-invariant and it is easy to show
(in agreement with the more general results of \cite{BuCaGu07}) that 
\begin{equation}\label{Ered}
\left.\frac{\mathcal{K} ^\perp}{\mathcal{K}}\right/G \; =\; \left.\frac{TM \oplus
\mathcal{V}^ \circ}{\mathcal{V} \oplus \{0\}}\right/G\; =\;
\left.\frac{TM}{\mathcal{V}}\oplus\mathcal{V}^ \circ\right/G
\end{equation}
is a Courant algebroid over $\bar M$ with the symmetric bilinear $2$-form that
descends
from the one on $\mathcal{K}^\perp/\mathcal{K}$ given by
\begin{equation}\label{paring_k_k_perp}
\langle(\widehat{X},\alpha),(\widehat{Y},\beta)\rangle_{\mathcal{K}^\perp/\mathcal{K}
}=\beta(X)+\alpha(Y)
\end{equation} 
for all $\alpha$, $\beta$ in $\Gamma(\mathcal{V}^\circ)$ and
$X$, $Y$ in $\mathfrak{X}(M)$; here
$\widehat{X}:=X(\operatorname{mod}\mathcal{V})$,
$\widehat{Y}:=Y(\operatorname{mod}\mathcal{V})$ denote local sections 
of $TM/ \mathcal{V} $ induced by local vector fields on $M$. 

\medskip

We have used above the following general fact that will be needed also in later
arguments.  The proof is straightforward.

\begin{lemma}
Let $\pi:E\to M$ be a smooth vector bundle over $M$. Assume that there are two free
proper $G $-actions on $E$ and $M $,  respectively, such that $\pi$ is equivariant
and the action on $E$ is linear on the fibers.
Then the induced map $\pi_G:E/G\to M/G$ defined by the commutative diagram
\begin{eqnarray*}
\begin{CD}
E @>\pi>> M\\
 @V{\pi_E}VV   @VV{\pi_M}V \\
E/G @>\pi_G>> M/G
\end{CD}
\end{eqnarray*}
is also a smooth vector bundle whose rank is equal to the rank of $E $.
\end{lemma}

In fact, with the identifications given
above of $\Gamma(\mathcal{V}^\circ)^G$ with $\Omega^1(\bar M)$ and
$\Gamma(TM/\mathcal{V})^G$ with $\mathfrak{X}(\bar M)$, it is obvious that the
$G$-equivariant sections of (\ref{Ered}) are in one-to-one correspondence with
those of $T\bar M\oplus T^*\bar M$. Note that this says that we have a vector bundle
isomorphism
\begin{equation}\label{isocou}
\left.\frac{\mathcal{K}^\perp}{\mathcal{K}} \right/G\; \simeq \;T\bar{M} \oplus
T^*\bar{M} 
\end{equation} over $\bar M=M/G$.
This vector bundle isomorphism preserves the symmetric pairing 
and it is easy to check that the Courant bracket on $T\bar{M} \oplus 
T^*\bar{M}$ also descends from the Courant bracket on  $TM\oplus 
T^*M$.

\medskip

We summarize here the approach for Dirac reduction
in \cite{BuCaGu07}. Assuming that $D \cap \mathcal{K}^\perp $ has constant rank, that is, $D
\cap \mathcal{K}^\perp $ is a smooth vector subbundle of $TM \oplus T^* M $, it
follows that $(D \cap \mathcal{K}^\perp)^\perp = D + \mathcal{K}$ and $D \cap
\mathcal{K}$ are vector subbundles of $TM \oplus T ^\ast M $. 
Form the pointwise quotient
\begin{equation}\label{Dred}
\frac{(D\cap \mathcal{K}^\perp)+ \mathcal{K}}{\mathcal{K}}
= \frac{\left(D\cap (TM\oplus \mathcal{V}^\circ) \right)+ (\mathcal{V}\oplus
  \{0\})}{\mathcal{V}\oplus\{0\} }
\end{equation}
with base $M $. At each point $m\in M$, one gets a subspace of the vector space
$(T_mM/\mathcal{V}(m)) \oplus
\mathcal{V}^\circ(m)\simeq \mathcal{K}^\perp(m) /\mathcal{K}(m)$ (see 
\eqref{Ered}).

\begin{proposition} \label{tilde_D}
Relative to the symmetric nondegenerate bilinear form \eqref{paring_k_k_perp} on
$\mathcal{K}^\perp/\mathcal{K}$, the vector subspace 
\[
\tilde{D}(m):=\frac{(D(m)\cap \mathcal{K}(m)^\perp)+
\mathcal{K}(m)}{\mathcal{K}(m)}\quad \text{ of }\quad \frac{\mathcal{K}(m)^\perp}{\mathcal{K}(m)} \] 
satisfies $\tilde{D}(m)=\tilde{D}(m)^\perp$. 
\end{proposition}

\begin{proof}
Let us prove that $\tilde{D}(m)\subseteq\tilde{D}(m)^\perp$. Let
$(\widehat{X}(m),\alpha(m)) \in \tilde{D}(m)$. If  $(\widehat{X},\alpha)\in
\Gamma(\tilde{D})$ are local sections about $m$, then $\alpha\in
\Gamma(\mathcal{V}^\circ)$ and there are $X\in \mathfrak{X}(M)$ and 
$V\in\Gamma(\mathcal{V})$ such that
$(X+V,\alpha)\in \Gamma(D)$ and $\widehat{X}=X(\operatorname{mod} \mathcal{V})$. For all $(\widehat{Y},\beta)\in
\Gamma(\tilde{D})$ we have analogously local vector fields $Y\in
\mathfrak{X}(M)$ and $W\in \Gamma(\mathcal{V})$ such that
$(Y+W,\beta)\in \Gamma(D)$ and $\widehat{Y}=Y(\operatorname{mod} \mathcal{V})$. This yields 
\begin{align*}
\langle(\widehat{X},\alpha),(\widehat{Y},\beta)\rangle_{\mathcal{K}^\perp/\mathcal{K}} &\overset{\eqref{paring_k_k_perp}}=\langle(X+V,\alpha),(Y+W,\beta)\rangle
=0,
\end{align*}
since $(X+V,\alpha), (Y+W,\beta)\in \Gamma(D)$.

To prove the inclusion, $\tilde{D}(m)^\perp\subseteq\tilde{D}(m)$ let
$(\widehat{X}(m),\alpha(m)) \in  \tilde{D}(m)^\perp $ be such that 
$(\widehat{X},\alpha) \in \Gamma(\mathcal{K}^\perp/\mathcal{K})$ and for
all $(\widehat{Y},\beta) \in \Gamma(\tilde{D})$ we have
$\langle(\widehat{X},\alpha),(\widehat{Y},\beta)
\rangle_{\mathcal{K}^\perp/\mathcal{K}}
=0$. Choose $X\in
\mathfrak{X}(M)$ such that $\widehat{X}=X(\operatorname{mod} \mathcal{V})$. For all $(Y,\beta)\in \Gamma(D\cap \mathcal{K}^\perp)$,  $(\widehat{Y},\beta)$ lies in $\Gamma(\tilde{D})$ and we get 
\begin{align*}
0=\langle(\widehat{X},\alpha),(\widehat{Y},\beta)\rangle_{\mathcal{K}^\perp/\mathcal{K}}=\langle(X,\alpha),(Y,\beta)\rangle=\alpha(Y)+\beta(X).
\end{align*}
This yields
$(X,\alpha)\in \Gamma((D\cap\mathcal{K}^\perp)^\perp)$. We have
$(D_q\cap\mathcal{K}_q^\perp)^\perp=
D_q^\perp+(\mathcal{K}_q^\perp)^\perp=D_q+\mathcal{K}_q$ for every $q $ in the domain
of definition of $(X, \alpha) $.
Thus, since $D$ and $\mathcal{K}$ are smooth vector bundles, there exists $X'\in \mathfrak{X}(M)$ and $W\in \Gamma(\mathcal{V})$ such that $(X',\alpha)\in \Gamma(D)$ and $X=X'+W$. Now recall that the $1$-form 
$\alpha$ is in fact in $\Gamma(\mathcal{V}^\circ)$ since
$(\widehat X,\alpha)$ was an element of $\Gamma(\mathcal{K}^\perp/\mathcal{K})$.  The pair $(X',\alpha)$ is consequently in 
$\Gamma(D\cap\mathcal{K}^\perp)$ and, since $\widehat{X}=(X'+W)(\operatorname{mod}\mathcal{V})=X'(\operatorname{mod}\mathcal{V})$, our $(\widehat{X},\alpha)$ is  a local section of $\tilde{D}$, as required.
\end{proof}

This proposition immediately implies that $\dim\tilde{D}(m)$ is constant on 
$M$ and equal to
\[
\frac{\dim \mathcal{K}^\perp(m) - \dim\mathcal{K}(m) }{2}=\frac{\dim M+(\dim M-\dim
G)-\dim G}{2}= \dim M-\dim G.
\] 
Thus $\tilde{D}$ is a smooth $G$-invariant subbundle of
$\mathcal{K}^\perp/\mathcal{K}$. Its image by the isomorphism \eqref{isocou}
gives a  subbundle $D_{\rm red}=\tilde{D}/G$ of
\[
\left.\frac{\mathcal{K}^\perp}{\mathcal{K}} \right/G\; \simeq \;T \bar{M}
\oplus T ^* \bar{M},
\]
whose rank is $(\dim M-\dim G)$, which is isotropic relative to the symmetric pairing on $T\bar M\oplus T^*\bar{M}$. Hence $D_{\rm red}$ is a
Dirac structure called the \textit{reduction of $ D $ by $ G $}.
This discussion and Proposition \ref{tilde_D} yield
that the sections of $D_{\rm red}$ are in one-to-one correspondence with the
$G$-equivariant sections of the quotient \eqref{Dred} via the isomorphisms
$\mathfrak{X}(\bar{M})\simeq \Gamma(TM/\mathcal{V})^G$ and
$\Omega^1(\bar{M})\simeq\Gamma(\mathcal{V}^\circ)^G$ given at the beginning of this subsection.

It is customary to denote the ``quotient'' Dirac structure on $M/G$ by
\[
D_{\rm red}=\left.\frac{(D\cap \mathcal{K}^\perp)+
\mathcal{K}}{\mathcal{K}}\right/G.
\]
It is easy to check that if the Dirac structure
$D$ is integrable, then the reduced Dirac structure is also integrable (see also \cite{StXu08}).

\medskip

There is an other  formulation in terms of smooth sections for the reduction of Dirac manifolds and 
implicit Hamiltonian systems, due to \cite{Blankenstein00} and \cite{BlvdS01} 
(see \cite{BlRa04} for the singular case). This  extension of the Poisson reduction
 was historically the first method to reduce Dirac structures.

In \cite{BlvdS01}, there is the additional assumption that
$\mathcal{V}+\mathsf{G}_0$ is constant dimensional on $M$. 
 It turns out that no hypothesis on 
$\mathsf{G_0}$ is needed; the proof of the theorem
in \cite{BlvdS01} works with the same assumptions as in \cite{BuCaGu07} if one just uses the results 
of \cite{JoRa11} applied to $D\cap \mathcal K^\perp$, assuming only that this intersection
is smooth (see \cite{JoRaZa11}).

The formula for the reduced Dirac structure given in \cite{BlvdS01} is
\begin{align}\label{bard}
\Gamma\left(D_{\rm red}\right)=\left\{\left(\bar X,\bar \alpha\right)\in \Gamma(T\bar{M}\oplus T^*\bar{M})
\mid \text{there is }
X\in \mathfrak{X}(M) \text{ such that } X\sim_\pi\bar{X} \text{ and
}(X,\pi^*\bar{\alpha})\in \Gamma(D)\right\}.
\end{align}

The description of $D_{\rm red}$ shows that the smooth distribution ${\mathsf{G}}_0/\mathcal{V} $ projects to $\mathsf{G}^{\rm red}_0$ and that the smooth codistribution $\pi_2 \left( D\cap(\mathcal{V}\oplus \mathcal{V}^\circ) \right)$ projects to $\mathsf{P}^{\rm red}_0$, where $\pi_2$ is the projection  $\pi_2:TM\oplus T^*M\to T^*M$. There is no analogous description as quotients of the distribution $\mathsf{G}^{\rm red}_1$, and $ \mathsf{P}^{\rm red}_1$;
they need to be computed from the definition on a case by case basis. 
\medskip

Depending on the example, one needs to choose which method of Dirac reduction is easier to implement. In the next section, we will present cases where we have global bases of sections for the Dirac structure and in that situation the first method is more convenient. 
\medskip

There is a  third method of reduction that is
due to \cite{MaYo06, MaYo07, MaYo09}. It is undergoing a major extension to encompass both the Lagrangian and Hamiltonian version of classical reduction (see \cite{CeMaRaYo08}). Since this work is still in progress we shall not comment on it here.


\section{Reduction of nonholonomic systems}\label{nonhol}

\subsection{Summary of the nonholonomic reduction method}
\label{sec:nonhol_red_summary}
\cite{BaSn93} propose a reduction method for constrained
Hamiltonian systems. They start with the configuration space $Q$, a
hyperregular Lagrangian $L:TQ \rightarrow  \mathbb{R}$  taken as the kinetic energy of a Riemannian metric minus a potential,
 and a
\emph{constraint distribution} $\mathcal{D}$  on $Q$ equal to the kernel of smooth  $1$-forms 
$\phi^1,\ldots,\phi^k \in  \Omega ^1 (Q)$ satisfying pointwise $\phi^1\wedge\ldots\wedge\phi^k\neq 0$, that is,
\[
\mathcal{D} : = \{v \in  TQ \mid \phi ^j (v) = 0 ,\,  j=1, \ldots ,k \}.
\]
The independence of the forms  ensures that $ \mathcal{D} $ is a smooth vector subbundle of $ TQ $.

Denote by  $\langle\cdot\,,\cdot\rangle: T ^\ast Q \times  T Q \rightarrow  \mathbb{R} $ the duality pairing between $1$-forms and tangent vectors. Let
 $ \mathbb{F} L : TQ \rightarrow  T ^\ast Q $,
\[
\left\langle \mathbb{F} L (v), w \right\rangle : = \left.\frac{d}{dt}\right|_{t=0} L(v+tw), \quad v, w \in  T_qQ, 
\]
be the Legendre transformation associated to $ L $ which is  a
diffeomorphism since the Lagrangian is hyperregular. If $ A(v) := \left\langle \mathbb{F} L (v), v \right\rangle $
denotes the action of $ L $, let $ H (p): = A ((\mathbb{F} L) ^{-1} (p)) -
L((\mathbb{F} L) ^{-1} (p)) $, $ p \in  T ^\ast  Q$, be the associated
Hamiltonian. The Hamiltonian vector field $ X $  determined by $ H $ and the
constraint forms $\phi^1,\ldots,\phi^k \in  \Omega ^1 (Q)$ is defined
classically by
\begin{align*}
\dot{q} =\frac{\partial H}{\partial p}, \qquad 
\dot{p}=\frac{\partial H}{\partial q} +\lambda_j\phi^j
\end{align*}
or 
\begin{equation}
\label{nonholonomic_hamiltonian_vector_field}
{\mathbf{i}}_{X}\omega_{\rm can}=\mathbf{d} H+\lambda_j{\pi_{T^*Q}}^*\phi^j
\end{equation}
where $\pi_{T^*Q}:T^*Q\to Q$ is the cotangent bundle projection and $ \lambda_1 , \ldots \lambda _k \in  C ^{\infty} (Q) $ are the Lagrange multipliers
associated to the constraint forms $\phi^1,\ldots,\phi^k$,
and the constraint equations
\begin{equation}\label{constraints}
\phi^j( \mathbb{F} L^{-1}\left(p)\right)=\phi^j(T\pi_{T^*Q} X)=0, \qquad  p\in T^*Q.
\end{equation}
The counterpart of the constraint distribution $\mathcal{D} $ in phase space is the
\emph{constraint manifold} 
\begin{equation}
\label{def_M}
M:=\mathbb{F}L(\mathcal{D}) 
= \left\{p \in T^*Q\mid \phi^j(\pi_{T^*Q}(p)) \left((\mathbb{F} L)^{-1}(p) \right)=0,\, j=1,\ldots,k\right\}
 \subset T ^\ast Q.
\end{equation}
Since we require that the solution is in the constraint submanifold $M $,
it follows that $X$ is tangent to $M$.

Set $\omega_M:=i^*\omega_{\rm can}$, where $i:M\hookrightarrow T^*Q$ is the inclusion and $\omega_{\rm can}$ is the canonical symplectic form on $T^*Q$. Define 
\begin{equation}
\label{def_F}
\mathcal{F}:=\{U\in TT^*Q\mid \pi_{T^*Q}^*\phi^j(U)=0, \, j=1,\ldots,k\} 
\end{equation}
and note that $\pi_{T^*Q}^*\phi^1 \wedge \ldots \wedge \pi_{T^*Q}^*\phi^k \neq 0 $ on $T ^\ast Q $. Therefore $\mathcal{F} \rightarrow T ^\ast Q $ is a vector subbundle of $T T ^\ast Q $. The \emph{nonholonomic horizontal distribution} is defined by
\begin{equation}
\label{def_cal_H} 
\mathcal{H}:= \mathcal{F}\cap TM \rightarrow M.
\end{equation}
\cite{BaSn93} prove that the restriction $\omega_\mathcal{H}$ of $\omega_M$ to
$\mathcal{H}\times\mathcal{H}$ is nondegenerate. (Their proof uses the fact
that the Lagrangian is the kinetic energy of a metric minus a potential.) They also show that $\mathcal{H}$ is a vector subbundle of $TM $. 
With the condition \eqref{constraints} on $X$, we get for  $j=1,\ldots,k$,
\[\pi_{T^*Q}^*\phi^j(X)=\phi^j(T\pi_{T^*Q}X)=0\]
and thus the vector field $X$ is a section of $\mathcal{H}$.
Hence it is easy to see  that the pull back to $M $ of
\eqref{nonholonomic_hamiltonian_vector_field} subject to the constraints \eqref{constraints} is
equivalent to $X\in\Gamma(\mathcal{H})$ and ${\mathbf{i}}_{X}\omega_\mathcal{H}= \mathbf{d} H \arrowvert_{\mathcal{H}}$.

\medskip

Let $G$ be a Lie group acting symplectically on $T^*Q$ (not necessarily the lift of an action on $Q $), 
that leaves $M$ invariant
and preserves the Hamiltonian $ H $. Furthermore, assume that   the quotient $\bar{M}=M/G$ is
 a smooth manifold with projection map $\pi:M\to\bar{M}$, a submersion. Since $G$ is a symmetry group of the
nonholonomic system,
all intrinsically defined vector fields and distributions push down to $\bar{M}$.

In particular,
the vector field $X$ on $M$ pushes down to a vector field $\bar{X}$ with
$X\sim_\pi \bar X$ and the distribution
$\mathcal{H}$ pushes down to a distribution $\mathcal{H}_{\rm red}$ on $\bar{M}$. However, $\omega_\mathcal{H}$ need not push down to a
$2$-form defined on $\mathcal{H}_{\rm red}$ on $\bar{M}$, despite the fact that $\omega_\mathcal{H}$ is $G$-invariant. This is because there
may be infinitesimal symmetries $\xi_M$  which are horizontal (that is, take values in $ \mathcal{H}$), but
such that ${\mathbf{i}}_{\xi_M}\omega_\mathcal{H}\neq 0$.
Let $\mathcal{V}$ be the distribution on $M$ tangent to the orbits of $G$, that is, its fibers are $ \mathcal{V} (m) : = \{ \xi _M (m) \mid \xi \in  \mathfrak{g} \} $ for all  $m \in  M \subseteq T ^\ast Q $. Define the \emph{horizontal annihilator} $\mathcal{U}$ of $\mathcal{V}$ by
\begin{equation}
\label{def_cal_U}
\mathcal{U}=(\mathcal{V}\cap \mathcal{H})^{\omega_M}\cap \mathcal{H} \subseteq TM \subseteq T T ^\ast Q,
\end{equation}
where the superscript $ \omega _M $ on a distribution denotes its fiberwise $ \omega _M $-orthogonal complement in $ TM $.
Clearly, $\mathcal{U}$  and $\mathcal{V}$ are both $G$-invariant, project down to $\bar{M}$, and the image of $\mathcal{V}$ is $\{0\}$. Define $\bar{\mathcal{H}} : = T \pi \left( \mathcal{U}\right) \subseteq T\bar{M}$ to be the projection of $\mathcal{U}$ to $ \bar{M}$.
\cite{BaSn93} show that $X$ takes values in $\mathcal{U}$ and that the restriction
$\omega_\mathcal{U}$ of $\omega_M$ to $\mathcal{U} \times \mathcal{U}$ pushes down to a
nondegenerate form $\omega_{\bar{\mathcal{H}}}$ on $\bar{\mathcal{H}}$, i.e.,
$\pi^*\omega_{\bar{\mathcal{H}}}=\omega_\mathcal{U}$.
In addition, the function $\bar{H}\in C^\infty(\bar{M})$ defined by $\pi^*\bar{H}=H\arrowvert_{M}$ and the induced vector field $\bar{X}$ on $ \bar{M}$ are related by
\begin{equation}
\label{def_nonhol_red_ham_equ}
{\mathbf{i}}_{\bar X}\omega_{\bar{\mathcal{H}}}=\mathbf{d} \bar H\arrowvert_{\bar{\mathcal{H}}}
\end{equation}
which can be interpreted as the definition of the reduced nonholonomic Hamiltonian vector field $ \bar{X}$.
\medskip

\begin{remark} Note that we have no information about the dimensions of the fibers of $\mathcal{U}$. 
In general, $\mathcal{U}$ is \textit{not} a  vector subbundle of $TM $.
In the following, we will often assume that
$\mathcal V\cap\mathcal H$ has constant rank on the manifold $M$. In this case, $\mathcal U
=(\mathcal{V}\cap \mathcal{H})^{\omega_M}\cap \mathcal{H}
=(\mathcal{V}\cap \mathcal{H})^{\omega_{\mathcal H}}$ also
has constant rank on $M$, since $\omega_{\mathcal H}$ is nondegenerate.
\end{remark}

\subsection{Link with Dirac reduction}\label{link-nonhol-dirac}
Let $M$, $\omega_M$, $\pi_{T^*Q}$, $\mathcal{H}$, $ \bar{M} $, and $ \pi : M \rightarrow  \bar{M}$  be as in the preceding subsection. An easy verification shows that 
\begin{equation}
\label{h_and_delta_q}
 \mathcal{H} = (T(\pi_{T^*Q}\arrowvert_{M}))^{-1}(\mathcal{D})
\subseteq TM \subseteq T T ^\ast Q,
\end{equation}  
where 
\[
\mathcal{D}: = \{v \in  TQ \mid \left\langle \phi ^j, v \right\rangle = 0,  \, j=1, \ldots, k\} \subseteq TQ
\]
 is the constraint distribution on $Q$.

We introduce the Dirac structure $ D $ on $M$ as in \cite{MaYo07}: 
\begin{equation}
\label{ym_dirac}
D(m)=\left\{(X(m),\alpha(m))\in TM\oplus T^*M\mid X\in \Gamma(\mathcal{H}),\, \alpha-{\mathbf{i}}_{X}\omega_M \in \Gamma(\mathcal{H}^\circ),\, \alpha\in\Omega^1(M)\right\}
\end{equation}
for all $m\in M$.

The Lie group $G$ acts on $M$ and leaves $\mathcal{H}$, $\omega_M$, and thus the Dirac structure $ D $ invariant. Define $\mathcal{K}: = \mathcal{V} \oplus \{0\} \subset TM \oplus T^* M $ and its orthogonal complement $\mathcal{K} ^\perp = TM \oplus \mathcal{V}^ \circ $ as in \S\ref{red}.
Assume, as in  \S\ref{red}, that $D\cap
\mathcal{K}^\perp$ is a vector subbundle of $ TM \oplus T ^\ast M $ and 
consider the reduced Dirac manifold $(\bar M,D_{\rm red})$. 
The next proposition shows that, if $\bar{\mathcal{H}}$ is constant
dimensional, the reduced Dirac structure is  given
by the formula
\[D_{\rm red} =\left\{(X,\alpha)\in\Gamma(T\bar M\oplus T^*\bar M)\mid
    X\in\Gamma(\bar{\mathcal{H}}), \,\alpha\in\Omega^1(\bar M),\,
    \alpha-{\mathbf{i}}_{X}\omega_{\bar{\mathcal{H}}}\in\Gamma\left(\bar{\mathcal{H}}^\circ\right)\right\}\]
where $\bar{\mathcal{H}}$ and $\omega_{\bar{\mathcal{H}}}$ are defined as in
the preceding subsection.

\begin{proposition}\label{prop_nonhol_dirac}
Let $D$ be as above and assume that $\mathcal H\cap \mathcal V$ has constant rank on $M$.
\begin{enumerate}
\item[{\rm (i)}] The associated generalized distribution $\mathsf{G}_0$ is trivial and the
  codistribution $\mathsf{P}_1$ is given by $\mathsf{P}_1=T^*M$.
\item[{\rm (ii)}] Let $\mathcal{U}=\mathcal{H}\cap(\mathcal{V}\cap
  \mathcal{H})^{\omega_M}$ {\rm (}see \eqref{def_cal_U}{\rm )}. Then
\begin{equation} \label{eqUDK} 
X\in \Gamma(\mathcal{U})
\Longleftrightarrow\text{ there exists }
\alpha\in \Gamma(\mathcal{V}^\circ) \text{ such that } (X,\alpha)\in
\Gamma(D\cap\mathcal{K}^\perp).
\end{equation}
With the additional assumption that $\mathcal{V}+\mathcal{H}=TM$, the
section $\alpha$ in \eqref {eqUDK} is unique.
\item[{\rm (iii)}] The reduced distributions $\mathsf{G}_1^{\rm red}$ and
  $\mathsf{G}_0^{\rm red}$ are given by
\begin{equation*}
 \mathsf{G}_1^{\rm red}=\bar{\mathcal{H}}\quad \text{ and }  \quad\mathsf{G}_0^{\rm red}=\{0\}.
\end{equation*}
\item[{\rm (iv)}] For each
  $\alpha\in\Gamma(\mathcal{V}^\circ)$ there exists exactly one section
  $X\in\Gamma(\mathcal{U})$ such that $(X,\alpha)\in\Gamma(D)$. Hence, we have
  $\pi_2(D\cap\mathcal{K}^\perp)=\mathcal{V}^\circ$ and the reduced
  codistribution $\mathsf{P}_1^{\rm red}$ is equal to $T^*(M/G)$.
\item[{\rm (v)}] Assume that $\mathsf{G}_1^{\rm red}=\bar{\mathcal{H}}$ is constant
  dimensional.  The $2$-form defined on $\mathsf{G}_1^{\rm red}=\bar{\mathcal{H}}$ by
  the Dirac structure $D_{\rm red}$ {\rm (}see \eqref{gone}{\rm )} is nondegenerate and is equal to
  $\omega_{\bar{\mathcal{H}}}$.
\end{enumerate}
\end{proposition}

\begin{proof}
(i)  If $X$ is a section of
$\mathsf{G}_0$, we have  ${\mathbf{i}}_{X}\omega_M\in \Gamma(\mathcal{H}^\circ)$  and $X\in \Gamma(\mathcal{H})$. Hence, since $\omega_\mathcal{H}$ is nondegenerate, the vector field $X$ has to be the zero section. Thus 
$\mathsf{G}_0 = \{0\}$.
Since the $2$-form $\omega_{\mathcal{H}}$ is nondegenerate an arbitrary 
  $\alpha\in\Omega^1(M)$ determines a unique section $X$ of $\mathcal{H}$ by the equation
  ${\mathbf{i}}_{X}\omega_{\mathcal{H}}=\alpha\arrowvert_{\mathcal{H}}$. Therefore,  $\mathsf{P}_1=T^*M$.

(ii) If $(X,\alpha)$ is a local section of  $D\cap \mathcal{K}^\perp$, then we
have $X \in \Gamma(\mathcal{H})$, $\alpha\in \Gamma(\mathcal{V}^\circ)$,  and $\alpha={\mathbf{i}}_{X}\omega_M$ on $\mathcal{H}$. Hence,
$\left( {\mathbf{i}}_{X}\omega_M\right) |_{\mathcal{H}\cap \mathcal{V}}=0$  and thus we have 
\[
X \in \Gamma(\mathcal{H}\cap(\mathcal{V}\cap \mathcal{H})^{\omega_M})=\Gamma(\mathcal{U}).
\] Conversely, if $X \in \Gamma(\mathcal{U})$, we have ${\mathbf{i}}_{X}\omega_M=0$ on
$\mathcal{V}\cap \mathcal{H}$. Since, by hypothesis, 
$\mathcal V\cap\mathcal H$ has constant rank, it is a vector subbundle so we can find a section
$\alpha\in\Gamma(\mathcal{V}^\circ)$ such that the restriction of $\alpha$ and
${\mathbf{i}}_{X}\omega_M$ to $\mathcal{H}$ are equal.

If, in addition, we make the usual assumption $\mathcal{V}+\mathcal{H}=TM$, we
have for each $X\in \Gamma(\mathcal{U})$ exactly one $\alpha\in \Omega^1(M)$
such that  $\alpha|_ {\mathcal{H} }={\mathbf{i}}_{X}\omega_M$  and
$\alpha|_{\mathcal{V}} = 0$.  

(iii) By construction, the constraint distribution $\mathsf{G}_1^{\rm red}$ associated to the Dirac structure $D_{\rm red}$ on $\bar{M}$ is given by
\[
\left.\frac{\mathcal{U}+\mathcal{V}}{\mathcal{V}}\right/G.
\]
This can obviously be identified with  
\[
\bar{\mathcal{H}}=T \pi (\mathcal{U}).\]

If we have $\bar X \in \Gamma(\mathsf{G}_0^{\rm red})$, then 
$(\bar X, 0)\in \Gamma(D_{\rm red})$ and there exists 
$X\in\mathfrak{X}(M)$ with $X\sim_\pi\bar X$ and $(X,0)\in\Gamma(D)$. Hence we have $X\in \Gamma(\mathsf{G}_0)$
and  since $\mathsf{G}_0=\{0\}$, we get $X = 0$. This shows that $\mathsf{G}_0^{\rm red}=\{0\}$.

(iv)  This follows directly from (i) and (ii).

(v) Let $\omega_{D_{\rm red}}$ be the 2-form defined on 
$\mathsf{G}_1^{\rm red}=\bar{\mathcal{H}}$ by the Dirac structure 
$D_{\rm red}$ (see \eqref{gone}).  If $X \in  \Gamma (\bar{\mathcal{H}})$ is
  such that $\omega_{D_{\rm red}}(X,Y)=0$ for all
  $Y\in\Gamma(\bar{\mathcal{H}})$, then $(X,0)$ is a section of $D_{\rm red}$
  and hence we find by (iii) that $X=0$. Thus  $\omega_{D_{\rm red}}$ is nondegenerate on $\bar{\mathcal{H}}$. Let $\bar X$ and $\bar Y$ be sections of $\bar{\mathcal{H}}$. We show now
  that $\omega_{D_{\rm
      red}}(\bar{X},\bar{Y})=\omega_{\bar{\mathcal{H}}}(\bar{X},\bar{Y})$.
Indeed, by definition, we have $\omega_{D_{\rm red}}(\bar X,\bar Y)=\bar\alpha(\bar Y)$, where
$\bar\alpha, \bar\beta\in\Omega^1(M/G)$ are such that $(\bar
X,\bar \alpha), (\bar Y,\bar\beta)\in\Gamma(D_{\rm
  red})$. Choose $X, Y \in \Gamma(\mathcal{U})$ with $X\sim_{\pi}\bar
X$, $Y\sim_{\pi}\bar Y$ and $(X,\pi^*\bar\alpha)$,
$(Y,\pi^*\bar\beta)\in\Gamma(D\cap\mathcal{K}^\perp)$. Then we have 
\begin{align*}
\omega_{D_{\rm red}}(\bar X,\bar Y)=\bar\alpha(\bar Y)=(\pi^*\bar\alpha)(Y)=\omega_{\mathcal{U}}(X,Y)
=\omega_{\bar{\mathcal{H}}}(\bar X,\bar Y),
\end{align*}
where the last equality follows simply from the definition of $\omega_{\bar{\mathcal{H}}}$.
\end{proof} 

We shall use part (ii) of this proposition to simplify certain computations in
the examples that follow. 

\begin{remark}\label{rem_dim_HplusV}
Note that if $\mathcal{H}\cap\mathcal{V}$ has constant rank on $M$, we have
automatically that $D\cap\mathcal{K}^\perp$ has constant dimensional fibers on $M$.

This is proved in the following way.
Since $\mathcal{H}$, $\mathcal{V}$, $\mathcal{H}\cap\mathcal{V}$ are vector
subbundles of $TM$, $\mathcal{H}+\mathcal{V}$ is also a subbundle of
$TM$. By the nondegeneracy of $\omega_{\mathcal{H}}$, we get that
$\mathcal{U}=(\mathcal{H}\cap\mathcal{V})^{\omega_M}\cap
\mathcal{H}=(\mathcal{H}\cap\mathcal{V})^{\omega_{\mathcal{H}}}$ has also
constant dimensional fibers on $M$ and is in particular a vector subbundle of
$\mathcal{H}$. Let $u$ be the dimension of the fibers of $\mathcal{U}$, $r$
the dimension of the fibers of $\mathcal{H}$. Then, if $n=\dim M$, $n-r$ is the rank of the
codistribution $\mathcal{H}^\circ$. Let finally $l$ be the rank of the
codistribution
$\mathcal{H}^\circ\cap\mathcal{V}^\circ=(\mathcal{V}+\mathcal{H})^\circ\subseteq\mathcal{H}^\circ$. Choose local
basis vector fields $H_1,\dots,H_r$ for $\mathcal{H}$ such that
$H_1,\dots,H_u$ are basis vector fields for $\mathcal{U}$. In the same way,
choose basis $1$-forms $\beta_1,\dots,\beta_{n-r}$ for $\mathcal{H}^\circ$  such that
$\beta_1,\dots,\beta_l$ are basis $1$-forms  for
$\mathcal{V}^\circ\cap\mathcal{H}^\circ$. Then a local basis of sections of $D$ is
\[\left\{(H_1,{\mathbf{i}}_{H_1}\omega_M), \dots,
    (H_r,{\mathbf{i}}_{H_r}\omega_M),(0,\beta_1), \dots, (0,\beta_{n-r})\right\}.\] 
The considerations above show that $D\cap\mathcal{K}^\perp$ is then spanned by
the sections
\[\left\{\left(H_1,{\mathbf{i}}_{H_1}\omega_M+\sum_{i=l+1}^{n-r}a_1^i\beta_i\right), \dots,
    \left(H_u,{\mathbf{i}}_{H_u}\omega_M+\sum_{i=l+1}^{n-r}a_u^i\beta_i\right),(0,\beta_1), \dots, (0,\beta_{l})\right\},\]
where $a_i^j$ are smooth functions chosen such that
${\mathbf{i}}_{H_j}\omega_M+\sum_{i=l+1}^{n-r}a_j^i\beta_i\in\Gamma(\mathcal{V}^\circ)$
for $j=1,\dots,u$. Since these sections are linearly independent, they are 
smooth local basis sections for  $D\cap\mathcal{K}^\perp$.
\end{remark}

\subsection{Examples}
\subsubsection{The constrained particle in space}\label{constrpart}
\cite{BaSn93} study the motion of the constrained particle in
space. The configuration space of this problem is $Q:=\mathbb{R}^3$ whose coordinates are denoted by $\mathbf{q}: = (x,y,z)$. They take the following concrete constraints on the velocities: 
\[
\mathcal{D}:=\ker(\mathbf{d} z-y\mathbf{d} x) = \left\{ \left. v_x \partial_x
+ v_y \partial_y + v_z \partial_z \, \right| \, v _z - y v _x = 0 \right\} \subset TQ.
\] 
The Lagrangian  and taken to be the kinetic energy of the Euclidean metric, that is,  
$L(\mathbf{q}, \mathbf{v}) : = \frac{1}{2} \| \mathbf{v}\|^2 $, and it is hyperregular. Hence the constraint manifold
 \eqref{def_M} is five dimensional and given by 
\[
M:=\{(x,y,z,p_x,p_y,p_z)\mid p_z=y p_x\}\subseteq T^*Q,
\]
where $(x,y,z,p_x,p_y, p_z)$ are the coordinates of $ T ^\ast Q $. The global coordinates on $ M $ are thus $(x,y,z,p_x,p_y) $. The pull back $ \omega_M $ of the canonical $2$-form $ \omega $ on $ T ^\ast Q $ to  $ M $ has hence the expression
\[
\omega_M= \mathbf{d}  x\wedge\mathbf{d} p_x+\mathbf{d} y\wedge\mathbf{d} p_y+\mathbf{d} z\wedge(p_x\mathbf{d} y+ y\mathbf{d} p_x).
\]
The Dirac structure $D$ modeling this problem is given by \eqref{ym_dirac}. Formula \eqref{h_and_delta_q} gives the vector subbundle
\[
\mathcal{H}:=(T(\pi_{T ^\ast Q}\arrowvert_{M}))^{-1}(\mathcal{D})=\operatorname{span}\{
\partial_x+y\partial_z, \partial_y, \partial_{p_x}, \partial_{p_y}\} \subset TM,
\]
and consequently
\[\mathcal{H}^\circ=\operatorname{span}\{\mathbf{d} z-y\mathbf{d} x\}.\]
A computation yields
\begin{align*}
{\mathbf{i}}_{\partial_x+y\partial_z}\omega_M&=(1+y^2)\mathbf{d} p_x+yp_x\mathbf{d} y, \qquad 
{\mathbf{i}}_{\partial_y}\omega_M=\mathbf{d} p_{y}-p_{x}\mathbf{d} z,\\
{\mathbf{i}}_{\partial_{p_y}}\omega_M&=-\mathbf{d} y, \qquad 
{\mathbf{i}}_{\partial_{p_x}}\omega_M=-y\mathbf{d} z-\mathbf{d} x.
\end{align*}
Hence
\begin{align}
\label{ex_particle_basis_D}
\big\{\left(\partial_x+y\partial_z,(1+y^2)\mathbf{d} p_x+yp_x\mathbf{d}
  y\right); \left(\partial_y, \mathbf{d} p_{y}-p_{x}\mathbf{d} z\right);
  \left(\partial_{p_y}, -\mathbf{d} y\right); \left(\partial_{p_x}, -y\mathbf{d}
  z-\mathbf{d} x\right); \left(0, \mathbf{d} z-y\mathbf{d} x\right)\big\} 
\end{align}
is a smooth global basis for $D$.

We consider the action of the Lie group $G=\mathbb{R}^2$ on $M$ given by
\[\Phi:G\times M\to M, \quad \Phi((r,s),m)=(x+r,y,z+s,p_x,p_y),\]
where ${\bf m}:=(x,y,z,p_x,p_y)\in M$.
This $ \mathbb{R} ^2 $-action is the restriction to $M $ of the cotangent
lift of the action $\phi:G\times Q\to Q$,\,
$\phi((r,s),(x,y,z))=(x+r,y,z+s)$.  It obviously leaves the Hamiltonian
$H({\bf m})=\frac{1}{2}((1+y^2)p_x^2+p_y^2)$ on $ M $ invariant.
Note that if $(X,\alpha)\in \Gamma (D)$
we have
\[({\boldsymbol{\pounds}}_{\xi_M}X,{\boldsymbol{\pounds}}_{\xi_M}\alpha)\in \Gamma(D) \quad \text{for all} \quad \xi \in \mathfrak{g} = \mathbb{R} ^2.\]
Since the vertical bundle in this example is $\mathcal{V}=\operatorname{span}\{\partial_x,\partial_z\}$, we have 
\[
\mathcal{K} =\mathcal{V}\oplus\{0\} =\operatorname{span}\{(\partial_x,0), (\partial_z, 0)\} \subset TM \oplus T ^\ast M
\]
and thus 
\begin{align}
\label{ex_particle_basis_K_perp}
\mathcal{K}^\perp&=TM \oplus\mathcal{V}^ \circ  
= \operatorname{span}\{(\partial_x,0), (\partial_y,0), (\partial_z,0), (\partial_{p_x},0), (\partial_{p_y},0), ( 0 , \mathbf{d}y ), ( 0, \mathbf{d}p _x) , (0, \mathbf{d}p _y) \}
\end{align}
A direct computation using \eqref{ex_particle_basis_D} and \eqref{ex_particle_basis_K_perp} yields
\begin{eqnarray*}
 D\cap \mathcal{K}^\perp=\operatorname{span}\left\{\left(\partial_{p_y},-\mathbf{d}
y\right),\left(\partial_x+y\partial_z,(1+y^2)\mathbf{d} p_x+yp_x\mathbf{d} y\right)\right., 
\left.\left((1+y^2)\partial_y-yp_x\partial_{p_x},(1+y^2)\mathbf{d} p_y\right)\right\}
\end{eqnarray*} 
and
\begin{align*}
(D\cap \mathcal{K}^\perp) +\mathcal{K}=\operatorname{span}\left\{\begin{array}{c}
\left(\partial_{p_y},
-\mathbf{d}y\right),\quad \left(\partial_x,0\right),\quad \left(\partial_z,0\right),\quad
\left(0,(1+y^2)\mathbf{d} p_x+yp_x\mathbf{d} y\right), \\
\left((1+y^2)\partial_y-yp_x\partial_{p_x},(1+y^2)\mathbf{d} p_y\right)
\end{array}
\right\}
\end{align*}
since in this case $( D \cap \mathcal{K}^ \perp) \cap \mathcal{K} = \{0\}. $

Note that there is an easier way to compute the spanning sections of $D\cap \mathcal{K}^\perp$ by using \eqref{eqUDK}. First, one determines spanning sections of $\mathcal{U}$. Second, for each spanning section $X\in\Gamma(\mathcal{U})$ we find $\lambda \in C ^{\infty}(M)$ such that 
\[
{\mathbf{i}}_{X}\omega_M +\lambda
(\mathbf{d} z-y\mathbf{d} x) \in \Gamma(\mathcal{V}^\circ).
\]
Third, setting $\alpha : = {\mathbf{i}}_{X}\omega_M +\lambda
(\mathbf{d} z-y\mathbf{d} x) $ we have found a spanning section 
$(X, \alpha) \in \Gamma(D\cap \mathcal{K}^\perp)$.
In the following examples, we will proceed like this.

We get the reduced Dirac structure
\begin{eqnarray} \label{ex1_D_red}
 D_{\rm red}=\left.\frac{(D\cap \mathcal{K}^\perp)+\mathcal{K}}{\mathcal{K}}\right/G \nonumber  
=\operatorname{span}\left\{\begin{array}{l}\left(\partial_{p_y},-\mathbf{d}
y\right),\left(0,(1+y^2)\mathbf{d} p_x+yp_x\mathbf{d} y\right), \\
\left((1+y^2)\partial_y-yp_x\partial_{p_x},(1+y^2)\mathbf{d} p_y\right)\end{array}\right\}
\end{eqnarray}
on the three dimensional manifold $\bar{M}: = M/G$ with global coordinates  $(y,p_y,p_x)$.

Since $\partial_x+y\partial_z$ is a spanning section of
$\mathcal{H}\cap\mathcal{V}$, the distribution $ \mathcal{U} \subset  TM $ (see \eqref{def_cal_U}) is given by 
\begin{align*}
\mathcal{U}& = ( \mathcal{V}\cap \mathcal{H}) ^ { \omega_M} \cap   \mathcal{H} 
=\ker\{{\mathbf{i}}_{\partial_x+y\partial_z}\omega_M\}\cap\mathcal{H}
=\ker\{(1+y^2)\mathbf{d} p_x+yp_x\mathbf{d} y\}\cap\mathcal{H}\\
&= \operatorname{span} \left\{(1+y^2)\partial_y-yp_x\partial_{p_x},\partial_x+y\partial_z, \partial_{p_y}  \right\}.
\end{align*}
Thus
\[
\bar{\mathcal{H}} = T \pi ( \mathcal{U} ) =\operatorname{span}\{\partial_{p_y},(1+y^2)\partial_y-yp_x\partial_{p_x}\}
\] 
recovering the result in \cite{BaSn93}. Note that, as discussed in 
\S\ref{link-nonhol-dirac}, the distribution $\bar{ \mathcal{H} } \subset  T \bar{M}$ coincides with the projection on the first factor of the reduced Dirac structure \eqref{ex1_D_red}. As in \cite{BaSn93}, $\bar{\mathcal{H}}$ is an integrable subbundle of $ T\bar{M}$; in fact $[ \partial_{p_y}, (1+y^2)\partial_y-yp_x\partial_{p_x}]= 0$. The $2$-form 
$\omega_{\bar{\mathcal{H}}}$ is easily computed to equal 
\begin{align*} 
\omega_{\bar{\mathcal{H}}}(\partial_{p_y},(1+y^2)\partial_y-yp_x\partial_{p_x})=
-\mathbf{d} y((1+y^2)\partial_y-yp_x\partial_{p_x})=-(1+y^2).
\end{align*} 
As predicted by the general theory in \S\ref{sec:nonhol_red_summary}, 
$\omega_{\bar{\mathcal{H}}}$ is nondegenerate.  

It is easy to check that the reduced manifold $\bar{M}$ is Poisson relative to the $2$-tensor
\[
-\partial_y\wedge\partial_{p_y}+\frac{yp_x}{1+y^2}\partial_{p_x}\wedge\partial_{p_y},
\]
or with Poisson bracket determined by $ \{y, p_y\} = -1$, $\{y, p_x\} = 0$, $\{p_y, p_x \} = yp_x/(1+y^2) $,
and that $D_{\rm red}$ given by \eqref{ex1_D_red} is the graph of the vector bundle homomorphism $\flat:T^*\bar{M} \rightarrow T\bar{M}$ associated to the Poisson structure.

\subsubsection{The vertical rolling disk}
\label{ex_vertical_disk}

This example is standard in the theory of nonholonomic mechanical systems; it can be found for example in \cite{Bloch03}.
Consider a vertical disk of zero width rolling on the $xy$-plane
and free to rotate about its vertical axis. Let $x$  and $y$ denote the position of
contact of the disk in the $xy$-plane. The remaining variables are $\theta$ and $\varphi$,
denoting the orientation of a chosen material point $P$ with respect to the
vertical and the ``heading angle'' of the disk.
Thus, the unconstrained configuration space for the vertical rolling disk
is $Q :=\mathbb{R}^2\times\mathbb{S}^1\times\mathbb{S}^1$. The Lagrangian for the problem is taken to
be the kinetic energy
\[L(x,y,\theta, \varphi,\dot{x},\dot{y},\dot{\theta},\dot{\varphi})=\frac{1}{2}\mu(\dot{x}^2+\dot{y}^2)+\frac{1}{2}I\dot{\theta}^2+\frac{1}{2}J\dot{\varphi}^2,\]
where $\mu$ is the mass of the disk, and $I$, $J$ are its moments of inertia. Hence, the Hamiltonian of the system is 
\[
H(x,y,\theta, \varphi,p_x,p_y,p_\theta,p_ \varphi)=\frac{1}{2\mu}(p_x^2+p_y^2)+\frac{1}{2I}p_\theta^2+\frac{1}{2J}p_ \varphi ^2.
\]
The rolling constraints may be written as $\dot{x}=R\dot{\theta} \cos \varphi $ and $\dot{y}=R\dot{\theta}\sin \varphi $, where $R$ is the radius of the disk, that is,
\[
\mathcal{D} : = \{(x,y,\theta, \varphi,R\dot{\theta} \cos \varphi
,R\dot{\theta}\sin \varphi ,\dot{\theta},\dot{\varphi}) \mid x, y \in  \mathbb{R}, \theta, \varphi \in  \mathbb{S}^1 \} \subset  TQ.
\]
Note that the $1$-forms defining this distribution $ \mathcal{D} $ are 
$ \phi_1 : = \mathbf{d} x - R \cos \varphi \mathbf{d} \theta $  and $ \phi _2 : = \mathbf{d} y - R \sin \varphi \mathbf{d} \theta $.

The constraint manifold \eqref{def_M} 
\[
M: =\left\{(x,y,\theta, \varphi,p_x,p_y,p_\theta,p_ \varphi)\in T^*Q\, \left|\,  p_x=\frac{\mu R}{I} p_\theta \cos \varphi,
\,p_y=\frac{\mu R}{I} p_\theta\sin \varphi  \right. \right\}\subseteq T^*Q\]
is in this example a graph over the coordinates $(x,y,\theta, \varphi,p_\theta,p_ \varphi)$ and is hence six dimensional.
The induced $2$-form $\omega_M=i^*\omega_{\rm can}$ is given by the formula
\begin{align*}
\omega_M=& \mathbf{d} x\wedge\left( \frac{\mu R\cos\varphi}{I}\mathbf{d} p_\theta-\frac{\mu R\sin\varphi}{I}p_\theta\mathbf{d} \varphi \right) +
\mathbf{d} y\wedge \left( \frac{\mu R\sin\varphi}{I}\mathbf{d} p_\theta+\frac{\mu R\cos\varphi}{I}p_\theta\mathbf{d} \varphi \right) \\
&+\mathbf{d}\theta\wedge\mathbf{d} p_\theta+\mathbf{d}\varphi\wedge\mathbf{d} p_\varphi
\end{align*}
and the distribution
$\mathcal{H}=\ker\{\mathbf{d} x-R\cos\varphi \mathbf{d} \theta,\,\,\mathbf{d} y-R \sin\varphi \mathbf{d}
\theta\}\subseteq TM$
is in this case
\begin{equation}
\label{ex2_h_span}
\mathcal{H}=\operatorname{span}\{\partial_\varphi,\partial_\theta+R\cos\varphi \,\partial_x+R\sin\varphi\,\partial_y,\partial_{p_\theta},\partial_{p_\varphi}\} \subset  TM.
\end{equation} 
Therefore its annihilator is
\[
\mathcal{H} ^ \circ  = \operatorname{span} \{ \mathbf{d} x - R \cos \varphi \mathbf{d} \theta , \mathbf{d} y - R \sin \varphi \mathbf{d} \theta \} \subset  T ^\ast M. 
\]
The Dirac structure on $M$ describing the nonholonomic mechanical system is again given by \eqref{ym_dirac}. Since
\begin{align*}
{\mathbf{i}}_{\partial_\varphi}\omega_M&=\mathbf{d} p_\varphi+\frac{\mu R\sin\varphi}{I}p_\theta\mathbf{d} x-\frac{\mu R\cos\varphi}{I}p_\theta\mathbf{d} y,\\
{\mathbf{i}}_{\partial_{p_\theta}}\omega_M&=-\frac{\mu R\cos\varphi}{I}\mathbf{d} x-\frac{\mu R\sin\varphi}{I}\mathbf{d} y-\mathbf{d} \theta,\\
{\mathbf{i}}_{\partial_{p_\varphi}}\omega_M&=-\mathbf{d}\varphi,\\
\end{align*}
and
\begin{align*}
&{\mathbf{i}}_{\partial_\theta+R\cos\varphi \partial_x+R\sin\varphi\partial_y}\omega_M\\
=\,&\mathbf{d} p_\theta+R\cos\varphi\left(\frac{\mu R\cos\varphi}{I}\mathbf{d}
p_\theta-\frac{\mu R\sin\varphi}{I}\mathbf{d}\varphi\right) +R\sin\varphi\left(\frac{\mu R\sin\varphi}{I}\mathbf{d}
p_\theta+\frac{\mu R\cos\varphi}{I}\mathbf{d}\varphi\right)\\
=\,&\left(1+\frac{\mu R^2}{I} \right) \mathbf{d} p_\theta,
\end{align*} 
we get again smooth global spanning sections of $D$:
\begin{equation}
\label{ex2_global_sections_D}
\begin{aligned}
\left(\partial_\varphi, \mathbf{d} p_\varphi+\frac{\mu
    R\sin\varphi}{I}p_\theta\mathbf{d} x-\frac{\mu
    R\cos\varphi}{I}p_\theta\mathbf{d} y\right)\\
\left(\partial_{p_\theta}, -\frac{\mu R\cos\varphi}{I}\mathbf{d} x-\frac{\mu
    R\sin\varphi}{I}\mathbf{d} y-\mathbf{d} \theta\right)\\
\left(\partial_{p_\varphi}, -\mathbf{d}\varphi\right),\quad    
\left(0, \mathbf{d} x-R\cos\varphi \mathbf{d} \theta\right), \quad \left(0,\mathbf{d} y-R \sin\varphi \mathbf{d}\theta\right)\\
  \left(\partial_\theta+R\cos\varphi \partial_x+R\sin\varphi\partial_y, \left(1+\frac{\mu R^2}{I} \right) \mathbf{d} p_\theta\right).
\end{aligned}
\end{equation} 
In this case, several groups of symmetries are studied in the literature.

\begin{enumerate}
\item \textit{The case $G=\mathbb{R}^2$} (\cite{CaDiLeMa98}).

The Lie group $\mathbb{R}^2$ acts on $M$ by
\[
(r,s)\cdot(x,y,\theta,\varphi,p_\theta,p_\varphi)=(x+r,y+s,\theta,\varphi,p_\theta,p_\varphi)
\]
and clearly leaves the Hamiltonian $ H $ invariant.
The distribution $\mathcal{V}$ on $M$ is in this case 
$\mathcal{V}=\operatorname{span}\{\partial_x,\partial_y\}$, so that $ \mathcal{V} \cap \mathcal{H} =\{0\}$ by \eqref{ex2_h_span}. Therefore, in this example, $\mathcal{U} = \mathcal{H}$. We have
\[
\mathcal{K}= \mathcal{V} \oplus \{0\} = \operatorname{span}\{(\partial_x,0), (\partial_y, 0)\} \subset  TM \oplus T ^\ast M
\]
and 
\begin{align*}
\mathcal{K}^\perp&=TM\oplus\mathcal{V} ^ \circ  = \operatorname{span}\left\{\begin{array}{c}
(\partial_x, 0),\quad  (\partial_y, 0),\quad (\partial_\theta, 0),\quad (\partial_\varphi, 0),\quad\\
 (\partial_{p_\theta}, 0),\quad  (\partial_{p_\varphi}, 0), \quad 
 (0, \mathbf{d} p_ \varphi), \\
(0, \mathbf{d} \varphi ), \quad (0, \mathbf{d} p_{ \theta }), \quad (0, \mathbf{d} \theta ) 
\end{array}\right\} .
\end{align*} 
By \eqref{eqUDK} and the fact that $\mathcal{V}+\mathcal{H}=TM$, we know that for each spanning
section $X$ of $\mathcal{H}$, there exists exactly one $\alpha\in\Gamma(\mathcal{V}^\circ)$
such that the pair $(X,\alpha)$ is a section of $D\cap\mathcal{K}^\perp$.
Using \eqref{ex2_global_sections_D} and the equalities
\begin{align}
&{\mathbf{i}}_{\partial_\varphi}\omega_M-\frac{\mu R\sin\varphi}{I}p_\theta(\mathbf{d} x-R\cos\varphi\mathbf{d} \theta)+\frac{\mu R\cos\varphi}{I}(\mathbf{d} y-R\sin\varphi\mathbf{d}
\theta)=\mathbf{d} p_\varphi\label{searching_lambda_1}\\
&{\mathbf{i}}_{\partial_{p_\theta}}\omega_M+\frac{\mu R\cos\varphi}{I}p_\theta(\mathbf{d}
x-R\cos\varphi\mathbf{d} \theta) \nonumber \\
&\qquad \quad \; +\frac{\mu R\sin\varphi}{I}(\mathbf{d} y-R\sin\varphi\mathbf{d}\theta)= - \left(1+\frac{\mu R^2}{I}\right)\mathbf{d}\theta
\label{searching_lambda_2}
\end{align}
we find
\begin{eqnarray*}
 D\cap \mathcal{K}^\perp=\operatorname{span}\left\{\begin{array}{c}\left(\partial_{\varphi},\mathbf{d}
     p_\varphi\right),\left(\partial_{p_\theta},- \left(1+\frac{\mu R^2}{I} \right)\mathbf{d}\theta\right),\left(\partial_{p_\varphi},-\mathbf{d}\varphi\right),\\
\left(\partial_\theta+R\cos\varphi \partial_x+R\sin\varphi\partial_y,
\left(1+\frac{\mu R^2}{I} \right)\mathbf{d} p_\theta\right)\end{array}\right\}.
\end{eqnarray*} 
Hence
\begin{eqnarray*}
 (D\cap \mathcal{K}^\perp)+\mathcal{K}=\operatorname{span}\left\{\begin{array}{c}
\left(\partial_{\varphi},\mathbf{d}
     p_\varphi\right),\left(\partial_{p_\theta},-\left(1+\frac{\mu R^2}{I} \right)\mathbf{d}\theta\right),\left(\partial_{p_\varphi},-\mathbf{d}\varphi\right),\\
\left(\partial_\theta,\left( 1+\frac{\mu R^2}{I}\right)\mathbf{d} p_\theta\right),\left(\partial_x,0\right),\left(\partial_y,0\right)\end{array}\right\}
\end{eqnarray*}
and finally we get the reduced Dirac structure
\begin{eqnarray} \label{ex2_reduced_D_first_action}
 D_{\rm red}=\left.\frac{(D\cap \mathcal{K}^\perp)+\mathcal{K}}{\mathcal{K}}\right/G
=\operatorname{span}\left\{\begin{array}{l}\left(\partial_{\varphi},\mathbf{d}
     p_\varphi\right),\left(\partial_{p_\theta},-(1+\frac{\mu R^2}{I})\mathbf{d}\theta\right),\left(\partial_{p_\varphi},-\mathbf{d}\varphi\right), \\
\left(\partial_\theta,(1+\frac{\mu R^2}{I})\mathbf{d} p_\theta\right)\end{array}\right\}
\end{eqnarray}
on the four dimensional manifold $\bar{M} = M/G $ with coordinates $(\varphi,\theta,p_\varphi,p_\theta)$.
Thus, $D_{\rm red}$ is the graph of the symplectic form on $ \bar{M} $ given by
$\omega_{\rm red}=\mathbf{d}\varphi\wedge\mathbf{d} p_\varphi+(1+\frac{\mu R^2}{I})\mathbf{d}\theta\wedge\mathbf{d} p_\theta$.

As already mentioned, in this example, $ \mathcal{U} = \mathcal{H} $ and hence $ \bar{ \mathcal{H} } = T \pi ( \mathcal{H} ) = \operatorname{span}   \{ \partial _ \varphi , \partial _ { p_ \varphi }, \partial _ \theta , \partial _{ p_ \theta } \}$ by \eqref{ex2_h_span} which coincides with the projection on the first factor of the reduced Dirac structure \eqref{ex2_reduced_D_first_action}. In this case $ \bar{\mathcal{H} } = T \bar{M} $ and so $ \omega _{\bar{ \mathcal{H} }} = \omega _{\rm red} $ is of course nondegenerate.

\item \textit{The case $G=\operatorname{SE}(2)$} (\cite{Bloch03}). 

The Lie group $\operatorname{SE}(2): =\mathbb{S}^1 \,\circledS\, \mathbb{R} ^2$ is the semidirect product of the circle $\mathbb{S}^1$ identified with matrices of the form
\[
\begin{bmatrix}
\cos \alpha  & - \sin \alpha \\
\sin \alpha & \cos \alpha 
\end{bmatrix}
\]
and acting on $ \mathbb{R} ^2 $ by usual matrix multiplication. Denote
elements of $\operatorname{SE}(2) $ by $(\alpha,r,s) $ where
$\alpha\in\mathbb{S}^1$ and $ r,s \in
\mathbb{R}$. Define the action of the Lie group 
$\operatorname{SE}(2)$ on $ M $ by 
\[
(\alpha,r,s)\cdot(x,y,\theta,\varphi,p_\theta,p_\varphi)=(x\cos\alpha-y\sin\alpha+r,x\sin\alpha+y\cos\alpha+s,\theta,\varphi+\alpha,p_\theta,p_\varphi)
\]
and note that the Hamiltonian $ H $ is invariant by this action.
The distribution $\mathcal{V}$ on $M$ is in this case
$\mathcal{V}=\operatorname{span}\{\partial_x,\partial_y,\partial_\varphi\}$ and we get 
\[
\mathcal{K}=\mathcal{V} \oplus \{0\} = \operatorname{span}\{(\partial_x,0), (\partial_y,0), (\partial_\varphi, 0)\}.
\]
Thus 
\begin{align*} \label{ex2_k_perp_second_action}
\mathcal{K}^\perp&=TM \oplus \mathcal{V} ^ \circ = \operatorname{span}\left\{(\partial_x, 0), (\partial_y, 0), (\partial_\theta, 0), (\partial_\varphi, 0), (\partial_{p_\theta}, 0), (\partial_{p_\varphi}, 0), 
\right.   \\
& \qquad \qquad\qquad\qquad\qquad \left.
 (0, \mathbf{d} p_ \varphi), (0, \mathbf{d} p_{ \theta }), (0, \mathbf{d} \theta ) \right\}. 
\end{align*} 
We have $\mathcal{V}\cap\mathcal{H}=\operatorname{span}\{\partial_\varphi\}$ (see \eqref{ex2_h_span}) and hence
\[
( \mathcal{V} \cap \mathcal{H} ) ^{ \omega _M} = \ker\left( \mathbf{d}p_\varphi+\frac{\mu R\sin\varphi}{I}p_\theta\mathbf{d} x-\frac{\mu R\cos\varphi}{I}p_\theta\mathbf{d} y \right) 
\]
so that 
\begin{align*}
\mathcal{U}&=\mathcal{H}\cap\ker\left(\mathbf{d}p_\varphi+\frac{\mu
  R\sin\varphi}{I}p_\theta\mathbf{d} x-\frac{\mu R\cos\varphi}{I}p_\theta\mathbf{d} y \right) \\
&=\operatorname{span}\{\partial_\varphi,\partial_\theta+R\cos\varphi \,\partial_x+R\sin\varphi\,\partial_y,\partial_{p_\theta}\}.\end{align*}

Using \eqref{ex2_global_sections_D}, \eqref{searching_lambda_1}, and \eqref{searching_lambda_2}, we get
\begin{align*}
 D\cap \mathcal{K}^\perp=&\operatorname{span}\left\{\begin{array}{c}
\left(\partial_{\varphi},\mathbf{d}
     p_\varphi\right),\left(\partial_{p_\theta},-\left(1+\frac{\mu R^2}{I} \right) \mathbf{d}\theta\right)\\
\left(\partial_\theta+R\cos\varphi \partial_x+R\sin\varphi\partial_y,
\left(1+\frac{\mu R^2}{I}\right)\mathbf{d} p_\theta\right)\end{array}\right\}.
\end{align*} 
Thus, 
\begin{align*}
 D_{\rm red}=\left.\frac{(D\cap \mathcal{K}^\perp)+\mathcal{K}}{\mathcal{K}}\right/G
=\operatorname{span}\left\{\begin{array}{l}\left(0,\mathbf{d}
     p_\varphi\right),\left(\partial_{p_\theta},-(1+\frac{\mu R^2}{I})\mathbf{d}\theta\right), \\
\left(\partial_\theta,(1+\frac{\mu R^2}{I})\mathbf{d} p_\theta\right)\end{array}\right\}
\end{align*}
is the graph of the Poisson tensor
\[
\frac{I}{\mu R^2+I}\partial_{p_\theta}\wedge\partial_\theta 
\]
defined on the manifold $\bar{M} : = M/G$ with coordinates $(\theta, p_\theta, p_\varphi)$. 

In addition, 
\[
\bar{ \mathcal{H} } = T \pi ( \mathcal{U} ) = \operatorname{span}\{ \partial _ \theta , \partial _{p_ \theta} \}
\]
is an integrable subbundle of $T\bar{M}$ (since $[ \partial _ \theta , \partial _{p_ \theta }] = 0$). Note that the projection on the first factor of 
$D_{\rm red}$ equals $\bar{ \mathcal{H}}$. Finally, the $2$-form 
$\omega _{\bar{\mathcal{H}}}$ is easily computed to be
\[
\omega _{ \bar{\mathcal{H}}} \left(\partial _\theta, \partial_{p _\theta} \right)  = 1 + \frac{\mu R ^2}{I}
\]
and, as predicted by the general theory, it is nondegenerate on $\bar{ \mathcal{H}}$.

\item \textit{The case $G= \mathbb{S}^1  \times \mathbb{R}^2$} (\cite{Bloch03}).
The direct product Lie group $\mathbb{S}^1  \times \mathbb{R}^2$ acts on $M$ by
\[
(\alpha,r,s)\cdot(x,y,\theta,\varphi,p_\theta,p_\varphi)=(x+r,y+s,\theta+\alpha,\varphi,p_\theta,p_\varphi).
\]
The distribution $\mathcal{V}$ on $M$ is in this case $\mathcal{V}=
\operatorname{span}\{\partial_x,\partial_y,\partial_\theta\}$,
\[
\mathcal{K}= \mathcal{V} \oplus \{0\} = 
\operatorname{span}\{(\partial_x,0), (\partial_y,0), (\partial_\theta, 0)\},
\]
and thus 
\begin{align*}
\mathcal{K}^\perp&= TM \oplus \mathcal{V} ^ \circ   \\
& = \operatorname{span}\left\{(\partial_x, 0), (\partial_y, 0), (\partial_\theta, 0), (\partial_\varphi, 0), (\partial_{p_\theta}, 0), (\partial_{p_\varphi}, 0), 
 (0, \mathbf{d} p_ \varphi), (0, \mathbf{d} p_{ \theta }), (0, \mathbf{d} \varphi ) \right\}.
\end{align*} 
Using \eqref{ex2_h_span} we get  $\mathcal{V}\cap\mathcal{H}=\operatorname{span}\{\partial_\theta+R\cos\varphi \partial_x+R\sin\varphi\partial_y\}$ and hence
\[
(\mathcal{V}\cap\mathcal{H})^{\omega_M} = \ker \left\{ \left( 1 + \frac{\mu R^2}{I} \right) \mathbf{d} p_ \theta  \right\}.
\]
Therefore, again by \eqref{ex2_h_span} we conclude
\begin{align*}
\mathcal{U} & = \mathcal{H} \cap (\mathcal{V}\cap\mathcal{H})^{\omega_M}
= \mathcal{H} \cap \ker \left\{ \left( 1 + \frac{\mu R^2}{I} \right) 
\mathbf{d} p_ \theta  \right\} \\
& = \operatorname{span}\{\partial_\varphi,\partial_\theta+R\cos\varphi \,\partial_x+R\sin\varphi\,\partial_y,\partial_{p_\varphi}\}.
\end{align*} 
Using  \eqref{ex2_global_sections_D} and \eqref{searching_lambda_1}, we obtain
\begin{align*}
 D\cap \mathcal{K}^\perp =\operatorname{span} 
 \left\{\left(\partial_{\varphi},\mathbf{d}p_\varphi\right),
 \left(\partial_{p_\varphi}, -\mathbf{d}\varphi\right), \left(\partial_\theta+R\cos\varphi \partial_x+R\sin\varphi\partial_y,\left(1+\frac{\mu R^2}{I}\right)\mathbf{d} p_\theta\right)\right\}
\end{align*} 
and hence
\begin{align*}
 D_{\rm red}=\left.\frac{(D\cap \mathcal{K}^\perp)+\mathcal{K}}{\mathcal{K}}\right/G
=\operatorname{span}\left\{\left(\partial_\varphi,\mathbf{d}
     p_\varphi\right),\left(\partial_{p_\varphi},-\mathbf{d}\varphi\right),\left(0,\mathbf{d} p_\theta\right)\right\},
\end{align*}
which is the graph of the Poisson tensor 
\[
\partial_{p_\varphi}\wedge\partial_\varphi
\]
on the three dimensional reduced manifold $\bar{M} = M/G$ with coordinates $(\varphi, p_\varphi, p_\theta)$. We have
\[
\bar{ \mathcal{H} } = T \pi ( \mathcal{U} ) 
= \operatorname{span}\{ \partial_\varphi, \partial _{p_\varphi} \}
\]
which is an integrable subbundle of $T\bar{M}$ (since $[\partial_\varphi, \partial _{p_\varphi}]=0$). As before, the projection on the first factor of 
$D_{\rm red}$ equals $\bar{\mathcal{H}}$. The $2$-form 
$\omega_{\bar{\mathcal{H}}}$ has the expression
\[
\omega_{\bar{\mathcal{H}}}\left(\partial_\varphi, \partial _{p_\varphi} \right) = 1 
\]
and, as the general theory states, it is nondegenerate on $\bar{ \mathcal{H}}$.
\end{enumerate}

\subsubsection{The Chaplygin skate}\label{chaplygin}
\textbf{The standard Chaplygin skate.}
This example can be found in \cite{Rosenberg77}. It describes the motion of a
hatchet on a hatchet planimeter, that behaves like a curved knife edge. It is
now commonly known under the name of ``Chaplygin skate''. Let the contact point of
the knife edge have the coordinates $x,y\in\mathbb{R}^2$, let its direction relative
to the positive $x$-axis be $\theta$, and let its center of mass be at distance
$s$ from the contact point. Denote the total mass of the knife edge by  $m$. Thus the moment of inertia about an axis  through the contact point normal to the $xy$
plane is $I=ms^2$.
The configuration space of this problem is the semidirect product
$Q:=\operatorname{SE}(2)=\mathbb{S}^1\,\circledS\, \mathbb{R}^2$ whose coordinates are denoted by $\mathbf{q}: =
(\theta,x,y)$. We have the following concrete constraints on the velocities: 
\[
\mathcal{D}:=\ker(\sin\theta\mathbf{d} x-\cos\theta\mathbf{d} y) = \operatorname{span}\left\{
 \cos\theta\partial_x
+ \sin\theta \partial_y, \partial_\theta \right\} \subset TQ.
\] 
The Lagrangian is hyperregular and taken to be the kinetic energy of the knife
edge, namely,
 \begin{align*}
L(\theta,x,y,\dot\theta,\dot x,\dot y)&=\frac{1}{2}m(\dot
x-s\dot\theta\sin\theta)^2+\frac{1}{2}m(\dot y+s\dot\theta\cos\theta)^2\\
&=\frac{1}{2}m(\dot
x^2+\dot y^2)+\frac{1}{2}ms^2\dot\theta^2+ms\dot\theta(\dot y\cos\theta-\dot x\sin\theta),
\end{align*} 
where we have used that the $x$ and $y$ components of the velocity of the 
center of mass are, respectively, 
\[\dot x-s\dot\theta\sin\theta \quad \text{ and } \quad \dot
y+s\dot\theta\cos\theta.\]
Compute 
\begin{align*}
p_x&=\frac{\partial L}{\partial\dot x}=m\dot x-ms\dot\theta\sin\theta\\
p_y&=\frac{\partial L}{\partial\dot y}=m\dot y+ms\dot\theta\cos\theta\\
p_\theta&=\frac{\partial L}{\partial\dot{\theta} }= ms^2\dot\theta+ms(\dot y\cos\theta-\dot x\sin\theta).
\end{align*}
In $\mathcal{D}$ we have $\dot y\cos\theta-\dot x\sin\theta=0$ and hence we
get for $(\theta,x,y,p_\theta,p_x,p_y)$ in the constraint submanifold $M\subseteq T^*Q$:
\begin{align*}
p_\theta=ms^2\dot\theta \quad \text{ and }\quad
p_x \sin\theta &=m\dot{x} \sin\theta-ms\dot{\theta}\sin^2\theta\\
               &=m\dot{y} \cos\theta-ms\dot{\theta}(1-\cos^2\theta)\\
              &= m\dot{y} \cos\theta+ms\dot{\theta}\cos^2\theta-ms\dot{\theta}\\
              &= p_y \cos\theta -\frac{1}{s}p_\theta.
\end{align*}
Hence the constraint manifold $M$ is five dimensional and given by 
\[
M:=\{(\theta,x,y,p_\theta,p_x,p_y)\mid p_\theta=sp_y \cos\theta -s p_x \sin\theta
\}\subseteq T^*Q.\]
 The global coordinates on $ M $ are thus $(\theta,x,y,p_x,p_y) $. The pull back $ \omega_M $ of the canonical $2$-form $ \omega $ on $ T ^\ast Q $ to  $ M $ has hence the expression
\begin{align*}
\omega_M&= \mathbf{d}  x\wedge\mathbf{d} p_x+\mathbf{d} y\wedge\mathbf{d}
p_y+\mathbf{d} \theta\wedge\mathbf{d}(s p_y \cos\theta -s p_x \sin\theta )\\
&=\mathbf{d}  x\wedge\mathbf{d} p_x+\mathbf{d} y\wedge\mathbf{d}
p_y +s\cos\theta \mathbf{d} \theta\wedge\mathbf{d}p_y-s\sin\theta\mathbf{d} \theta\wedge\mathbf{d}p_x.
\end{align*}
The Dirac structure $D$ modeling this problem is given by \eqref{ym_dirac}. Formula \eqref{h_and_delta_q} gives the vector subbundle
\[
\mathcal{H}:=(T(\pi_{T ^\ast Q}\arrowvert_{M}))^{-1}(\mathcal{D})=\operatorname{span}\{
\cos\theta\partial_x+\sin\theta\partial_y, \partial_\theta, \partial_{p_x}, \partial_{p_y}\} \subset TM,
\]
or equivalently 
\[\mathcal{H}^\circ=\operatorname{span}\{\sin\theta\mathbf{d}
x-\cos\theta\mathbf{d} y\}.\]
A computation yields
\begin{align*}
{\mathbf{i}}_{\cos\theta\partial_x+\sin\theta\partial_y}\omega_M
& =\cos\theta\mathbf{d} p_x+\sin\theta\mathbf{d} p_y \\
{\mathbf{i}}_{\partial_\theta}\omega_M&=s\cos\theta\mathbf{d} p_y-s\sin\theta\mathbf{d} p_x\\
{\mathbf{i}}_{\partial_{p_y}}\omega_M&=-\mathbf{d} y-s\cos\theta\mathbf{d}\theta\\
{\mathbf{i}}_{\partial_{p_x}}\omega_M&=-\mathbf{d} x+s\sin\theta\mathbf{d}\theta.
\end{align*}
Hence
\begin{align*}
\big\{\left(\cos\theta\partial_x+\sin\theta\partial_y,\cos\theta\mathbf{d} p_x+\sin\theta\mathbf{d} p_y\right); \left(\partial_\theta,s\cos\theta\mathbf{d} p_y-s\sin\theta\mathbf{d} p_x \right);\\
  \left(\partial_{p_y}, -\mathbf{d} y-s\cos\theta\mathbf{d}\theta\right);
 \left(\partial_{p_x}, -\mathbf{d} x+s\sin\theta\mathbf{d}\theta\right); \left(0, \sin\theta\mathbf{d} x-\cos\theta\mathbf{d} y\right)\big\} 
\end{align*}
is a smooth global basis for $D$.

We consider the action of the Lie group $G=\operatorname{SE}(2)$ on $Q$, given by
\[\phi:G\times Q\to Q, \quad \phi((\alpha,r,s),(\theta,x,y))=(\theta+\alpha, \cos\alpha
x-\sin\alpha y+r, \sin\alpha x+\cos\alpha y+s).\]
Thus, the induced action on $\Phi:G\times T^*Q\to T^*Q$ is given by 
\begin{align*}
&\Phi((\alpha,r,s),(\theta,x,y,p_\theta,p_x,p_y))\\
=&(\theta+\alpha,\cos\alpha
x-\sin\alpha y+r, \sin\alpha x+\cos\alpha y+s,p_\theta,\cos\alpha
p_x-\sin\alpha p_y, \sin\alpha p_x+\cos\alpha p_y).
\end{align*}

The action on $Q$ obviously leaves the Lagrangian invariant.
We show that the induced action on $T^*Q$ leaves the manifold $M$ invariant: we
denote with $\theta',x', y', p_x',p_y',p_\theta'$ the coordinates of
$\Phi((\alpha,r,s),(\theta,x,y,p_\theta,p_x,p_y))$ and compute
\begin{align*}
s\cos\theta' p_y'-s\sin\theta' p_x'=&s\cos(\theta+\alpha)(\sin\alpha p_x+\cos\alpha p_y)-s\sin(\theta+\alpha)(\cos\alpha
p_x-\sin\alpha p_y)\\
=&s(\cos\theta\cos\alpha-\sin\theta\sin\alpha)(\sin\alpha p_x+\cos\alpha
p_y)\\
&-s(\sin\theta\cos\alpha+\cos\theta\sin\alpha)(\cos\alpha
p_x-\sin\alpha p_y)\\
=&s\cos\theta p_y-s\sin\theta p_x=p_\theta=p_\theta'.
\end{align*}

Since the vertical bundle in this example is
$\mathcal{V}=\operatorname{span}\{\partial_\theta,\partial_x,\partial_y\}$, we have
$\mathcal{V}\cap\mathcal{H}=\operatorname{span}\{\partial_\theta,\cos\theta\partial_x+\sin\theta\partial_y\}$
and $(\mathcal{V}\cap\mathcal{H})^{\omega_M}=\ker\{\cos\theta\mathbf{d}
p_x+\sin\theta\mathbf{d} p_y,s\cos\theta\mathbf{d} p_y-s\sin\theta\mathbf{d} p_x
\}=\ker\{\mathbf{d} p_x,\mathbf{d} p_y\}$.
Hence the distribution
$\mathcal{U}=(\mathcal{V}\cap\mathcal{H})^{\omega_M}\cap\mathcal{H}$ is given
by $\operatorname{span}\{\partial_\theta,\cos\theta\partial_x+\sin\theta\partial_y\}$
and
\begin{eqnarray*}
D\cap \mathcal{K}^\perp=\operatorname{span}\left\{
\left(\cos\theta\partial_x+\sin\theta\partial_y,\cos\theta\mathbf{d} p_x+\sin\theta\mathbf{d} p_y\right),
\left(\partial_\theta,s\cos\theta\mathbf{d} p_y-s\sin\theta\mathbf{d} p_x\right)\right\}.
\end{eqnarray*}
We get the reduced Dirac structure
\begin{eqnarray*} 
 D_{\rm red}&=&\left.\frac{(D\cap \mathcal{K}^\perp)+\mathcal{K}}{\mathcal{K}}\right/G 
=\operatorname{span}\left\{\begin{array}{l}
\left(0,\cos\theta\mathbf{d} p_x+\sin\theta\mathbf{d} p_y\right),
\left(0,s\cos\theta\mathbf{d} p_y-s\sin\theta\mathbf{d} p_x\right)
\end{array}\right\}\\
&=&\operatorname{span}\left\{\left(0,\mathbf{d} p_x\right),
\left(0,\mathbf{d} p_y\right)\right\}
\end{eqnarray*}
on the two dimensional manifold $\bar{M}: = M/G$ with global coordinates
$(p_x,p_y)$. Note that this is the graph of the trivial Poisson tensor  on $\bar
M$.

\medskip

\noindent \textbf{The Chaplygin skate with a rotor on it.}
We propose here a variation of the previous example by considering the Chaplygin
skate with a disk attached to the center of mass of the skate that is  free to rotate about the vertical
axis.  Again, let the contact point of
the knife edge have the coordinates $x,y\in\mathbb{R}^2$, let its direction relative
to the positive $x$-axis be $\theta$, and let its center of mass be at distance
$s$ from the contact point. Denote by $m $ the mass of the knife edge. Thus  its
moment of inertia about an axis  through the contact point normal to the $xy$
plane is $I=ms^2$. Let $\phi$ be the angle between a fixed
point on the disk and the positive $x$-axis and $J$ be the moment of inertia
of the disk about the vertical axis.
The configuration space of this problem is 
$Q:=\mathbb{S}^1\times\mathbb{S}^1\times \mathbb{R}^2$ whose points are denoted by $\mathbf{q}: =
(\phi,\theta,x,y)$. We have again the following concrete constraints on the velocities: 
\[
\mathcal{D}:=\ker(\sin\theta\mathbf{d} x-\cos\theta\mathbf{d} y) = \operatorname{span}\left\{
 \cos\theta\partial_x
+ \sin\theta \partial_y, \partial_\theta \right\} \subset TQ.
\] 
The Lagrangian is the kinetic energy of the knife
edge: \begin{align*}
L(\phi,\theta,x,y,\dot\phi,\dot\theta,\dot x,\dot y)&=\frac{1}{2}m(\dot
x-s\dot\theta\sin\theta)^2+\frac{1}{2}m(\dot y+s\dot\theta\cos\theta)^2+\frac{1}{2}J(\dot\theta+\dot\phi)^2\\
&=\frac{1}{2}m(\dot
x^2+\dot y^2)+\frac{1}{2}(I+J)\dot\theta^2+ms\dot\theta(\dot y\cos\theta-\dot x\sin\theta)+\frac{1}{2}J\dot\phi^2+J\dot\phi\dot\theta.
\end{align*} 
Compute 
\begin{align*}
p_x&=m\dot x-ms\dot\theta\sin\theta\\
p_y&=m\dot y+ms\dot\theta\cos\theta\\
p_\theta&=(I+J)\dot\theta+ms(\dot y\cos\theta-\dot x\sin\theta)+J\dot\phi\\
p_\phi&=J(\dot\phi+\dot\theta).
\end{align*}
Again, if we have $\dot y\cos\theta-\dot x\sin\theta=0$, we compute:
\begin{align*}
p_\theta=(I+J)\dot\theta+J\dot\phi=I\dot\theta+J(\dot\theta+\dot\phi)=I\dot\theta+p_\phi=ms^2\dot\theta+p_\phi\end{align*}
and
\begin{align*}
p_x \sin\theta &=m\dot{x} \sin\theta-ms\dot{\theta}\sin^2\theta
=m\dot{y} \cos\theta -ms\dot{\theta}(1-\cos^2\theta)\\
              &= m\dot{y} \cos\theta+ms\dot{\theta}\cos^2\theta-ms\dot{\theta}
              = p_y \cos\theta +\frac{1}{s}(p_\phi-p_\theta).
\end{align*}
Hence the constraint manifold $M$ is seven dimensional and given by 
\[
M:=\{(\phi, \theta, x,y, p_ \phi, p_\theta,p_x,p_y)\mid p_\theta=s p_y \cos\theta -s p_x \sin\theta
+p_\phi\}\subseteq T^*Q.\]
 The global coordinates on $ M $ are thus $(\phi,\theta,x,y,p_\phi, p_x,p_y) $. The pull back $ \omega_M $ of the canonical $2$-form $ \omega $ on $ T ^\ast Q $ to  $ M $ has hence the expression
\begin{align*}
\omega_M&= \mathbf{d}  x\wedge\mathbf{d} p_x+\mathbf{d} y\wedge\mathbf{d}
p_y+\mathbf{d} \theta\wedge\mathbf{d}(s p_y \cos\theta -s p_x \sin\theta +
p_\phi)+\mathbf{d}
\phi\wedge\mathbf{d} p_\phi\\
&=\mathbf{d}  x\wedge\mathbf{d} p_x+\mathbf{d} y\wedge\mathbf{d}
p_y+s\cos\theta \mathbf{d} \theta\wedge\mathbf{d}p_y-s\sin\theta\mathbf{d}
\theta\wedge\mathbf{d}p_x+(\mathbf{d}\theta+\mathbf{d}
\phi)\wedge\mathbf{d} p_\phi.
\end{align*}
The Dirac structure $D$ modeling this problem is given by \eqref{ym_dirac}. Formula \eqref{h_and_delta_q} gives the vector subbundle
\[
\mathcal{H}:=(T(\pi_{T ^\ast Q}\arrowvert_{M}))^{-1}(\mathcal{D})=\operatorname{span}\{\partial_\phi, \partial_\theta,
\cos\theta\partial_x+\sin\theta\partial_y, \partial_{p_\phi},
\partial_{p_x}, \partial_{p_y}\} \subset TM,
\]
or equivalently 
\[\mathcal{H}^\circ=\operatorname{span}\{\sin\theta\mathbf{d} x-\cos\theta\mathbf{d} y\}.\]
A computation yields
\begin{align*}
{\mathbf{i}}_{\partial_\phi}\omega_M&=\mathbf{d} p_\phi\\
{\mathbf{i}}_{\partial_\theta}\omega_M&=s\cos\theta\mathbf{d} p_y-s\sin\theta\mathbf{d}
p_x+\mathbf{d} p_\phi\\
{\mathbf{i}}_{\cos\theta\partial_x+\sin\theta\partial_y}\omega_M&=\cos\theta\mathbf{d}
p_x+\sin\theta\mathbf{d} p_y\\
{\mathbf{i}}_{\partial_{p_\phi}}\omega_M&=-\mathbf{d}\theta-\mathbf{d}
\phi\\
{\mathbf{i}}_{\partial_{p_y}}\omega_M&=-\mathbf{d} y-s\cos\theta\mathbf{d}\theta\\
{\mathbf{i}}_{\partial_{p_x}}\omega_M&=-\mathbf{d} x+s\sin\theta\mathbf{d}\theta.
\end{align*}
We get
\begin{align*}
D=\operatorname{span}\big\{\left(\partial_\phi,\mathbf{d} p_\phi\right);\left(\partial_\theta,s\cos\theta\mathbf{d} p_y-s\sin\theta\mathbf{d}
p_x+\mathbf{d} p_\phi\right);\left(\cos\theta\partial_x+\sin\theta\partial_y,\cos\theta\mathbf{d} p_x+\sin\theta\mathbf{d} p_y\right); \\
 \left(\partial_{p_\phi}, -\mathbf{d}\theta-\mathbf{d}
\phi\right); \left(\partial_{p_y}, -\mathbf{d} y-s\cos\theta\mathbf{d}\theta\right);
 \left(\partial_{p_x}, -\mathbf{d} x+s\sin\theta\mathbf{d}\theta\right); \left(0, \sin\theta\mathbf{d} x-\cos\theta\mathbf{d} y\right)\big\}. 
\end{align*}

We consider the action of the Lie group $G=\mathbb{S}^1\times \operatorname{SE}(2)$ on $Q$, given by
\[\phi:G\times Q\to Q, \quad \phi((\beta,\alpha,r,s),(\phi,\theta,x,y))=(\phi+\beta,\theta+\alpha, \cos\alpha
x-\sin\alpha y+r, \sin\alpha x+\cos\alpha y+s).\]
Thus, the induced action $\Phi:G\times T^*Q\to T^*Q$ on $T ^\ast Q $ is given by 
\begin{align*}
&\Phi((\beta,\alpha,r,s),(\phi,\theta,x,y,p_\phi,p_\theta,p_x,p_y))\\
&=(\phi+\beta,\theta+\alpha,\cos\alpha
x-\sin\alpha y+r, \sin\alpha x+\cos\alpha y+s,p_\theta,\cos\alpha
p_x-\sin\alpha p_y, \sin\alpha p_x+\cos\alpha p_y).
\end{align*}

The Lagrangian is invariant under the lift to $TQ$ of $\phi$ and it is easy to see,
with the considerations in the previous example, that the induced action $\Phi$ on $T^*Q$ leaves the manifold $M$ invariant.

Since the vertical bundle in this example is
$\mathcal{V}=\operatorname{span}\{\partial_\phi,\partial_\theta,\partial_x,\partial_y\}$, we have
$\mathcal{V}\cap\mathcal{H}=\operatorname{span}\{\partial_\phi,\partial_\theta,\cos\theta\partial_x+\sin\theta\partial_y\}$
and $(\mathcal{V}\cap\mathcal{H})^{\omega_M}=\ker\{\mathbf{d} p_\phi,\cos\theta\mathbf{d}
p_x+\sin\theta\mathbf{d} p_y,s\cos\theta\mathbf{d} p_y-s\sin\theta\mathbf{d} p_x+\mathbf{d} p_\phi
\}=\ker\{\mathbf{d} p_\phi,\mathbf{d} p_x,\mathbf{d} p_y\}$.
Hence the distribution
$\mathcal{U}=(\mathcal{V}\cap\mathcal{H})^{\omega_M}\cap\mathcal{H}$ is given
by $\operatorname{span}\{\partial_\phi,\partial_\theta,\cos\theta\partial_x+\sin\theta\partial_y\}$
and
\begin{eqnarray*}
D\cap \mathcal{K}^\perp=\operatorname{span}\left\{\left(\partial_\phi,\mathbf{d} p_\phi\right),
\left(\cos\theta\partial_x+\sin\theta\partial_y,\cos\theta\mathbf{d} p_x+\sin\theta\mathbf{d} p_y\right),
\left(\partial_\theta,s\cos\theta\mathbf{d} p_y-s\sin\theta\mathbf{d} p_x+\mathbf{d} p_\phi\right)\right\}.
\end{eqnarray*}
We get the reduced Dirac structure
\begin{eqnarray*}
 D_{\rm red}&=&\left.\frac{(D\cap \mathcal{K}^\perp)+\mathcal{K}}{\mathcal{K}}\right/G \\
&=&\operatorname{span}\left\{\begin{array}{l}\left(0,\mathbf{d} p_\phi\right),
\left(0,\cos\theta\mathbf{d} p_x+\sin\theta\mathbf{d} p_y\right),
\left(0,s\cos\theta\mathbf{d} p_y-s\sin\theta\mathbf{d} p_x+\mathbf{d} p_\phi\right)
\end{array}\right\}\\
&=&\operatorname{span}\left\{\left(0,\mathbf{d} p_\phi\right),\left(0,\mathbf{d} p_x\right),
\left(0,\mathbf{d} p_y\right)\right\}
\end{eqnarray*}
on the three dimensional manifold $\bar{M}: = M/G$ with global coordinates
$(p_\phi,p_x,p_y)$. This is again the graph of the trivial Poisson tensor on $\bar
M$.

\medskip

In these six examples we get integrable Dirac structures after
reduction. We shall come back to this remark in the last section of the paper.

\section{Optimal reduction for nonholonomic systems}
\label{nonholoptred}
Recall the setting of \S \ref{nonhol}: $Q$ is a configuration space which is a smooth
Riemannian manifold, $\mathcal{D}\subseteq TQ$ is the constraints distribution given
as the intersection of the kernels of $k $ linearly independent $1$-forms on $Q $ (and
$\mathcal{D}$ is hence a vector subbundle of $TQ$), $L $ is a hyperregular classical Lagrangian equal to the
kinetic energy of the given Riemannian metric on $Q $ minus a potential, $M : =
\mathbb{F}L( \mathcal{D}) \subset T^*Q$ is a submanifold and represents the
constraints in phase space $T ^\ast Q $, and 
$\omega_M:=i^*\omega_{\rm can} \in \Omega^2(M)$ is the induced $2$-form on $M $,
where $i:M\hookrightarrow  T^*Q$ is the inclusion and
$\omega_{\rm can}$ the canonical symplectic form on $T^*Q$. The
vector bundle
 $\mathcal{H}:=TM\cap (T\pi_{T^*Q})^{-1}(\mathcal{D})$ is not integrable
but has the property that the restriction $\omega_\mathcal{H}$ of $\omega_M$ to
$\mathcal{H} \times \mathcal{H}$ is
nondegenerate. The Dirac structure $D$ associated to this nonholonomic system has
fibers 
\[
D(m) = \{(X(m),\alpha_m)\in T_m M \oplus T^\ast_mM \mid X\in \Gamma(\mathcal{H}),\,
\alpha-{\mathbf{i}}_{X}\omega_M\in \Gamma(\mathcal{H}^\circ)\}
\] 
for all $m\in M$ and is, in general, not integrable. Recall from Proposition
\ref{prop_nonhol_dirac}(i) that $\mathsf{G}_0=\{0\}$,
$\mathsf{P}_1= T^*M$ and hence \emph{all} functions
are admissible.
\medskip

Consider a $G$-action $\phi:G\times Q\to Q$ on $Q$ that leaves the constraints and
the Lagrangian invariant. The lift $\Phi:G\times T^*Q\to
T^*Q$ of the action is defined by  $\Phi_g=(T\phi_{g^{-1}})^*$; this is a symplectic
action on $T^*Q$ that leaves $M$ invariant. Thus we get a canonical $G$-action on the
Dirac manifold $(M,D)$ and we have for all $g\in G$,
 \[
 \Phi^*_g\omega_M=\Phi^*_g \left(i^*\omega_{\rm can} \right) =
 i^*\left(\Phi^*_g\omega_{\rm can} \right) =i^*\omega_{\rm can}=\omega_M
 \] 
since the $G$-action commutes with the inclusion. Note that in this section the $G
$-action on $T ^\ast Q $ is a lift, whereas in \S\ref{nonhol} we needed only that it
is a symplectic action.
\medskip

In this section we shall define a distribution on $M$ that yields the equations of
motion and the conserved quantities given by the Nonholonomic Noether Theorem
(see \cite{BaCuKeSn95}, Theorem 2 and also \cite{Bloch03}, Chapter 5 and the corresponding internet supplement). If this
distribution is integrable, we will prove a Marsden-Weinstein reduction theorem that
gives a reduced Dirac structure which is the graph of a nondegenerate $2$-form (not
necessarily closed). This reduction procedure is done from an ``optimal'' point of
view as in \cite{JoRa10b}.

\subsection{The nonholonomic Noether theorem}\label{sec:nonholonomic_Noether}

We recall in this subsection the Hamiltonian formulation of the Nonholonomic Noether
Theorem. Let  $\mathbf{J}:T^*Q\to\mathfrak{g}^*$ be the canonical momentum map
associated to the action of
$G$ on $T^*Q$ (see, e.g., \cite{MaRa99}) 
\begin{equation}\label{jxi1}
\mathbf{J}(p)(\xi)=\langle p\,,\xi_Q(\pi(p))\rangle
\end{equation} for all $p\in T^*Q$, 
 where
$\pi:T^*Q\to Q$ is the projection. For all $\xi\in\mathfrak{g}$, the $\xi$-component
of $\mathbf{J}$ is the map
$\mathbf{J}^\xi:T^*Q\to\mathbb{R}$ defined by 
\begin{equation}\label{jxi2}
\mathbf{J}^\xi(p):=
\mathbf{J}(p)(\xi)
\end{equation} 
for all $p\in T ^\ast Q$. We shall denote by the same symbol $\mathbf{J}^\xi$ its restriction
to the manifold $M$. For an
arbitrary $\xi\in\mathfrak{g}$ we have therefore
\begin{equation}\label{noether}
{\mathbf{i}}_{\xi_{T^*Q}}\omega_{\rm can}=\mathbf{d}
\mathbf{J}^\xi.
\end{equation}
 Since the action of $G$ on $T^*Q$ leaves the submanifold $M$
invariant, we have $\xi_{T^*Q}(m)\in T_mM$ for all $m\in M$ and hence the
fundamental vector field $\xi_{T^*Q}$ is $i$-related to $\xi_M$,
i.e., $Ti\circ\xi_M=\xi_{T^*Q}\circ i$.
Choosing for each vector field $X\in \mathfrak{X}(M)$ an arbitrary extension
$X'\in\mathfrak{X}(T^*Q)$ (and hence $X\sim_{i}X'$) we get for all $m\in M$,
\begin{align*}
{\mathbf{i}}_{\xi_M}\omega_M(X)(m)&={\mathbf{i}}_{\xi_M}(i^*\omega_{\rm
can})(X)(m)={\mathbf{i}}_{\xi_{T^*Q}}\omega_{\rm can}(X')(i(m))
=(\mathbf{d} \mathbf{J}^\xi(X'))(i(m)) \\
&=(i^*\mathbf{d} \mathbf{J}^\xi)(X)(m)
=(\mathbf{d} \mathbf{J}^\xi(X))(m)
\end{align*}
which shows that \eqref{noether} naturally restricts to $M$ 
\begin{equation}\label{resnoe}
{\mathbf{i}}_{\xi_M}\omega_M=\mathbf{d}
\mathbf{J}^\xi.
\end{equation}

Define for all $p\in M$ the vector subspace
$\mathfrak{g}^p:=\{\xi\in \mathfrak{g} \mid \xi_M(p)\in
(\mathcal{V}\cap\mathcal{H})(p)\}\subseteq\mathfrak{g}$. Then  
\[
\mathfrak{g}^{\mathcal{H}}:=\bigcup_{p\in M}\mathfrak{g}^p
\]
is a smooth (not necessarily trivial) vector subbundle of the trivial bundle
$M \times \mathfrak{g}$ if and only if $\mathcal{H}+\mathcal{V}$ has constant
rank on $M $, for instance if $\mathcal{H}+\mathcal{V}=TM$. Indeed, note first
that $\mathfrak{g}^{\mathcal{H}} = \Lambda^{-1}( \mathcal{V}\cap \mathcal{H})
$, where $\Lambda: M \times \mathfrak{g} \rightarrow \mathcal{V}$ is the
vector bundle isomorphism over $ M $ given by  $\Lambda(m , \xi): =
\xi_M(m)$. 
However, since $\mathcal{H}+\mathcal{V}$, $\mathcal{H}$, and
$\mathcal{V}$ are subbundles of $TM$, it follows that $\mathcal{V}\cap \mathcal{H} $ is also a 
subbundle of both $TM$ and $\mathcal{V}$.  Consequently, $\mathfrak{g}^\mathcal{H}= \Lambda^{-1}( \mathcal{V}\cap \mathcal{H})
$ is a subbundle of the trivial vector bundle $M\times \mathfrak{g}$. 
If $\mathfrak{g}^\mathcal{H}$ is a vector bundle over $M$, then
$\mathcal{V}\cap \mathcal{H} $ is also a vector bundle and hence its fibers
have constant dimension on $M$. It follows immediately that the rank of
$\mathcal{H}+\mathcal{V}$ is also constant on $M$.

\textit{For the rest of this subsection we assume that $\mathcal{H}+ \mathcal{V}$ has constant rank on 
$M$ and hence that $\mathfrak{g}^ \mathcal{H} $ is a vector subbundle of the trivial vector bundle $M 
\times  \mathfrak{g}$.}
If $\xi^{\mathcal{H}}$ is a smooth section of $\mathfrak{g}^{\mathcal{H}}$, then
$\boldsymbol{\xi}(p):=(\xi^{\mathcal{H}}(p))_M(p)$ defines a smooth section of
$\mathcal{V}\cap\mathcal{H}$.
Conversely, if $\{\xi^1,\dots,\xi^k\}$
is a chosen basis for the Lie algebra $\mathfrak{g}$, then the vector fields
$\xi^1_M,\dots, \xi^k_M$ are global vector fields on $M$ that don't vanish
and are everywhere linearly independent. Hence, $\xi^1_M,\dots, \xi^k_M$
are smooth basis vector fields for the bundle $\mathcal{V}$.
Every section
$\boldsymbol{\xi}$ of $\mathcal{V}\cap\mathcal{H}$ can hence be written
$\boldsymbol{\xi}=\sum_{i=1}^kf_i\xi^i_M$ with smooth (local) functions
$f_1,\dots,f_k$, and corresponds exactly to the section
$\xi^{\mathcal{H}}=\sum_{i=1}^kf_i\xi^i$ of $\mathfrak{g}^\mathcal{H}$. 

Since  $\mathcal{V}\cap \mathcal{H} $ is a
subbundle of $TM$ with
\[[\Gamma(\mathcal{V}\cap\mathcal{H}),\Gamma(\mathcal{V})]\subseteq\Gamma(\mathcal{V})=
\Gamma((\mathcal{V}\cap\mathcal{H})+\mathcal{V}),\]
we that for each
$p\in M$ there exists a neighborhood $U$ of $p$ and 
 spanning sections of  $\mathcal{V}\cap \mathcal{H} $ on
$U$ that descend to the quotient $M/G$ (see \cite{JoRaZa11}). 

Let $\xi^{\mathcal{H}}$ be a smooth section of $\mathfrak{g}^{\mathcal{H}}$.  For all $p\in M$
and all $X\in \mathfrak{X}(M)$ the definition of the corresponding $\boldsymbol{\xi}$ and
\eqref{resnoe} yield
\begin{equation}
\label{xi_h_omega_equ}
\omega_M(p)\left(\boldsymbol{\xi}(p),X(p)\right)=\mathbf{d}
\mathbf{J}^{\xi^{\mathcal{H}}(p)}(p)\left(X\right).
\end{equation}
As above, write $\xi^{\mathcal{H}}=\sum_{i=1}^kf_i\xi^i$ with  smooth
functions $f_1,\dots,f_k$ and the chosen basis
$\{\xi^1,\dots,\xi^k\}$ of $\mathfrak{g}$.
 Define the smooth map 
\begin{align*}
\begin{array}{llcl}
\mathbf{J}^{\xi^{\mathcal{H}}}:&M&\to & \mathbb{R}\\
& p&\mapsto &\mathbf{J}^{\xi^{\mathcal{H}}(p)}(p)=\langle
\mathbf{J}(p),\xi^{\mathcal{H}}(p)\rangle.
\end{array}
\end{align*}
Using \eqref{jxi1} and \eqref{jxi2} we get
\[
\mathbf{J}^{\xi^{\mathcal{H}}}(p)=\mathbf{J}^{\sum_{i=1}^kf_i(p)\xi^i}(p)=\sum_{i=1}^
kf_i(p)\mathbf{J}^{\xi^i}(p).
\]
If $c:(-\varepsilon;\varepsilon)\to M$ is a solution curve of a vector field
$X\in \mathfrak{X}(M)$ with $c(0)=p$, we have 
\begin{align*}
\mathbf{d} \mathbf{J}^{\xi^{\mathcal{H}}}(p)(X)
&=\left.\frac{d}{dt}\right\arrowvert_{t=0}\mathbf{J}^{\xi^{\mathcal{H}}}(c(t))
=\left.\frac{d}{dt}\right\arrowvert_{t=0}\sum_{i=1}^kf_i(c(t))\mathbf{J}^{\xi^i}(c(t))\\
&=\sum_{i=1}^k\mathbf{d} f_i(p)(\dot{c}(0))\mathbf{J}^{\xi^i}(p)
+\sum_{i=1}^kf_i(p)\mathbf{d} \mathbf{J}^{\xi^i}(p)(\dot{c}(0)) \\
&=\mathbf{J}^{\sum_{i=1}^k{\mathbf{d} f_i}(p)(X)\xi^i}(p)
+\mathbf{d} \mathbf{J}^{\xi^{\mathcal{H}}(p)}(p)(X)\\
&=\mathbf{J}^{X\left[\xi^{\mathcal{H}}\right]}(p)
+\mathbf{d} \mathbf{J}^{\xi^{\mathcal{H}}(p)}(p)(X),
\end{align*} 
where we write $X\left[\xi^{\mathcal{H}}\right]: = \sum_{i=1}^kX[f_i]\xi^i$. Thus
\eqref{xi_h_omega_equ} becomes for all $p\in M$ and all $X\in \mathfrak{X}(M)$,
\begin{align*}
\omega_M(p)\left(\boldsymbol{\xi}(p),X(p)\right)&=\mathbf{d}
\mathbf{J}^{\xi^{\mathcal{H}}}(p)(X)
- \mathbf{J}^{X\left[\xi^{\mathcal{H}}\right]}(p).
\end{align*} 
Hence, if the one form $\alpha^{\boldsymbol{\xi}}\in \Omega^1(M)$ is defined by
\begin{equation*}
\alpha^{\boldsymbol{\xi}}(X) : =\mathbf{d} \mathbf{J}^{\xi^{\mathcal{H}}}(X)
- \mathbf{J}^{X\left[\xi^{\mathcal{H}}\right]}
\end{equation*}
for all $X\in \mathfrak{X}(M)$, we have
${\mathbf{i}}_{\boldsymbol{\xi}}\omega_M=\alpha^{\boldsymbol{\xi}}$ and so the pair
$(\boldsymbol{\xi},\alpha^{\boldsymbol{\xi}})$ is a section of $D$. 

Let $h $ be a $G$-invariant Hamiltonian and
$X_h\in\Gamma(\mathcal{H})$ the solution of the implicit Hamiltonian system 
$(X,\mathbf{d} h)\in\Gamma(D)$. Then 
\begin{align*}
\mathbf{d} \mathbf{J}^{\xi^{\mathcal{H}}}(X_h)-\mathbf{J}^{X_h\left[\xi^{\mathcal{H}}
\right]}=\alpha^{\boldsymbol{\xi}}(X_h)=\omega_M(\boldsymbol{\xi}, X_h)=-\mathbf{d} h
(\boldsymbol{\xi})=0
\end{align*}
since $\mathbf{d}h (\boldsymbol{\xi})(p) = \left\langle \mathbf{d}h(p),
  \boldsymbol{\xi}(p) \right\rangle = \left\langle \mathbf{d}h(p), \left(\xi^
    {\mathcal{H}}(p) \right)_M(p) \right\rangle = 0 $ by $G $-invariance of $h
$.
Thus, we have proved the following result.
\begin{theorem}\label{nonholonomic_theorem}
Let $\xi^{\mathcal{H}}$ be a section of
  $\mathfrak{g}^{\mathcal{H}}$ and $X_h\in\Gamma(\mathcal{H})$ the solution of the
implicit Hamiltonian system 
$(X,\mathbf{d} h)\in\Gamma(D)$, where  $h $ is a $G$-invariant
Hamiltonian. Then $X_h$ satisfies the \emph{Nonholonomic Noether Momentum Equation}:
\begin{equation}
\label{nonhol_Noether_eq}
\mathbf{d}
\mathbf{J}^{\xi^{\mathcal{H}}}(X_h)-\mathbf{J}^{X_h\left[\xi^{\mathcal{H}}
\right]}=0.
\end{equation}
\end{theorem}

Recall from \eqref{def_M} and \eqref{def_cal_H} that $\mathcal{H}$ is defined
in terms of the given Lagrangian $L : TQ\rightarrow \mathbb{R}$ and hence, only the dynamics defined by the 
corresponding Hamiltonian $H$ is of interest. Thus, for each other Lagrangian $L'$ we obtain another 
distribution $\mathcal{H}'$.

\begin{remark} 
In \cite{Bloch03}, Theorem 5.5.4, the Nonholonomic Noether Theorem is formulated in terms of a Lagrangian of a classical mechanical systems (hence equal to the kinetic energy of a metric minus a potential). Let
$\mathcal{V}_Q\subseteq TQ$ be the vertical subbundle of the action 
$\phi: G \times  Q \rightarrow  Q$. Under the \emph{Dimension Assumption} $\mathcal{D}+\mathcal{V}_Q= TQ$,
the distribution $\mathcal{D}\cap\mathcal{V}_Q$ is a smooth subbundle of
$TQ$. Note that this assumption leads automatically to
$\mathcal{F}+\mathcal{V}_{T^*Q}=TT^*Q$ and hence to $\mathcal{H}+\mathcal{V}=TM$
 (see \eqref{def_F} for the definition of $ \mathcal{F} $).

The smooth vector bundle $\mathfrak{g}^{\mathcal{D}}:=\bigcup_{p\in M}\mathfrak{g}^{\mathcal{D}}(q)$ is defined pointwise by
\[\mathfrak{g}^{\mathcal{D}}(q):=\{\xi\in \mathfrak{g} \mid \xi_Q(q)\in
(\mathcal{V}_Q\cap\mathcal{D})(q)\}\subseteq\mathfrak{g}.\]
Let $(\mathfrak{g}^{\mathcal{D}})^*$ be the dual bundle, that is, its fibers are  $(\mathfrak{g}^{\mathcal{D}})^*(q):=(\mathfrak{g}^{\mathcal{D}}(q))^*$ for all $q\in Q$.
The \emph{nonholonomic momentum map} $J^{\rm nhc}: TQ \rightarrow  (\mathfrak{g}^{\mathcal{D}})^*$ is the vector bundle map over $ Q $ defined by 
\[\langle J^{\rm nhc}(v_q),\xi\rangle=
\langle
\mathbb{F}L(v_q),\xi_Q(q)\rangle=
\frac{\partial L}{\partial\dot{q}^i}(\xi_Q)^i(q)=:J^{\rm nhc}(\xi)(v_q)\]
where $\xi\in \mathfrak{g}^{\mathcal{D}}(q)$. Let $\xi^{\mathcal{D}}$ be a section
of the bundle $\mathfrak{g}^{\mathcal{D}}$. Theorem 5.5.4 in \cite{Bloch03} states
that  any solution $c(t)=(q(t),\dot{q}(t))$ of the
Lagrange-d'Alembert equations for a nonholonomic system must satisfy, in
addition to the given kinematic constraints, the \emph{momentum equation} 
\begin{equation}\label{lagrange_momentum_eq}
\frac{d}{dt} J^{\rm nhc}\left(\xi^{\mathcal{D}}(q(t))\right)(c(t))=
\frac{\partial L
}{\partial\dot{q}^i}\left[\frac{d}{dt}(\xi^{\mathcal{D}}(q(t))\right]^i_Q.
\end{equation} 
Since for all $\xi\in \mathfrak{g}$ the vector field $\xi_{T^*Q}$ is the cotangent lift of $\xi_Q$, we have in local charts
\[\xi_{T^*Q}=\xi_Q^i\frac{\partial}{\partial q^i} -
\frac{\partial\xi_Q^i}{\partial q^j}p_i\frac{\partial}{\partial p_j}.\]
Hence, if $\xi_Q(q)\in\mathcal{D}(q)$, we get
$\xi_{T^*Q}(\alpha_q)\in\mathcal{F}(\alpha_q)$ for all $\alpha_q\in T_q^*Q$.

Consequently $\boldsymbol{\xi}(p):=\left(\xi^{\mathcal{D}}(\pi_{T^*Q}(p))\right)_M(p)$ for all $p\in M$ defines a smooth section $\boldsymbol{\xi}$ of $\mathcal{V}\cap\mathcal{H}$ and hence a smooth section 
$\xi^{\mathcal{H}}$ of $\mathfrak{g}^{\mathcal{H}}$. Note that
$\xi^{\mathcal{H}}(p)=\xi^{\mathcal{D}}(\pi_{T^*Q}(p))$ for all $p\in
M$ and if  $\xi^{\mathcal{D}}=\sum_{i=1}^kf_i\xi^i_Q$ with smooth functions $f_1,\dots,f_k$, then $\xi^{\mathcal{H}}=\sum_{i=1}^kF_i\xi^i_M$ with the smooth functions $F_i$ defined by  
$F_i=i_M^*\pi_{T^*Q}^*f_i$, where $i_M:M \hookrightarrow T ^\ast Q $ is the inclusion. Let $X_H$ be a solution of the implicit Hamiltonian system 
$(X,\mathbf{d} H)\in\Gamma(D)$, where $H$ is  the $G$-invariant Hamiltonian on $M$ associated to the Lagrangian $L$ by the Legendre transformation, and $p(t)$ an integral curve of $X_H$. Then
$c(t):=(q(t),\dot{q}(t))=(\mathbb{F}L)^{-1}(p(t))$ is a solution of the
Lagrange-d'Alembert equations. We have for all $t$
\begin{align*}
0&=\left(\mathbf{d} \mathbf{J}^{\xi^{\mathcal{H}}}(X_H)-\mathbf{J}^{X_H\left[\xi^{\mathcal{H}}\right]}\right) (p(t))=
\frac{d}{dt}\mathbf{J}^{\xi^{\mathcal{H}}}(p(t))
-\sum_{i=1}^k\frac{d}{dt}(F_i(p(t)))\mathbf{J}^{\xi_i}(p(t))\\
&=\frac{d}{dt}\left\langle
  p(t),\left(\xi^{\mathcal{H}}(p(t))\right)_Q(q(t))\right\rangle-\left\langle
  p(t),\left(\frac{d}{dt}\big(F_i(p(t))\big)\xi^i\right)_Q(q(t))\right\rangle\\
&= \frac{d}{dt}\left\langle
 \mathbb{F}L(c(t)) ,\left(\xi^{\mathcal{D}}(q(t))\right)_Q(q(t))\right\rangle-\left\langle
 \mathbb{F}L(c(t))
 ,\left(\frac{d}{dt}\big(f_i(q(t))\big)\xi^i\right)_Q(q(t))\right\rangle\\
&= \frac{d}{dt}\left\langle
 \mathbb{F}L(c(t)) ,\left(\xi^{\mathcal{D}}(q(t))\right)_Q(q(t))\right\rangle-\left\langle
 \mathbb{F}L(c(t))
 ,\left(\frac{d}{dt}\xi^{\mathcal{D}}(q(t))\right)_Q(q(t))\right\rangle\\
&=\frac{d}{dt} J^{\rm nhc}\left(\xi^{\mathcal{D}}(q(t))\right)(c(t))-
\frac{\partial L
}{\partial\dot{q}^i}\left[\frac{d}{dt}(\xi^{\mathcal{D}}(q(t))\right]^i_Q.
\end{align*}
Hence our  Nonholonomic Noether Theorem \ref{nonholonomic_theorem} is the Hamiltonian version of 
Theorem 5.5.4 in \cite{Bloch03}, that is, \eqref{nonhol_Noether_eq} and \eqref{lagrange_momentum_eq} 
are equivalent.
\end{remark}

\begin{proposition} \label{prop_73}
Assume that  $\mathcal{V}+\mathcal{H}$ has constant rank on $M$. Let
$\xi^{\mathcal{H}}$ be
  a $G$-equivariant section of $\mathfrak{g}^{\mathcal{H}}$. Then the corresponding
section
$(\boldsymbol{\xi},\alpha^{\boldsymbol{\xi}})$ of $D$ is also $G$-equivariant.
 There
are two possibilities:
\begin{itemize}
\item[{\rm (i)}]
  $\alpha^{\boldsymbol{\xi}}={\mathbf{i}}_{\boldsymbol{\xi}}\omega_M=0$ on
$\mathcal{V}\cap\mathcal{H}$. Then there exist $\alpha'\in\Gamma(\mathcal{V}^\circ)$
such that
$(\boldsymbol{\xi},\alpha')$ is a $G$-equivariant section of
$D\cap\mathcal{K}^\perp$ and  exactly one section
$\bar\alpha\in\Gamma(\mathsf{P}_0^{\rm red})$ such
that $\pi^*\bar\alpha=\alpha'$. Conversely, each section of
$\Gamma(\mathsf{P}_0^{\rm red})$ pulls back to a section $\alpha'$ defined as
above and satisfying this condition.
\item[{\rm (ii)}] $\mathcal{V}\cap\mathcal{H}\nsubseteq
  \boldsymbol{\xi}^{\omega_M}$ and hence $\alpha^{\boldsymbol{\xi}}\neq0$ on
$\mathcal{V}\cap\mathcal{H}$. Then
$\alpha^{\boldsymbol{\xi}}$ leads to a momentum
equation  that doesn't appear in the reduced implicit Hamiltonian system.
\end{itemize}
\end{proposition}

\begin{proof}
If $\xi^{\mathcal{H}}$ is $G$-equivariant, we have $\xi^{\mathcal{H}}(g\cdot
p)=\operatorname{Ad}_g\xi^{\mathcal{H}}(p)$ for all $p\in M$ and hence we get for the corresponding
$\boldsymbol{\xi}$ 
\begin{align*}
(\Phi_g^*(\boldsymbol{\xi}))(p)&=T_{g\cdot
  p}\Phi_{g^{-1}}\boldsymbol{\xi}(g\cdot p)
=T_{g\cdot
  p}\Phi_{g^{-1}}\left(\xi^{\mathcal{H}}(g\cdot p)\right)_M(g\cdot p)\\
&=T_{g\cdot  p}\Phi_{g^{-1}}\left(\operatorname{Ad}_g(\xi^{\mathcal{H}}(p))\right)_M(g\cdot p)\\
&=\left(\operatorname{Ad}_{g^{-1}}\circ\operatorname{Ad}_g(\xi^{\mathcal{H}}(p))\right)_M(p)=\left(\xi^{\mathcal{
H}}(p)\right)_M(p)=\boldsymbol{\xi}(p).
\end{align*}
Note that conversely, if $\boldsymbol{\xi}$ is equivariant, then the
corresponding section $\xi^{\mathcal{H}}$ of $\mathfrak{g}^{\mathcal{H}}$ is
$G$-equivariant. Since $\Phi_g^*\omega_M=\omega_M$ for all $g\in G$, the
section $(\boldsymbol{\xi},\alpha^{\boldsymbol{\xi}})$ is
$G$-equivariant. Since $\mathcal{V}+\mathcal{H}=TM$, if $\alpha^{\boldsymbol{\xi}}=0$
on
$\mathcal{V}\cap\mathcal{H}$ there exists as in Proposition
\ref{prop_nonhol_dirac}{\rm (ii)} a  section
$\beta\in\Gamma(\mathcal{H}^\circ)$ such that
$\alpha^{\boldsymbol{\xi}}+\beta\in\Gamma(\mathcal{V}^\circ)$
and hence
$({\boldsymbol{\xi}},\alpha^{\boldsymbol{\xi}}+\beta)\in\Gamma(D\cap\mathcal{K}^\perp)$.
Since $\boldsymbol{\xi}$, $\alpha^{\boldsymbol{\xi}}$, $\mathcal H^\circ$, $\mathcal V^\circ$ 
are all $G$-invariant and the action is free, we can take $\beta$ $G$-invariant.
 Hence the first statement of {\rm (i)} holds
with $\alpha':=\alpha^{\boldsymbol{\xi}}+\beta$. But because
$\boldsymbol{\xi}\in\Gamma(\mathcal{V})$, the section of $D_{\rm red}$
corresponding to $(\boldsymbol{\xi},\alpha')$ will be $(0,\bar\alpha)$ with
$\bar\alpha\in\Omega^1(\bar M)$ such that $\pi^*\bar\alpha=\alpha'$.

On the
other hand, if we choose a non-zero section $\bar\alpha$ of $\mathsf{P}_0^{\rm
red}$,
the codistribution associated to the reduced Dirac structure on $\bar{M}$, we
have $(0,\bar\alpha)\in\Gamma(D_{\rm red})$ and we find
$X\in\Gamma(\mathcal{H})$ such that $X\sim_{\pi}0$ and
$(X,\pi^*\bar\alpha)\in\Gamma(D\cap\mathcal{K}^\perp)$. If $X=0$, then we
have $\pi^*\bar\alpha=0$ on $\mathcal{H}+\mathcal{V}=TM$, contradicting the fact that
$\bar{ \alpha} $ is a non-zero section of $\mathsf{P}_0^{\rm red}$. Therefore $X$ is
a non-zero vector field lying in $\Gamma(\mathcal{H}\cap\mathcal{V})$ with
${\mathbf{i}}_{X}\omega_M=0$ on $\mathcal{H}\cap\mathcal{V}$. We conclude from this that
the sections of $\mathsf{P}_0^{\rm red}$ pull back exactly to the
$G$-equivariant 
sections $\alpha^{\boldsymbol{\xi}}+\beta\in \Gamma(\mathcal{V})^\circ$ induced by sections
$\boldsymbol{\xi}$ of
$(\mathcal{V}\cap\mathcal{H})\cap(\mathcal{V}\cap\mathcal{H})^{\omega_M}$.

If $\mathcal{V}\cap\mathcal{H}\nsubseteq
  \boldsymbol{\xi}^{\omega_M}$, i.e., $\boldsymbol{\xi}\not\in\Gamma(\mathcal U)$,
 and hence $\alpha^{\boldsymbol{\xi}}\neq0$ on
  $\mathcal{V}\cap\mathcal{H}$, then there is no
  $\beta\in\Gamma(\mathcal{H}^\circ)$ such that
  $(\boldsymbol{\xi},\alpha^{\boldsymbol{\xi}}+\beta)\in
  \Gamma(D\cap\mathcal{K}^\perp)$ and  {\rm (ii)} follows immediately (see also Proposition 
\ref{prop_nonhol_dirac}\rm{(ii)}).
\end{proof}

\begin{definition}\label{def_nonhol_noether}
In the conditions of the preceding proposition, we will call \emph{Nonholonomic Noether Equation} a section
$\alpha^{\boldsymbol{\xi}}$ corresponding to a smooth section
$\boldsymbol{\xi}$ of $\mathcal{V}\cap\mathcal{H}$.
A \emph{$\mathcal{H}$-modified Nonholonomic Noether Equation} is a $1$-form $\alpha'\in \Omega^1(M)$ that can be written $\alpha'=\alpha^{\boldsymbol{\xi}}+\beta$
with a nonholonomic Noether equation $\alpha^{\boldsymbol{\xi}}$ and
$\beta\in\Gamma(\mathcal{H}^\circ)$. A \emph{Descending ($\mathcal{H}$-Modified) Nonholonomic Noether
Equation} is a ($\mathcal{H}$-modified) nonholonomic Noether
equation as in Proposition \ref{prop_73} {\rm (i)}. 
\end{definition}

Note that because of the $\beta$-part of a descending $\mathcal{H}$-modified nonholonomic 
Noether equation, sections of $\mathsf{P}_0^{\rm red}$ don't pull back exactly to sections 
$\alpha^{\boldsymbol{ \xi}}$ associated to sections $\xi^{ \mathcal{H}} $ as
in Theorem \ref{nonholonomic_theorem} (the nonholonomic Noether equations). It is possible that they pull back to one-forms
that coincide only on $\mathcal{H}$ with some  $\alpha^{\boldsymbol{ \xi}}$.

\begin{proposition}
The codistribution spanned by the Noether equations which descend to the quotient $M/G$ is given by 
\begin{equation}\label{d_projection_nonholonomic}
\pi_2(D\cap
(\mathcal{V}\oplus\mathcal{V}^\circ))  
= (\flat(\mathcal{V}\cap \mathcal{H})+\mathcal{H}^\circ)\cap
\mathcal{V}^\circ\end{equation}
where $\flat:TM \rightarrow T ^\ast M $ is associated to $\omega_M$.
\end{proposition}

\begin{proof}
We have seen that a  descending ($\mathcal{H}$-modified) nonholonomic Noether equation $\alpha'$  is a  $G$-invariant
section of $\mathcal{V}^\circ$ such that there exists a $G$-equivariant
section $X$ of $\mathcal{V}\cap\mathcal{H}$ with
$(X,\alpha)\in\Gamma(D\cap\mathcal{K}^\perp)$. So we only have to show
equality \eqref{d_projection_nonholonomic}.
Let $\alpha$ be a section of the left-hand side. Then there exists
$X\in\Gamma(\mathcal{V}\cap\mathcal{H})$ such that $(X,\alpha)\in\Gamma(D\cap
(\mathcal{V}\oplus\mathcal{V}^\circ))$ and hence there exists
$\beta\in\Gamma(\mathcal{H}^\circ)$ such that
$\alpha={\mathbf{i}}_{X}\omega_M+\beta$. Therefore
$\alpha\in\Gamma\left((\flat(\mathcal{V}\cap
  \mathcal{H})+\mathcal{H}^\circ)\cap\mathcal{V}^\circ\right)$.
Conversely, let $\alpha$ be a section of $(\flat(\mathcal{V}\cap
  \mathcal{H})+\mathcal{H}^\circ)\cap\mathcal{V}^\circ$; then
  $\alpha={\mathbf{i}}_{X}\omega_M+\beta\in\Gamma(\mathcal{V}^\circ)$ with $X\in\Gamma(\mathcal{V}\cap
  \mathcal{H})$ and $\beta\in\Gamma(\mathcal{H}^\circ)$. But this means
  that $(X,\alpha)$ is a section of $D\cap
(\mathcal{V}\oplus\mathcal{V}^\circ)$.
\end{proof}

\begin{example} \label{ex:sec_71}
We compute $\boldsymbol{\xi}$ and $\alpha^{\boldsymbol{\xi}}$ for the constrained
particle (see \S\ref{constrpart}). In this example, $Q = \mathbb{R}^3$,
$M:=\{(x,y,z,p_x,p_y,p_z)\mid p_z=y p_x\}\subseteq T^*Q = \mathbb{R}^3\times
\mathbb{R}^3$, and $\mathcal{V}\cap \mathcal{H} = \operatorname{span}\{\partial_x + y
\partial_z\} $.
 If $\xi^1:=(1,0)$ and $\xi^2:=(0,1)$ is the chosen basis of
$\mathfrak{g}=\mathbb{R}^2$, then $(1,0)_M = \partial_x$ and $(0,1)_M =
\partial_z$
so that $\mathfrak{g}^{(x,y,z,p_x,p_y)} := \operatorname{span}\{(1,0) + y (0,1) \}$ is
the fiber of the vector bundle $\mathfrak{g}^ \mathcal{H}$ at the point
$(x,y,z,p_x,p_y) \in M$. Therefore, any section $\xi^ \mathcal{H}$ of $\mathfrak{g}^
\mathcal{H}$ has the form $\xi^ \mathcal{H}(x,y,z,p_x,p_y) = f(x,y,z,p_x,p_y)
\left((1,0) + y (0,1) \right)$, where $f \in C ^{\infty}(M) $. Consequently 
\[
\boldsymbol{ \xi}(x,y,z,p_x,p_y) = \left(\xi^{\mathcal{H}}(x,y,z,p_x,p_y)\right)_M
(x,y,z,p_x,p_y) = f(x,y,z,p_x,p_y) \left( \partial_x + y\partial_z \right).
\]

The components of the momentum map $\mathbf{J}: T ^\ast Q \rightarrow
\mathfrak{g}^\ast$ are $\mathbf{J}^{(1,0)} = p _x$ and $\mathbf{J}^{(0,1)} = p _z$ so
that the restrictions to $M $ of these functions are $\mathbf{J}^{(1,0)} = p _x$ and 
$\mathbf{J}^{(0,1)} = yp _x$. Therefore 
$\mathbf{J}^{\xi^{\mathcal{H}}} (x,y,z,p_x,p_y) = f(x,y,z,p_x,p_y) p _x( 1+y^2)$
and if $X \in \mathfrak{X}(M) $, then $X[f(1,0) + yf (0,1)] = X[f](1,0) + X[yf](0,1)
$ and so 
\[
\mathbf{J}^{ X[f(1,0) + yf (0,1)]} = X[f]p _x + X[yf] y p _x
= p _x(1+y^2) X[f] + y p _x f X[y].
\]
The $1$-form on $M $ which applied to $X $ yields the  right hand
side is $p _x(1+y^2) \mathbf{d}f + yp _x f \mathbf{d}y$ and hence 
\begin{align*}
\alpha^{\boldsymbol{ \xi}}(x,y,z,p_x,p_y) 
&= \mathbf{d} \mathbf{J}^{\xi^{ \mathcal{H}}}  (x,y,z,p_x,p_y)
- p _x(1+y^2) \mathbf{d}f(x,y,z,p_x,p_y) - yp _x f(x,y,z,p_x,p_y) \mathbf{d}y\\
& = f(x,y,z,p_x,p_y) \left((1+y^2) \mathbf{d}p _x + y p _x\mathbf{d}y  \right).
\end{align*}
So the section spanning the codistribution  $\mathsf{P}_0^{\rm red}$ in this example
is $(1+y^2) \mathbf{d} p _x  + yp _x \mathbf{d} y $, as \eqref{ex1_D_red}. It is easy
to see that in this case
$\mathcal{V}\subseteq(\mathcal{V}\cap\mathcal{H})^{\omega_M}$ and hence the
nonholonomic Noether equation descends to the quotient.
\end{example}

\subsection{The reaction-annihilator distribution}\label{sec:reac_annhili}

In this section also, we assume that $\mathcal V\cap\mathcal H$ has constant rank on $M$.
An important problem is to decide when the nonholonomic Noether momentum equation
gives a \emph{constant of motion}
rather than an \emph{equation of motion}. We have to distinguish between two cases: 
\begin{enumerate}
\item[{\rm (i)}] The section $\xi^{\mathcal{H}}$ is constant,
i.e, $\xi^{\mathcal{H}}(p)=\xi$ for all $ p \in  M$, where $\xi\in\mathfrak{g}$. Then
$\boldsymbol{ \xi} (p) = \xi_M(p) $ and so we have $\alpha^{\boldsymbol{\xi}}
=\mathbf{d} \mathbf{J}^\xi$, so $\mathbf{J}^\xi$ is a constant of the motion. We will
see below that sometimes one can find $\eta \in \mathfrak{g}$ such that $\mathbf{J}
^\eta$ is a constant of motion for all solutions of $G $-invariant Hamiltonians, but
$\eta _M$ is not a section of $\mathcal{V}\cap \mathcal{H}$ (see also
\cite{FaRaSa07}).

\item[{\rm (ii)}] The other case is that of \emph{gauge symmetries}, that is,
non-constant sections of $\mathfrak{g}^{\mathcal{H}}$ that yield
constants of motion (see \cite{FaRaSa07}). Note that if
$\xi^{\mathcal{H}}=\sum_{i=1}^kf_i\xi_i$ then it leads to  a constant of
motion if one of the corresponding forms $\alpha^{\boldsymbol{ \xi}} +  \beta$
is exact, that is, we can find $f\in C^\infty(M)$ such that $\mathbf{d}f=
\beta + \sum_{i=1}^k f_i\mathbf{d}J(\xi_i)$. 
However, we do not know of any other
characterization of the section so that the momentum equation gives constants of motion
rather than an equation of motion.
\end{enumerate} 

In the reduction method for nonholonomic systems, the first step is to compute the
horizontal annihilator $\mathcal{U}$ of $\mathcal{V}$, that is, the subbundle $\mathcal{U}=(\mathcal{V}\cap 
\mathcal{H})^{\omega_M}\cap \mathcal{H} \subseteq TM
\subseteq T T ^\ast Q$ (see \eqref{def_cal_U}). 
We have seen in Proposition \ref{prop_nonhol_dirac}(ii) that any section of
$\mathcal{U}$ corresponds to a section of $D\cap \mathcal{K}^\perp $:
for each $X\in \Gamma(\mathcal{U})$ there exists
$\alpha\in\Gamma(\mathcal{V}^\circ)$ such that $(X,\alpha)\in\Gamma(D)$ and hence
$\alpha-{\mathbf{i}}_{X}\omega_M\in\Gamma(\mathcal{H}^\circ)$. So the method of finding a
section $\alpha \in \Gamma( \mathcal{V}^\circ)$ associated to $X \in \Gamma(
\mathcal{U}) $ is the same as determining $\beta\in \Gamma(\mathcal{H}^ \circ ) $
such that $\mathbf{i}_ X \omega _M + \beta = : \alpha \in \Gamma( \mathcal{V}^ \circ
) $. As we have seen in \S\ref{ex_vertical_disk}, case 3, sometimes not the whole of
$\mathcal{H}^\circ$ is needed in this construction. This is why we introduce the new
codistribution $\mathcal{R}$ on $ M $ whose fiber at $p\in M$ equals
\begin{equation}\label{defR}
\mathcal{R}(p)=\{\beta(p)\mid \beta\in\Gamma(\mathcal{H}^\circ)\text{ and there is
some }
X\in\Gamma(\mathcal{U}) \text{ such that }
\beta+{\mathbf{i}}_{X}\omega_M\in\Gamma(\mathcal{V}^\circ)\}\subseteq\mathcal{H}^\circ.
\end{equation} 
In general, $ \mathcal{R} $ is strictly included in $\mathcal{H}^\circ$.

If $h\in C^\infty(M)^G$ is an admissible function, then there exists
$X_h\in\Gamma(\mathcal{H})$ such that $(X_h,\mathbf{d} h)\in \Gamma(D)$. Recall that
$X _h$ is unique since $\mathsf{G}_0 = \{0\} $. In addition, since $\mathbf{d}
h\in\Gamma(\mathcal{V}^\circ)$, we have $(X_h,\mathbf{d} h)\in\Gamma(D\cap K^\perp)$
and hence $X_h\in\Gamma(\mathcal{U})$. Thus, there exists
$\beta\in\Gamma(\mathcal{H}^\circ)$ such
that $\mathbf{d} h={\mathbf{i}}_{X_h}\omega_M+\beta$. This is exactly the Hamilton equation
for
the given
nonholonomic system (with the Hamiltonian $h$) and  we have
$\beta\in\Gamma(\mathcal{R})$, often interpreted as the \emph{reaction force}.
In fact $\mathcal{R}^\circ\subset  TM$ is the analogue of the
\emph{reaction-annihilator distribution} of \cite{FaRaSa07}.

\begin{proposition}\label{RplusV}
We have
\begin{equation*}
\flat(\mathcal{U})\oplus\mathcal{R}=\mathcal{V}^\circ+\mathcal{R},
\end{equation*}
where $\flat: TM \rightarrow T ^\ast M $ corresponds to $\omega_M$.
\end{proposition}
\begin{proof}
The sum on the left hand side is direct since if
$X\neq 0\in\Gamma(\mathcal{U})$, then $X\in\Gamma( \mathcal{H})$ and hence
${\mathbf{i}}_{X}\omega_M\not\in\Gamma(\mathcal{H}^\circ)$ because
$\omega_M\arrowvert_{\mathcal{H}\times\mathcal{H}}$ is nondegenerate. Thus
${\mathbf{i}}_{X}\omega_M \notin \Gamma(\mathcal{R})\subseteq
\Gamma(\mathcal{H}^\circ)$.
Second, recall that $\mathsf{P}_1=T^*M$ (see Proposition \ref{prop_nonhol_dirac}(i))
so for all $\alpha\in\Gamma(\mathcal{V}^\circ)$ we find $X\in\Gamma(\mathcal{H})$
(actually $X\in\Gamma(\mathcal{U})$) such that $(X,\alpha)\in\Gamma(D)$. Thus,
$\pi_2(D\cap\mathcal{K}^\perp)=\mathcal{V}^\circ$. 

Now we are ready to prove the formula in the statement. 
If $X\in\Gamma(\mathcal{U})$, the considerations in 
\S\ref{link-nonhol-dirac} show that there exists 
$\beta\in\Gamma(\mathcal{H}^\circ)$ such
that ${\mathbf{i}}_{X}\omega_M+\beta\in\Gamma(\mathcal{V}^\circ)$.
The definition \eqref{defR} of $\mathcal{R}$ yields directly that
$\beta\in\Gamma(\mathcal{R})$. This shows
$\flat(\mathcal{U})\oplus\mathcal{R}\subseteq\mathcal{V}^\circ+\mathcal{R}$.
For the other inclusion, choose $\alpha\in\Gamma(\mathcal{V}^\circ)$ and
$X\in\Gamma(\mathcal{U})$ such that
$(X,\alpha)\in\Gamma(D\cap\mathcal{K}^\perp)$. Then the definition
of $D$ yields $\beta:=\alpha-{\mathbf{i}}_{X}\omega_M\in\Gamma(\mathcal{H}^\circ)$ and
again, using \eqref{defR}, we conclude that  $\beta\in\Gamma(\mathcal{R})$.
\end{proof}
The last lemma leads directly to the equality
\begin{equation*}
\mathcal{U}^{\omega_M}\cap\mathcal{R}^\circ=\mathcal{V}\cap\mathcal{R}^\circ.
\end{equation*}
Note that
$\mathcal{U}^{\omega_M}=((\mathcal{H}\cap\mathcal{V})^{\omega_M}\cap\mathcal{H})^{\omega_M}=(\mathcal{H}\cap\mathcal{V})+\mathcal{H}^{\omega_M}$
since the kernel of $\omega_M$ lies in $\mathcal{H}^{\omega_M}$. 

Now we are able to state the main theorem of this subsection which is the
Hamiltonian analogue of the main statement of \cite{FaRaSa07}.
\begin{theorem}
Let $\xi\in\mathfrak{g}$. Then the function $\mathbf{J}^\xi$ is a constant of motion
for
every $G$-invariant Hamiltonian $h$ if and only if
$\xi_M\in\Gamma(\mathcal{V}\cap\mathcal{R}^\circ)$.
\end{theorem}

\begin{proof}
Choose $\xi\in \mathfrak{g}$ such that
$\xi_M\in\Gamma(\mathcal{V}\cap\mathcal{R}^\circ)$. We have seen in the
preceding section that ${\mathbf{i}}_{\xi_M}\omega_M=\mathbf{d} \mathbf{J}^\xi$. For an
arbitrary
$X\in\Gamma(\mathcal{U})$ choose $\beta\in\Gamma(\mathcal{R})$ with
${\mathbf{i}}_{X}\omega_M+\beta=:\alpha\in\Gamma(\mathcal{V}^\circ)$ and get 
\begin{align*}
\mathbf{d} \mathbf{J}^\xi(X)=\omega_M(\xi_M,X)=\beta(\xi_M)-\alpha(\xi_M)=0.
\end{align*}
This yields the statement since for all $G$-invariant Hamiltonian $h$ the (unique)
solution $X_h$ of the implicit Hamiltonian system $(X,\mathbf{d}
h)\in\Gamma(D)$ is a section of $\mathcal{U}$ (with $\alpha=\mathbf{d} h$ the
corresponding section of $\mathcal{V}^\circ$ and
$\beta=\mathbf{d} h-{\mathbf{i}}_{X_h}\omega_M$). 
For the converse implication, choose $\xi\in\mathfrak{g}$ such that $\mathbf{J}^\xi$ is
a
constant of the motion for the solution curves of every $G$-invariant
Hamiltonian. Note that since $\mathcal{V}$ is an involutive subbundle of
$TM$, the exterior derivatives of all $G$-invariant functions span pointwise
$\mathcal{V}^\circ$ and hence the corresponding solutions span
$\mathcal{U}$. This yields $\mathbf{d} \mathbf{J}^\xi=0$ on $\mathcal{U}$. If we
choose
$\beta\in\Gamma(\mathcal{R})$, there exists
$X\in\Gamma(\mathcal{U})$ such that
${\mathbf{i}}_{X}\omega_M+\beta\in\Gamma(\mathcal{V}^\circ)$. Hence we get 
\begin{align*}
0=({\mathbf{i}}_{X}\omega_M+\beta)(\xi_M)=\omega_M(X,\xi_M)+\beta(\xi_M)=-\mathbf{d} \mathbf{J}^\xi(X)+\beta(\xi_M)=0+\beta(\xi_M)
\end{align*}
and therefore $\xi_M\in\Gamma(\mathcal{R}^\circ\cap\mathcal{V})$.
\end{proof}

\begin{corollary}
Assume that $\mathcal{H}+ \mathcal{V}$ has constant rank on $M $ and choose $\xi\in \mathfrak{g}$.
If $\mathbf{d} \mathbf{J}^\xi=0$ on
$\mathcal{V}\cap\mathcal{H}$ there exist $\beta\in\Gamma(\mathcal{H}^\circ)$
and $\boldsymbol{\eta}\in\Gamma(\mathcal{V}\cap\mathcal{H})$ such that
$\alpha^{\boldsymbol{\eta}}=\mathbf{d}
\mathbf{J}^\xi+\beta$. 
\end{corollary}

\begin{proof}
Since $\mathsf{P}_1=T^*M$ (see Proposition \ref{prop_nonhol_dirac}{\rm (i)})
there exists $X\in\Gamma(\mathcal{H})$ such that $(X,\mathbf{d}
\mathbf{J}^\xi)\in\Gamma(D)$.
Hence we have $\mathbf{d} \mathbf{J}^\xi={\mathbf{i}}_{X}\omega_M+\beta'$ with
$\beta'\in\Gamma(\mathcal{H}^\circ)$ and since $\mathbf{d} \mathbf{J}^\xi=0$
on $\mathcal{U}$ we get $X\in\mathcal{U}^{\omega_M}=
(\mathcal{H}\cap\mathcal{V})+\mathcal{H}^{\omega_M}$. Write $X=V+Y$ with
$V\in\Gamma(\mathcal{H}\cap\mathcal{V})$ and $Y\in
\Gamma(\mathcal{H}^{\omega_M})$. Since $X$ and $V$ are sections of
$\mathcal{H}$, then so is $Y$. But since
$\mathcal{H}\cap\mathcal{H}^{\omega_M}=\{0\}$, this yields $Y=0$ and hence
$X\in\Gamma(\mathcal{H}\cap\mathcal{V})$. We find
$\eta^{\mathcal{H}}\in\Gamma(\mathfrak{g}^{\mathcal{H}})$
such that the corresponding section
$\boldsymbol{\eta}\in\Gamma(\mathcal{V}\cap\mathcal{H})$ is equal to $X$ and
therefore
$(\eta,\mathbf{d} \mathbf{J}^\xi)\in\Gamma(D)$. We get
$\alpha^{\boldsymbol{\eta}}=\mathbf{d}
\mathbf{J}^\xi+\beta$ with $\beta\in\Gamma(\mathcal{H}^\circ)$, a nonholonomic
Noether equation corresponding to the section
$\eta^{\mathcal{H}}\in\Gamma(\mathfrak{g}^{\mathcal{H}})$.
\end{proof}

\subsection{Optimal momentum map for nonholonomic mechanical systems}

\textit{In this and the next subsection we assume that $\mathcal{H}+ \mathcal{V}$ has
constant rank on $M $.} Recall from Remark \ref{rem_dim_HplusV} that this
implies that $\mathcal{V}\cap \mathcal{H}$ and $\mathcal{U}$ also have constant rank on $M $. 

\medskip

We show in this subsection that, under certain integrability assumptions, it is
possible to restrict the system to ``level sets'' given by the nonholonomic momentum
equations and then perform reduction.  

Consider the distribution where all $\alpha^{\boldsymbol{\xi}},\alpha'$ as in Definition \ref{def_nonhol_noether}
 vanish, namely
\begin{align*}
\mathcal{D}_G:= \left[\pi_2(D\cap
(\mathcal{V}\oplus\mathcal{V}^\circ)) \right]^ \circ 
\overset{\eqref{d_projection_nonholonomic}}=\left[(\flat(\mathcal{V}\cap \mathcal{H})+\mathcal{H}^\circ)\cap
\mathcal{V}^\circ\right]^ \circ 
=\left((\mathcal{V}\cap \mathcal{H})^{\omega_M}\cap
\mathcal{H}\right)+\mathcal{V}=\mathcal{U}+\mathcal{V},
\end{align*}
where $\mathcal{U}:= (\mathcal{V}\cap \mathcal{H})^{\omega_M}\cap \mathcal{H}
\subseteq TM \subseteq T T ^\ast Q$ is the horizontal annihilator  of $\mathcal{V}$
(see \eqref{def_cal_U}). Note that 
$\mathcal{D}_G= (\mathcal{V}\cap \mathcal{H})^{\omega_{\mathcal{H}}}+\mathcal{V}$.
Since $\mathcal{U}\subseteq \mathcal{H} $ it follows easily that 
$\mathcal{D}_G\cap\mathcal{H}= \mathcal{U}+(\mathcal{V}\cap
\mathcal{H})$. If $\mathcal{D}_G$ is integrable, its leaves are the level sets of
the constants of motion and equations of motion  given by the Nonholonomic Noether
Theorem \ref{nonholonomic_theorem} for sections $\boldsymbol{\xi}$ of
$(\mathcal{V}\cap\mathcal{H})\cap(\mathcal{V}\cap\mathcal{H})^{\omega_M}$: 
the fiber at  $m\in M$ of the distribution
$\pi_2(D\cap(\mathcal{V}\oplus\mathcal{V}^\circ))$  equals
\[
\{\alpha(m)\mid \alpha \in
\Gamma(\mathcal{V}^\circ) \text{ and there exists } X \in \Gamma(\mathcal{V})
\text{ such that } (X,\alpha)\in \Gamma(D)\}.
\]
Note that if this distribution is spanned by closed $1$-forms, hence locally exact
$1$-forms,  then it can be written as 
\[
\{(\mathbf{d} f)(m)\mid f\in
C^\infty(M)^G \text{ and there exists } X_f \in \Gamma(\mathcal{V})
\text{ such that } (X_f,\mathbf{d} f)\in \Gamma(D)\}.
\]
For every $m \in M $, we have
\begin{align*}
\dim \left[(\mathcal{V}(m)\cap
\mathcal{H}(m))^{\omega_{\mathcal{H}}}\cap(\mathcal{V}(m)\cap
\mathcal{H}(m))\right] &= 
\dim \left[\mathcal{U}(m)\cap(\mathcal{V}(m)\cap
\mathcal{H}(m))\right]
= \dim\left[\mathcal{U}(m)\cap\mathcal{V}(m) \right]\\
&=\dim \mathcal{U}(m)+\dim\mathcal{V}(m)-
\dim \mathcal{D}_G(m).
\end{align*}
Recall that $\mathcal{U}$ and $\mathcal{H} \cap \mathcal{V}$ are vector subbundles of $TM $
by hypothesis. If, in addition, $\mathcal{D}_G$ is integrable, then its fibers $\mathcal{D}_G(m)$
have constant dimension along the leaves of the generalized foliation determined by
$\mathcal{D}_G $ and so the computation above shows that the fibers of 
$(\mathcal{V}\cap \mathcal{H})^{\omega_{\mathcal{H}}}\cap(\mathcal{V}\cap
\mathcal{H})$ along a leaf of $\mathcal{D}_G$ are constant. Thus, the same is true for
the fibers of $\mathcal{D}_G \cap \mathcal{H} = \mathcal{U} + ( \mathcal{V} \cap
\mathcal{H}) $ since $\mathcal{U} \cap \mathcal{V} \cap \mathcal{H} =
(\mathcal{V}\cap \mathcal{H})^{\omega_{\mathcal{H}}}\cap(\mathcal{V}\cap
\mathcal{H})$.  We shall use
this fact in the next subsection where we  describe the induced Dirac structure on a
leaf. 

In order to restrict the system to the leaves of the distribution
$\mathcal{D}_G$ and then perform reduction, we have to show several
statements, the analogues of those needed for the Dirac optimal reduction. Since $\Phi_g ^\ast \omega_M = \omega_M $ for all $g \in G $, the proof of the following proposition follows easily.
 
\begin{proposition}
The distribution $(\mathcal{V}\cap \mathcal{H})^{\omega_M}$ is $G$-invariant in
the sense that 
\[
\Phi_g^*\left((\mathcal{V}\cap \mathcal{H})^{\omega_M}\right)
=(\mathcal{V}\cap \mathcal{H})^{\omega_M}
\] 
for all $g\in G$. Since $\mathcal{V}$ and
$\mathcal{H}$ are also $G$-invariant, it follows that
the distribution $\mathcal{D}_G=\left((\mathcal{V}\cap
\mathcal{H})^{\omega_M}\cap\mathcal{H}\right)+\mathcal{V}$ is $G$-invariant.
\end{proposition}

If $\mathcal{D}_G$ is integrable, define  the 
\textit{nonholonomic optimal momentum map}
\[
\mathcal{J}:M\to M/\mathcal{D}_G.
\]
\begin{lemma}
If $m$ and $m'$ are in the same leaf of $\mathcal{D}_G$,  then  $\Phi_g(m)$ and $\Phi_g(m')$ are in the same
leaf of $\mathcal{D}_G$ for all $g \in G$. Hence there is a well defined action of
$G$ on $M/\mathcal{D}_G$:
\begin{align*}
\bar{\Phi}:G\times M/\mathcal{D}_G&\to M/\mathcal{D}_G\\
\bar{\Phi}_g(\mathcal{J}(m))&=\mathcal{J}(g\cdot m)
\end{align*}
For all $\rho\in M/\mathcal{D}_G$, the isotropy subgroup of $\rho$ contains
$G^\circ$ (the connected component of the identity in $G$). 
\end{lemma}

\begin{proof}
Let $g \in G $, $m,m' \in M$ be in the same leaf of $\mathcal{D}_G$, i.e.,  there exists without loss
of generality one vector field
$X\in \Gamma(\mathcal{D}_G)$ with flow $F^{X}$ such that  $F^{X}_t(m)=m'$ for some
$t>0$. We have for all $s \in [0,t]$:
\begin{align*}
\frac{d}{ds}\left(\Phi_g\circ F^{X}_s\right)(m)&=T_{F^{X}_s(m)}\Phi_g
\left(X(F^{X}_s(m))\right)
=(\Phi_{g^{-1}}^*X)(g\cdot F^{X}_s(m))\in \mathcal{D}_G(g\cdot F^{X}_s(m)).
\end{align*} 
Hence the curve $c(s)=\left(\Phi_g\circ F^{X}_s\right)(m)$ connecting
$c(0)=\Phi_g(m)$ and $c(t)=\Phi_g(m')$ lies entirely in the leaf of $\mathcal{D}_G$
through the point $\Phi_g(m)$ and the assertion follows.

The Lie group $G^\circ$ is generated as a group by the exponential of an open
neighborhood of $0\in\mathfrak{g}$. Thus, we can assume without loss of generality,
that for any $g\in G^\circ$ and $m\in M$, there exists some $\xi\in \mathfrak{g}$
such that the curve $\gamma:[0,t]\to M$, $\gamma(s)=\Phi_{\exp(s\xi)}(m)$, has
endpoints $m$ and $g\cdot m$ (in reality, the points $m$ and $g\cdot m$ can be joined
with finitely many  such curves). 
 For all $s\in [0,t]$, we have
$\dot{\gamma}(s)=\xi_M(\gamma(s))\in \mathcal{D}_G(\gamma(s))$ and, arguing as above,
we conclude that the whole curve $\gamma([0,t])$ lies in the leaf of $\mathcal{D}_G$
through $m$. Hence, if $\rho= \mathcal{J}(m)$, the equality $\Phi_g(\mathcal{J}(m))=
\mathcal{J}(g\cdot m) = \mathcal{J}(m)$ proves the statement.
\end{proof}

\begin{remark}\label{rem}
The last statement shows that for all $\rho\in M/\mathcal{D}_G$, the isotropy
subgroup $G_\rho$ is the union of connected components of $G$ and is therefore 
closed
in $G$. This implies that the Lie group $G_\rho$ acts properly on the leaf
$\mathcal{J}^{-1}(\rho)$. It is obvious that this action is also free. 
In the Optimal reduction results in \cite{OrRa04}, \cite{JoRa10b},
the induced action on the leaves of the distribution is not necessarily proper.
The reason why  the action is here always proper is the inclusion $\mathcal{V} \subset \mathcal{D}_G$.
\end{remark}

\begin{remark}
Note that if the nonholonomic system satisfies $\mathcal{H}\oplus\mathcal{V}=TM$, then the bundle $\mathcal{U}$ is given by
$\mathcal{U}=\{0\}^{\omega_M}\cap\mathcal{H}=\mathcal{H}$ and hence
$\mathcal{D}_G=\mathcal{U}+\mathcal{V}=TM$ is trivially integrable with the
connected components of $M$ as integral leaves.
Hence, if $M $ is connected,  the method of reduction presented in the next subsection leads to the
same reduced Dirac manifold as the Dirac reduction method of  
\S\ref{link-nonhol-dirac}.
\end{remark}

\subsection{Optimal reduction for nonholonomic systems}
\label{sec:optimal_reduction_nonholonomic}
Assume as above that $\mathcal D_G$ is integrable and choose
$\rho\in M/\mathcal D_G$. Since  the isotropy subgroup $G_\rho$
contains $G^\circ$, the distribution $\mathcal{V}_\rho$ spanned by the
fundamental vector fields of 
the action of $G$ on $\mathcal{J}^{-1}(\rho)$ is
$\mathcal{V}_\rho=\mathcal{V}|_{ \mathcal{J} ^{-1} 
(\rho) }$. Since $\mathcal U\subseteq \mathcal D_G$, 
$\mathcal U$ also restricts to $\mathcal U_\rho$ on  the manifold $\mathcal{J}^{-1}(\rho)$. 
Let $\mathcal H_\rho$ be the intersection of  $\mathcal{H}:=TM\cap (T\pi_{T^*Q})^{-1}(\mathcal{D})$
with $T\mathcal{J}^{-1}(\rho)$.

Since the distribution $(\mathcal{H}\cap (\mathcal{V}\cap \mathcal{H})^{\omega_{M}}+\mathcal{V})\cap\mathcal{H}
=\mathcal{U}+(\mathcal{V}\cap \mathcal{H}) \subseteq \mathcal{D}_G$ is constant dimensional  on the leaves of  
$\mathcal{D}_G$, the Dirac structure on a leaf $\mathcal{J}^{-1}(\rho)$ of 
$\mathcal{D}_G$ is given by
\begin{align*}
D_{\mathcal{J}^{-1}(\rho)}(m)=\left\{(X(m),\alpha_m)\in
T\mathcal{J}^{-1}(\rho)\oplus T^*\mathcal{J}^{-1}(\rho)\left|
\begin{array}{c}
 X\in
\Gamma\left(\mathcal{U}_\rho+(\mathcal{V}_\rho\cap
  \mathcal{H}_\rho)\right),\\
\alpha-{\mathbf{i}}_{X}\omega_{\mathcal{J}^{-1}(\rho)}\in
\Gamma\left(\left(\mathcal{U}_\rho+(\mathcal{V}_\rho\cap
  \mathcal{H}_\rho)\right)^{\circ}\right)
\end{array}
\right.\right\}
\end{align*} 
for all $m\in \mathcal{J}^{-1}(\rho)$ (see \cite{BlvdS01}); here $
i_ \rho : \mathcal{J} ^{-1} (\rho) \hookrightarrow M $ is the inclusion and $\omega_{\mathcal{J}^{-1}(\rho)}:=i_\rho^*\omega_M$. 

\begin{lemma}\label{lem:dimDcapKperp}
Let $\mathcal{K}_\rho=\mathcal{V}_\rho\oplus\{0\}$ and
$\mathcal{K}_\rho^\perp=T\mathcal{J}^{-1}(\rho)\oplus\mathcal{V}_\rho^\circ$ as in \S\ref{red} . Then
$D_{\mathcal{J}^{-1}(\rho)}\cap\mathcal{K}_\rho^\perp$ is a vector bundle over
$\mathcal{J}^{-1}(\rho)$. 
\end{lemma}

\begin{proof}
Since $\mathcal{H}+\mathcal{V}$ has constant rank on the $n $-dimensional manifold $M$, recall from Remark \ref{rem_dim_HplusV} that $D\cap\mathcal{K}^\perp$ is a vector bundle on $M$. We denote $r = \operatorname{rank} \mathcal{H}$, $n-r = \operatorname{rank} \mathcal{H}^ \circ $, $l = \operatorname{rank}\mathcal{V}^ \circ \cap \mathcal{H}^ \circ $, $u = \operatorname{rank} \mathcal{U}$, and $s = \operatorname{rank} ( \mathcal{U} + ( \mathcal{V}\cap  \mathcal{H}))|_{ \mathcal{J}^{-1}(\rho)}$. Let $m\in\mathcal{J}^{-1}(\rho)$. As in Remark \ref{rem_dim_HplusV}, choose local
basis fields $H_1,\dots,H_r$ for $\mathcal{H}$ and local basis $1$-forms
$\beta_1,\dots,\beta_{n-r}$ for $\mathcal{H}^\circ$ defined on a neighborhood $U$ of $m$ in $M$. Assume that  $H_1,\dots,H_u$ are local basis fields for
$\mathcal{U}$, $H_1,\dots,H_s$, with $u\leq s\leq r$, are basis fields
for $\mathcal{U}+(\mathcal{V}\cap\mathcal{H})$ on $\mathcal{J}^{-1}(\rho)\cap U$, and $\beta _1, \ldots , \beta_l$ a basis of $\mathcal{V}^ \circ \cap \mathcal{H}^ \circ = ( \mathcal{V}+ \mathcal{H}) ^ \circ$. Note that
the $1$-forms $\beta_1,\dots,\beta_l$ vanish on
$\mathcal{U}+\mathcal{V}\subseteq\mathcal{H}+\mathcal{V}$ and that
$\beta_{l+1},\dots,\beta_{n-r}$ don't vanish on $\mathcal{U}+\mathcal{V}$
(otherwise we would have
$\beta_j\in\Gamma(\mathcal{U}^\circ\cap\mathcal{V}^\circ\cap\mathcal{H}^\circ)=\Gamma(\mathcal{V}^\circ\cap\mathcal{H}^\circ)$
for $j=l+1,\dots,n-r$,
in contradiction to the choice of $\beta_1,\dots,\beta_{n-r}$). The Dirac
structure $D_{\mathcal{J}^{-1}(\rho)}$  is then given on $U \cap \mathcal{J}^{-1}(\rho)$ by 
(see  \cite{Blankenstein00})
\begin{equation}\label{formula2Durb}
\operatorname{span}\left\{(\tilde H_1,i_\rho^*{\mathbf{i}}_{ H_1}\omega_M), \dots,
    (\tilde H_s,i_\rho^*{\mathbf{i}}_{H_s}\omega_M),(0, i_\rho^*\beta_{l+1}), \dots, (0,
    i_\rho^*\beta_{n-r})\right\},
\end{equation}
where $\tilde{H}_1,\dots,\tilde{H}_s$ are vector fields on $U\cap\mathcal{J}^{-1}(\rho)$ such
that $\tilde{H}_i\sim_{i_\rho} H_i$ for $i=1,\dots,s$. Note that
$i_\rho^*{\mathbf{i}}_{H_i}\omega_M={\mathbf{i}}_{\tilde
  H_i}\omega_{\mathcal{J}^{-1}(\rho)}$ for $i=1,\dots,s$.

If $(X,\alpha)\in D_{\mathcal{J}^{-1}(\rho)}\cap (T\mathcal{J}^{-1}(\rho)\oplus
\mathcal{V}^\circ_\rho)$ then 
$X\in \Gamma(\mathcal{U}_\rho+(\mathcal{V}_\rho\cap\mathcal{H}_\rho))$, 
$\alpha\in \Gamma(\mathcal{V}^\circ_\rho )$, and ${\mathbf{i}}_{X}\omega_{\mathcal{J}^{-1}(\rho)}-\alpha\in
\Gamma((\mathcal{U}_\rho+(\mathcal{V}_\rho\cap\mathcal{H}_\rho))^\circ)$. This is only possible if ${\mathbf{i}}_{X}\omega_{\mathcal{J}^{-1}(\rho)}=0$ on
\begin{equation*}
\mathcal{V}_\rho\cap\left(\mathcal{U}_\rho+(\mathcal{V}_\rho\cap\mathcal{H}_\rho)\right)=(\mathcal{V}_\rho\cap\mathcal{U}_\rho)+(\mathcal{V}_\rho\cap\mathcal{H}_\rho)=\mathcal{L}_\rho \quad
\text{where} \quad
\mathcal{L}:=\left[(\mathcal{V}\cap\mathcal{H})\cap(\mathcal{V} \cap\mathcal{H})^{\omega_M} \right]+(\mathcal{V}\cap\mathcal{H}).
\end{equation*}
We have two different cases. First, if $X\in\Gamma(\mathcal{U}_\rho)$ then  for all $m\in\mathcal{J}^{-1}(\rho)$ and $V(m)\in \mathcal{L}(m)$ we have necessarily
\[
({\mathbf{i}}_{X}{\omega_{\mathcal{J}^{-1}(\rho)}})(m)(V(m))
=(i_\rho^*\omega_M)(m)(X(m),V(m))
={\omega_M}(i_\rho(m))(T_mi_\rho X(m),T_mi_\rho V(m))
=0
\]
where we have used  
\[
T_mi_\rho X(m)\in\left(\mathcal{H}\cap(\mathcal{V}
\cap \mathcal{H})^{\omega_M}\right)(i_\rho(m)), \quad T_mi_\rho V(m)\in\mathcal{L}(i_\rho(m))
\] 
and the definition of $\mathcal{L}$. Hence, for all $X\in \Gamma( 
\mathcal{U}_\rho)$ we have ${\mathbf{i}}_{X}\omega_{\mathcal{J}^{-1}(\rho)}\arrowvert_{\mathcal{L}_\rho}=0$ and hence we
find $\alpha\in\Gamma(\mathcal{V}_\rho^\circ)$ such that
$(X,\alpha)\in\Gamma\left(D_{\mathcal{J}^{-1}(\rho)}\right)$. Second, 
for a section $X$ of $\mathcal{V}\cap\mathcal{H}$ that doesn't take values in
$\mathcal{U}_\rho$, the $1$-form  ${\mathbf{i}}_{X}\omega_{\mathcal{J}^{-1}(\rho)}$ doesn't vanish on $\mathcal{V}_\rho\cap\mathcal{H}_\rho$ and thus neither on $\mathcal{L}_\rho$. Consequently,
 the
sections of $D_{\mathcal{J}^{-1}(\rho)}\cap\mathcal{K}_\rho^\perp$ have as first component a section of
$\mathcal{U}_\rho$. Since for $i=l+1,\dots,n-r$ we have
$i_\rho^*\beta_i\not\in\Gamma(\mathcal{V}_\rho^\circ)$, we get 
\[D_{\mathcal{J}^{-1}(\rho)}\cap\mathcal{K}_\rho^\perp=\operatorname{span}\left\{\left(\tilde H_1,{\mathbf{i}}_{\tilde H_1}\omega_{\mathcal{J}^{-1}(\rho)}+\sum_{i=l+1}^{n-r}a_1^ii_\rho^*\beta_i\right), \dots,
    \left(\tilde H_u,{\mathbf{i}}_{\tilde
        H_u}\omega_{\mathcal{J}^{-1}(\rho)}+\sum_{i=l+1}^{n-r}a_u^ii_\rho^*\beta_i\right)\right\},\]
where the functions $a_j^i$ are chosen such that ${\mathbf{i}}_{\tilde
        H_j}\omega_{\mathcal{J}^{-1}(\rho)}+\sum_{i=l+1}^{n-r}a_j^ii_\rho^*\beta_i$ are sections
      of $\mathcal{V}^\circ_\rho$ for $j=1,\dots,u$ and
      $i=l+1,\dots,n-r$. Since the vector fields
      $\tilde{H}_1,\dots,\tilde{H}_u$ are linearly independent, we have found
      basis fields for $D_{\mathcal{J}^{-1}(\rho)}$ on $U$.
\end{proof}

Hence, the reduced Dirac structure $D_\rho$ on $ \mathcal{J} ^{-1}
(\rho)/{G_\rho}$ is given, 
according to the general considerations in 
\S\ref{red} (or \cite{BuCaGu07}) by                      
\begin{equation}\label{drec}
D_\rho
=\left.\frac{\left[D_{\mathcal{J}^{-1}(\rho)}\cap \left(T\mathcal{J}^{-1}(\rho)\oplus
\mathcal{V}^\circ_\rho\right)\right]+(\mathcal{V}_\rho\oplus\{0\})}
{\mathcal{V}_\rho\oplus\{0\} }\right/G_\rho
\end{equation}
The next theorem gives an easier description of this reduced Dirac structure.

We write in the following 
$i_\rho:\mathcal{J}^{-
1}(\rho)\hookrightarrow M$ for the inclusion and
$\pi_\rho:\mathcal{J}^{-1}(\rho)\to M_\rho$ for the projection.

\begin{theorem}[Nonholonomic optimal point reduction by Dirac actions]
As\-su\-me that $\mathcal{H}+ \mathcal{V}$ has
constant rank on $M $ and that the Lie group $G$ acts freely and properly on $M$ by Dirac
actions. If 
$\mathcal{D}_G=\mathcal U+\mathcal V$
 is an integrable subbundle of $TM$, then for any $\rho \in
M/\mathcal{D}_G$ we have 
the following results.
\begin{enumerate}
\item[{\rm (i)}] The orbit space $M_\rho=\mathcal{J}^{-1}(\rho)/{G_\rho}$ is a
smooth regular Dirac 
quotient manifold whose Dirac structure $D_\rho$ is given by the graph of a
\emph{nondegenerate} 
(not necessarily closed) $2$-form $\omega_\rho$. 

\item[{\rm (ii)}] Let $h \in C^\infty(M)^G$ be an admissible and $G$-invariant
Hamiltonian and $X_h$ 
the (unique) solution of the implicit Hamiltonian system $(X_h,\mathbf{d}h)\in
  \Gamma(D)$. Then $X_h\in  \Gamma(\mathcal{U})$ and we have
$(X_h|_{ \mathcal{J} ^{-1}(\rho)},i_\rho^*\mathbf{d} h)\in \Gamma
\left(D_{\mathcal{J}^{-
1}(\rho)}\right)$.
\item[{\rm (iii)}] The flow $F_t$ of $X_h$ leaves $\mathcal{J}^{-1}(\rho)$
invariant, commutes with the $G$-action, and therefore induces a flow
$F_t^\rho$ on $M_\rho$ uniquely determined by the relation $\pi_\rho\circ F_t\circ i_\rho=F_t^\rho\circ \pi_\rho$.
\item[{\rm (iv)}] The flow $F_t^\rho$ is the flow of a vector field 
$X_{h_\rho}$ in $\mathfrak{X}(M_\rho)$ that is the solution of the Hamiltonian
system 
${\mathbf{i}}_{X_{h_\rho}}\omega_\rho=\mathbf{d} h_\rho$, where the function $h_\rho\in
C^\infty(M_\rho)$ is 
given by the equality
  $h_\rho\circ\pi_\rho=h\circ i_\rho$. 
\end{enumerate}
\end{theorem}

\begin{proof}
According to Remark \ref{rem}, the $G_\rho$-action on $\mathcal{J}^{-1}(\rho)$
is free and proper. Thus, the quotient $\mathcal{J}^{-1}(\rho)/{G_\rho}$ is a
regular quotient manifold and the projection
$\pi_\rho:\mathcal{J}^{-1}(\rho)\to M_\rho$ is a smooth surjective
submersion. We denote 
from now on by 
$\omega_{\mathcal{J}^{-1}(\rho)}:=i_\rho^*\omega_M$ the
pull back of $\omega_M$  to $\mathcal{J}^{-1}(\rho)$.

(i) With Lemma \ref{lem:dimDcapKperp}, we get
\begin{align}\label{durbk}
\frac{\left(D_{\mathcal{J}^{-1}(\rho)}\cap \mathcal{K}_\rho^\bot\right) +
\mathcal{K}_\rho}{\mathcal{K}_\rho}
=&\left\{(\widehat{X},\alpha)\in
\Gamma\left( \left(T\mathcal{J}^{-1}(\rho)/\mathcal{V}_\rho\right) \oplus T^*\mathcal{J}^{-1}(\rho) \right)\left| 
\begin{array}{c} \widehat{X}=X(\operatorname{mod}\mathcal{V}_\rho)  \text{ with }\\ X\in
\Gamma(\mathcal{U}_\rho), 
\alpha\in
\Gamma(\mathcal{V}^\circ_\rho),\\
\text{ and } \alpha-{\mathbf{i}}_{X}\omega_{\mathcal{J}^{-1}(\rho)}\\\in
\Gamma\left((\mathcal{U}_\rho+(\mathcal{V}_\rho\cap\mathcal{H}_\rho))^\circ\right)
\end{array}\right.\right\}.
\end{align}
The $G _\rho$-quotient of this bundle defines the reduced Dirac structure $D_\rho$ on $M _\rho$.

Note that the fibers 
\[
\left(\mathcal{U}_\rho/\mathcal{V}_\rho\right)(m):=(\mathcal{U}_\rho+\mathcal{V}_\rho)(m)/\mathcal{V}(m)=T_m\mathcal{J}^{-1}(\rho)/
\mathcal{V}_\rho(m),  \qquad m \in \mathcal{J}^{-1}(\rho)
\]
of the vector bundle $\mathcal{U}_\rho/\mathcal{V}_\rho$,  
project surjectively to  $T_{\pi_\rho(m)}M_\rho$. Like in \S\ref{red},  for each
$G$-invariant $X\in\Gamma(\mathcal{U}_\rho)$  we
identify  $\widehat{X}=X(\operatorname{mod}\mathcal{V_\rho})$ with the section $\bar{X}$ of $M_\rho$ such that $T\pi_\rho \circ 
X=\bar{X}\circ\pi_\rho$. Write each $G$-invariant
$\alpha\in \Gamma(\mathcal{V}_\rho^\circ)$  as $\alpha=\pi_\rho^*\bar{\alpha}$ for some $\bar{\alpha} \in \Omega^1\left(M_\rho\right) $. 

Next we show that $D_\rho$ is the graph of a nondegenerate $2$-form. We begin by giving a formula for this $2$-form $\omega_\rho$. Let $\bar{X}, \bar{Y}\in \mathfrak{X}(M_\rho)$ and choose $G$-invariant $X, Y\in \mathfrak{X}(\mathcal{J}^{-1}(\rho))$ that are $\pi_\rho$-related to $\bar{X}$ and $\bar{Y}$, respectively. Write $X=\tilde{X}+V$ and $Y=\tilde{Y}+W$ with $\tilde{X}, \tilde{Y}\in \Gamma(\mathcal{U}_\rho)^G$ and $V,W\in\Gamma(\mathcal{V}_\rho)^G$. Then $\tilde{X}$ and $\tilde{Y}$ are also $\pi_\rho$-related to $\bar{X}$ and $\bar{Y}$ and we can write,
using the existence of $\bar\alpha\in \Omega^1(M_\rho)$ such that $(\bar
X,\bar\alpha)\in\Gamma(D_\rho)$,
\[
\omega_\rho(\bar X,\bar Y)=\bar\alpha(\bar Y)
                          =(\pi_\rho^*\bar\alpha)(\tilde Y)
                          =\omega_{\mathcal{U}_\rho}(\tilde X,\tilde Y),
\]
since $\tilde X$ has to be the (unique) section of $\mathcal{U}_\rho$ associated to the $1$-form $\pi_\rho^*\bar\alpha$ (see \eqref{lem:dimDcapKperp}) and where $\omega_{\mathcal{U}_\rho} $ is the restriction of $\omega_{\mathcal{J}^{-1}(\rho)}$ to $\mathcal{U}_\rho \times \mathcal{U}_\rho$.

We prove that $\omega_\rho$ is nondegenerate.
Let $\bar{X}\in\mathfrak{X}(M_\rho)$ with $\omega_\rho(\bar{X},\bar{Y})=0$ for all $\bar{Y}\in\mathfrak{X}(M_\rho)$. Choose a $G$-invariant section $\tilde{X} \in\Gamma(\mathcal{U}_\rho)$ as above. Extend $\tilde{X}$ to a local vector field $X $ on $M $, that is, $X\in \Gamma(\mathcal{U})\subseteq\mathfrak{X}(M)$ satisfies $\tilde X\sim_{i_\rho}X$. For $m\in\mathcal{J}^{-1}(\rho)$ and $v\in \mathcal{U}(m)\subseteq T_mM$ we have
\[
\omega_{\mathcal{H}}(m)(X(m),v)
=\omega_{\mathcal{U}}(m)(\tilde X(m),v)
={\omega_\rho}\left(\pi_\rho(m)\right)(\bar X({\pi_\rho(m)}),T_m\pi_\rho(v))
=0.
\]
Thus, the vector $X(m)$ is an element of
$\mathcal{U}(m)=(\mathcal{H}\cap\mathcal{V})^{\omega_{\mathcal{H}}}(m)$ that is 
$\omega_{\mathcal{H}}(m)$-orthogonal to all
  $v\in\mathcal{U}(m)$
  and hence lies in
  $((\mathcal{H}\cap\mathcal{V})^{\omega_{\mathcal{H}}})^{\omega_{\mathcal{H}}}(m)$.
Since $\omega_{\mathcal{H}}$ is nondegenerate, we have
$((\mathcal{H}\cap\mathcal{V})^{\omega_{\mathcal{H}}})^{\omega_{\mathcal{H}}}
= \mathcal{H}\cap\mathcal{V}$.
This yields $\tilde X(m)=X(m)\in(\mathcal{H}\cap\mathcal{V})(m)$ and thus the vector
$\bar X(m)$ is zero in $T_mM_\rho$.

(ii) Recall that, since $\mathsf{G}_0=\{0\}$, the solution $X_h$ of the implicit
  Hamiltonian system $(X,\mathbf{d} h)\in \Gamma(D)$ is unique: if $Y$ is another solution, then $Y-X_h\in\Gamma(\mathsf{G}_0)=\{0_M\}$. 

We know already that $X_h \in \Gamma(\mathcal{U}_\rho)$.
Furthermore, we have for  all $Y\in \Gamma(\mathcal{U}_\rho)$,
$V\in\Gamma(\mathcal{V}_\rho\cap\mathcal{H}_\rho)$ and all $m\in\mathcal{J}^{-1}(\rho)$ 
\begin{align*}
\omega_{\mathcal{J}^{-1}(\rho)}(m)(X_h(m),Y(m)+V(m))
&=\omega_M(m)(X_h(m),Y(m)+V(m))=\mathbf{d} h_m(Y(m)+V(m))\\
&= (i_\rho^*\mathbf{d} h)(m)(Y(m)+V(m))
\end{align*}
and the assertion follows. 

(iii)  The fact that the flow of $X_h$ leaves $\mathcal{J}^{-1}(\rho)$ invariant follows from the preceding statement since we have $X_h\in\Gamma(\mathcal{D}_G)$. By $G $-invariance of $D$ we have
$(\Phi_g^*X_h,\Phi_g^*\mathbf{d} h)\in \Gamma(D)$ for all $g\in G$. Since $h$ is $G$-invariant, the equality $\Phi_g^*\mathbf{d}
h=\mathbf{d}
\Phi_g^*h=\mathbf{d} h$ holds and thus we have
$\Phi_g^*X_h-X_h\in\Gamma(\mathsf{G}_0)=\{0_M\}$. The vector
field $X_h$ is consequently $G$-equivariant and its flow commutes with the
$G$-action.

(iv) Since $X_h\in \Gamma(\mathcal{U}_\rho)$ and $i_\rho^*\mathbf{d} h\in
\mathcal{V}^\circ_\rho$, we have
  \[
  (X_h,\mathbf{d} h)\in \Gamma\left(D_{\mathcal{J}^{-1}(\rho)}\cap
(T\mathcal{J}^{-1}(\rho)\oplus
    \mathcal{V}^\circ_\rho)\right).
    \] 
    
The flow $F_t^\rho$ on $M _\rho$ induces a vector field $X_{h_\rho}\in\mathfrak{X}(M_\rho)$. Therefore, taking the $t $-derivative of the relation in (iii) we get
\[
X_{h_\rho}(\pi_\rho(m))
=\left.\frac{d}{dt}\right|_{t=0}F_t^\rho(\pi_\rho(m))
=\left.\frac{d}{dt}\right|_{t=0}\left(\pi_\rho\circ F_t\right)(m)
=T_{m}\pi_\rho X_h(m), 
\]
that is, $X_h\sim_{\pi_\rho}X_{h_\rho}$. Choose $\bar{Y}\in \mathfrak{X}(M_\rho)$, $Y\in \Gamma(\mathcal{U}_\rho)^G$, and $V\in
\Gamma(\mathcal{V}_\rho)^G$ such that
$T\pi_\rho \circ (Y+V)=\bar{Y}\circ \pi_\rho$. Then, for all 
$m\in \mathcal{J}^{-1}(\rho)$ we get
\begin{align*}
\omega_\rho(\pi_\rho&(m))\left(X_{h_\rho}(\pi_\rho(m)),
\bar{Y}(\pi_\rho(m))\right)
=(\pi_\rho^*\omega_\rho)(m)(X_h(m),Y(m)+V(m))
=(\pi_\rho^*\omega_\rho)(m)(X_h(m),Y(m))\\
&=\omega_{\mathcal{J}^{-1}(\rho)}(m)(X_h(m),Y(m))
= (i_\rho^*\mathbf{d} h)_m(Y(m))
= (\pi_\rho^*\mathbf{d} h_\rho)_m(Y(m)+V(m))
=(\mathbf{d} h_\rho)_{\pi_\rho(m)}(\bar Y(\pi_\rho(m))),
\end{align*}
so we have ${\mathbf{i}}_{X_{h_\rho}}\omega_\rho=\mathbf{d} h_\rho$, as claimed.
\end{proof}

\subsection{Examples of optimal reduction for nonholonomic systems}
\label{sec:exaples_optimal_reduction_nonholonomic}

\subsubsection{The constrained particle in space}
We return to the example treated in \S\ref{constrpart} and use the same notations and conventions. The distribution $\mathcal{V}\cap
\mathcal{H}$ is pointwise the span of
the vector field $\partial_x+y\partial_z$. Since $\mathcal{V}^\circ$ is
spanned by the covector fields $\mathbf{d} y$, $\mathbf{d} p_x$, and $\mathbf{d} p_y$, the considerations in \S\ref{constrpart} yield ${\mathbf{i}}_{\partial_x+y\partial_z}\omega_M=-(1+y^2)\mathbf{d}p_x-yp_x\mathbf{d} y\in\Gamma(\mathcal{V}^\circ)$ and hence $(\partial_x+y\partial_z,
-(1+y^2)\mathbf{d} p_x-yp_x\mathbf{d} y)\in\Gamma(D\cap(\mathcal{V}\oplus\mathcal{V}^\circ))$. Hence the distribution $\mathcal{D}_G$ is in this case  $\ker \{-(1+y^2)\mathbf{d} p_x-yp_x\mathbf{d} y\}=\ker \{\mathbf{d} f\}$, where 
$f(x,y,z,p_x,p_y)=\sqrt{1+y^2}p_x$ is the constant of motion (in agreement with \cite{BaSn93}). Note that by Example \ref{ex:sec_71},  $\mathbf{d}f $ is the $1$-form giving the Nonholonomic Noether Theorem. Hence  
\begin{equation*}
\mathcal{D}_G=\operatorname{span}
\{\partial_{p_y},\partial_x,\partial_z,yp_x\partial_{p_x}-(y^2+1)\partial_y\}
\end{equation*}
is obviously involutive and constant dimensional (and consequently integrable). This shows that $M/ \mathcal{D}_G = \mathbb{R}$.
The Dirac structure on a leaf $f^{-1}(\mu)$, $\mu\in M/\mathcal{D}_G = \mathbb{R}$,  of this distribution
is given by
\begin{align*}
D_{f^{-1}(\mu)}=\left\{(X,\alpha)\in\Gamma(Tf^{-1}(\mu)\oplus T^*f^{-1}(\mu))\mid
X\in\Gamma(\mathcal{H}\cap\mathcal{D}_G),\alpha-{\mathbf{i}}_{X}i_{\mu}^*\omega_M\in\Gamma((\mathcal{H}\cap\mathcal{D}_G)^\circ)\right\}
\end{align*}
and a computation yields
\begin{align*}
D_{f^{-1}(\mu)}=\operatorname{span}\left\{\left(\partial_{p_y},-\mathbf{d}
y\right),\left(\partial_x+y\partial_z,0\right)\left(0,\mathbf{d} z-y\mathbf{d}
x\right),\left((1+y^2)\partial_y-yp_x\partial_{p_x},(1+y^2)\mathbf{d}
    p_y\right)\right\}
\end{align*}
because the $1$-form $(1+y^2)\mathbf{d} p_x+yp_x\mathbf{d} y$ vanishes on
$T\left(f^{-1}(\mu)\right)$.
Since $G$ is in this case connected, we have
$G_\mu=G$ (see Remark \ref{rem}). Consider the codistribution $\mathcal{V}^\circ$ on $f^{-1}(\mu)$ and get
\begin{align*}
D_{f^{-1}(\mu)}\cap (Tf^{-1}(\mu)\oplus
\mathcal{V}^\circ)=\operatorname{span}&\left\{\left(\partial_{p_y},-\mathbf{d}
y\right),\left(\partial_x+y\partial_z,0\right) \left((1+y^2)\partial_y-yp_x\partial_{p_x},(1+y^2)\mathbf{d}
    p_y\right)\right\}.
\end{align*}
Hence the reduced Dirac structure $D_{\mu}$ on $M_\mu=f^{-1}(\mu)/G$ is given by
the formula
\begin{align*}
D_{\mu}&=\left.\frac{\left[D_{f^{-1}(\mu)}\cap \left(T\left(f^{-1}(\mu)\right)\oplus \mathcal{V}^\circ\right) \right] +
  \mathcal{V}\oplus \{0\}}{\mathcal{V}\oplus \{0\}}\right/G\\
&= \operatorname{span}\left\{\left(\partial_{p_y},-\mathbf{d}
y\right),\left((1+y^2)\partial_y-yp_x\partial_{p_x},(1+y^2)\mathbf{d}
    p_y\right)\right\}.
\end{align*}
This corresponds exactly to a symplectic leaf (with its associated Dirac structure) of the Poisson structure \eqref{ex1_D_red} obtained in the first part of this example (see \S \ref{constrpart}).

Finally we compute $\mathcal{R}$ for this example. Since $\mathcal{H}^\circ$ is one-dimensional, we get  $\mathcal{R}=\mathcal{H}^\circ$ or
$\mathcal{R}=\{0\}$. Recall that $D$ is the span of  
\begin{eqnarray*}
\left\{\left(\partial_{p_y},-\mathbf{d}
y\right),\left(\partial_x+y\partial_z,(1+y^2)\mathbf{d} p_x+yp_x\mathbf{d}
y\right),\left(0, \mathbf{d} z-y\mathbf{d} x\right), 
\left(\partial_y, \mathbf{d} p_{y}-p_{x}\mathbf{d}
    z\right),\left(\partial_{p_x}, -y\mathbf{d} z-\mathbf{d} x\right)\right\}
\end{eqnarray*} 
where we have computed:
\begin{eqnarray*}
{\mathbf{i}}_{\partial_x+y\partial_z}\omega_M&=&(1+y^2)\mathbf{d} p_x+yp_x\mathbf{d} y\\
{\mathbf{i}}_{\partial_y}\omega_M&=&\mathbf{d} p_{y}-p_{x}\mathbf{d} z\\
{\mathbf{i}}_{\partial_{p_y}}\omega_M&=&-\mathbf{d} y\\
{\mathbf{i}}_{\partial_{p_x}}\omega_M&=&-y\mathbf{d} z-\mathbf{d} x.\\
\end{eqnarray*}
Since $\mathcal{U} = \operatorname{span}\{\partial_{p_y},(1+y^2)\partial_y-yp_x\partial_{p_x}, 
\partial_x+y\partial_z \}$, we conclude from 
\begin{eqnarray*}
{\mathbf{i}}_{(1+y^2)\partial_y-yp_x\partial_{p_x}}\omega_M&=&(1+y^2)(\mathbf{d}
p_{y}-p_{x}\mathbf{d} z)-yp_x(-y\mathbf{d} z-\mathbf{d} x)\\
&=&(1+y^2)\mathbf{d} p_y-p_x(\mathbf{d} z-y\mathbf{d} x),
\end{eqnarray*}
that the distribution $\mathcal{R}$ is equal to $\mathcal{H}^\circ$. The
constant of motion that we have found above is a \emph{gauge constant of motion}.

\subsubsection{The vertical rolling disk}
In this subsection we shall determine the Nonholonomic Momentum Equations for the example of the vertical rolling disk studied in \S\ref{ex_vertical_disk}.
The Dirac structure for this nonholonomic system is given by
\begin{eqnarray*}
D=\operatorname{span}\left\{\left(\partial_\phi,\mathbf{d} p_\phi+\frac{\mu
R\sin\phi}{I}p_\theta\mathbf{d} x-\frac{\mu R\cos\phi}{I}p_\theta\mathbf{d} y
\right),\left(\partial_{p_\phi},-\mathbf{d}\phi \right),\right.\\
\left(\partial_{p_\theta},-\frac{\mu R\cos\phi}{I}\mathbf{d} x-\frac{\mu
R\sin\phi}{I}\mathbf{d}
  y-\mathbf{d} \theta\right),\left(0, \mathbf{d} x-R\cos\phi\mathbf{d} \theta
\right), \\
\left.\left(0, \mathbf{d} y-R\sin\phi\mathbf{d}
\theta\right),\left(\partial_\theta+R\cos\phi
\partial_x+R\sin\phi\partial_y,(1+\frac{\mu R^2}{I})\mathbf{d} p_\theta
\right)\right\}
\end{eqnarray*} 
where we have computed
\begin{align*}
{\mathbf{i}}_{\partial_\phi}\omega_M&=\mathbf{d} p_\phi+\frac{\mu
R\sin\phi}{I}p_\theta\mathbf{d} x-\frac{\mu R\cos\phi}{I}p_\theta\mathbf{d} y\\
{\mathbf{i}}_{\partial_{p_\theta}}\omega_M&=-\frac{\mu R}{I}\cos\phi\mathbf{d} x-\frac{\mu
R}{I}\sin\phi\mathbf{d} y-\mathbf{d} \theta\\
{\mathbf{i}}_{\partial_{p_\phi}}\omega_M&=-\mathbf{d}\phi\\
{\mathbf{i}}_{\partial_\theta+R\cos\phi \partial_x+R\sin\phi\partial_y}\omega_M&=\mathbf{d}
p_\theta+R\cos\phi\left(\frac{\mu R\cos\phi}{I}\mathbf{d}
p_\theta-\frac{\mu R\sin\phi}{I}\mathbf{d}\phi\right) +R\sin\phi\left(+\frac{\mu
R\sin\phi}{I}\mathbf{d}
p_\theta+\frac{\mu R\cos\phi}{I}\mathbf{d}\phi\right)\\
&=\left(1+\frac{\mu R^2}{I} \right)\mathbf{d} p_\theta.
\end{align*}
We consider again the three possible Lie group actions:
\begin{enumerate}
\item \emph{The case $G=\mathbb{R}^2$} (\cite{CaDiLeMa98})

Here, $\mathcal{V}^\circ=\operatorname{span}\{\mathbf{d}
p_\phi,\mathbf{d} \phi,\mathbf{d} p_\theta,\mathbf{d} \theta\}$ but there are no
nontrivial horizontal 
symmetries and hence the distribution $\mathcal{D}_G$ is simply the
whole bundle $TM$. We next compute $\mathcal{R}$. The vector bundle  $ D\cap
\mathcal{K}^\perp$ is
given in this case by \begin{eqnarray*}
\operatorname{span}\left\{\left(\partial_{\phi},\mathbf{d}
     p_\phi\right),\left(\partial_{p_\theta},-\left(1+\frac{\mu
R^2}{I}\right)\mathbf{d}\theta\right),\left(\partial_{p_\phi},-\mathbf{d}\phi\right)
\right., 
\left.\left(\partial_\theta+R\cos\phi \partial_x+R\sin\phi\partial_y,\left(1+\frac{\mu
R^2}{I}\right)\mathbf{d} p_\theta\right)\right\}.
\end{eqnarray*}
To get this, we have added 
\begin{equation}\label{add1}
\frac{\mu R}{I}p_\theta\cos\phi(\mathbf{d} x-R\cos\phi\mathbf{d} \theta)+\frac{\mu
R}{I}\sin\phi(\mathbf{d} y-R\sin\phi\mathbf{d}
\theta)
\end{equation}
 to ${\mathbf{i}}_{\partial_{p_\theta}}\omega_M$ and
\begin{equation}\label{add2}
-\frac{\mu R}{I}p_\theta\sin\phi(\mathbf{d} x-R\cos\phi\mathbf{d}
\theta)+\frac{\mu R}{I}\cos\phi(\mathbf{d} y-R\sin\phi\mathbf{d} \theta)
\end{equation}
 to
${\mathbf{i}}_{\partial_\phi}\omega_M$. 
This yields $\mathcal{R}=\mathcal{H}^\circ$ and thus the distribution 
$\mathcal{R}^\circ\cap\mathcal{V}$ is equal to
$\mathcal{H}\cap\mathcal{V}$ and hence trivial.

\item \emph{The case $G=\operatorname{SE}(2)$ } (\cite{Bloch03})

In this case, we have $\mathcal{V}^\circ=\operatorname{span}\{\mathbf{d}
p_\phi,\mathbf{d} p_\theta,\mathbf{d} \theta\}$ and $\mathcal{H}\cap\mathcal{V} = \operatorname{span}\{\partial_\phi\}$. We get
 $\mathcal{H}\cap\mathcal{V}\subseteq
(\mathcal{H}\cap\mathcal{V})^{\omega_M}$. A direct computation gives
\[
{\mathbf{i}}_{\partial_\phi}\omega_M=\mathbf{d} p_\phi+\frac{\mu
R\sin\phi}{I}p_\theta\mathbf{d}
x-\frac{\mu R\cos\phi}{I}p_\theta\mathbf{d} y.
\] 
Adding
\begin{align*}
\frac{\mu R}{I}p_\theta(\cos\phi\mathbf{d}
  y-\sin\phi\mathbf{d} x) d 
=\frac{\mu R}{I}p_\theta\left(\cos\phi(\mathbf{d}
  y-R\sin\phi\mathbf{d}\theta)-\sin\phi(\mathbf{d}
x-R\cos\phi\mathbf{d}\theta)\right)\in\Gamma(\mathcal{H}^\circ)
\end{align*}
to this expression we see that $\Gamma(D\cap(\mathcal{V}\oplus\mathcal{V}^\circ))$ is
spanned by $(\partial_{\phi},\mathbf{d}
p_\phi)$ and thus the
distribution $\mathcal{D}_G=\ker \mathbf{d} p_{\phi}$ is obviously integrable.
For a value $\rho \in \mathbb{R}$ of the map $p_\phi$, the reduced Dirac
structure  on $M_\rho$ is spanned by 
\[
\left(\partial_\theta,\left(1+\frac{\mu
R^2}{I}\right)\mathbf{d} p_\theta\right) \quad \text{and} \quad
\left( \partial_{p_\theta},-\left(1+\frac{\mu R^2}{I}\right)\mathbf{d}\theta
\right).
\]

The Nonholonomic Noether Theorem yields a constant of motion but this constant
doesn't arise from an element of $\mathfrak{g}$ whose corresponding fundamental
vector field is lying in $\Gamma(\mathcal{V}\cap\mathcal{R}^\circ)$:
we have computed in \S \ref{nonhol} that, in this case, $\mathcal{U}$ is the span of
the
three vector fields $\partial_{\phi}$, $\partial_{p_\theta}$, and
$\partial_\theta+R\cos\phi\partial_x+R\sin\phi\partial_y$. Again, we have to
add \eqref{add1} to ${\mathbf{i}}_{\partial_{p_\theta}}\omega_M$ and
\eqref{add2} to
${\mathbf{i}}_{\partial_\phi}\omega_M$ in order to get sections of
$\mathcal{V}^\circ$. Thus, we need the whole
of $\mathcal{H}^\circ$ in the construction of $D\cap\mathcal{K}^\perp$.

\item \emph{The case $G=\mathbb{S}^1\times\mathbb{R}^2$ } (\cite{Bloch03})

Here, we have $\mathcal{V}^\circ=\operatorname{span}\{\mathbf{d}
p_\phi,\mathbf{d} p_\theta,\mathbf{d} \phi\}$ and $\mathcal{H}\cap\mathcal{V}$ is
again 
one-dimensional: this time it is the span of the vector field
$\partial_\theta+R\cos\phi\partial_x+R\sin\phi\partial_y$. Thus,
$\Gamma(D\cap(\mathcal{V}\oplus\mathcal{V}^\circ))$ is spanned by
\[\left(\partial_\theta+R\cos\phi\partial_x+R\sin\phi\partial_y,(1+\frac{\mu
R^2}{I})\mathbf{d} p_\theta\right)\]
 and the
distribution $\mathcal{D}_G=\ker \{(1+\frac{\mu R^2}{I})\mathbf{d} p_\theta\} $ is
again integrable. 
For a value $\rho \in \mathbb{R}$ of the map $p_\theta$, the reduced Dirac
structure  on $M_\rho$ is spanned by $(\partial_\phi,\mathbf{d} p_\phi)$ and
$(\partial_{p_\phi},-\mathbf{d} \phi)$.

In this case, we have $\mathcal{U}=\operatorname{span}\{\partial_{\phi},\partial_{p_\phi},
\partial_\theta+R\cos\phi\partial_x+R\sin\phi\partial_y\}$.
We get $\mathcal{R}=\operatorname{span}\{\frac{\mu R}{I}p_\theta(\sin\phi\mathbf{d}
x-\cos\phi\mathbf{d}
  y)\}$ from the considerations for the second case. Thus we have
  $\mathcal{R}^\circ=\operatorname{span}\{\partial_{\phi},\partial_{p_\phi},\partial_{p_\theta},
\partial_\theta,\cos\phi\partial_x+\sin\phi\partial_y\}$ and our constant of motion
$p_\theta$ really arises from a
  fundamental vector field lying in $\mathcal{R}^\circ$.

\item \emph{The case $G=\operatorname{SE}(2)\times\mathbb{S}^1$ } (\cite{Bloch03})

In this last case, we have $\mathcal{V}^\circ=\operatorname{span}\{\mathbf{d}
p_\phi,\mathbf{d} p_\theta\}$ and $\mathcal{H}\cap\mathcal{V}$ is this time
two-dimensional: it is the span of the vector fields
$\partial_\theta+R\cos\phi\partial_x+R\sin\phi\partial_y$ and $\partial_\phi$. Thus,
$\Gamma(D\cap(\mathcal{V}\oplus\mathcal{V}^\circ))$ is spanned by
$(\partial_\theta+R\cos\phi\partial_x+R\sin\phi\partial_y,(1+\frac{\mu
R^2}{I})\mathbf{d}
p_\theta)$ and $(\partial_\phi,\mathbf{d} p_\phi)$ and the
distribution $\mathcal{D}_G=\ker \{(1+\frac{\mu R^2}{I})\mathbf{d} p_\theta,
\mathbf{d} p_\phi\} $ is
integrable. Here, the reduced manifolds are single points.

We have $\mathcal{U}=\operatorname{span}\{\partial_{\phi},
\partial_\theta+R\cos\phi\partial_x+R\sin\phi\partial_y\}$ and we get as above 
$\mathcal{R}^\circ=\operatorname{span}\{\partial_{\phi},\partial_{p_\phi},\partial_{p_\theta},
\partial_\theta,\cos\phi\partial_x+\sin\phi\partial_y\}$.
\end{enumerate}

\subsubsection{The Chaplygin skate}
We continue here the examples of \S\ref{chaplygin}.  
\medskip

\noindent \textbf{The standard Chaplygin skate.} We have seen that
$\mathcal{V}\cap\mathcal{H}=\operatorname{span}\{\partial_\theta,\cos\theta\partial_x+\sin\theta\partial_y\}$.
If we choose the basis $\xi^1:=(1,0,0)$, $\xi^2:=(0,1,0)$, $\xi^3:=(0,0,1)$
of the Lie algebra $\mathfrak {se}(2)$ we get
$\xi^1_M=\partial_\theta-y\partial_x+x\partial_y$, $\xi^2_M=\partial_x$, and
$\xi^3_M=\partial_y$. Hence, the sections $\xi^1+y\xi^2-x\xi^3$ and
$\cos\theta\xi^2+\sin\theta\xi^3\in\Gamma(\mathfrak{g}^{\mathcal{H}})$ are spanning
sections of $\mathfrak{g}^{\mathcal{H}}$ and the corresponding Nonholonomic Noether
equations are $s\cos\theta\mathbf{d} p_y-s\sin\theta\mathbf{d} p_x$ and $\cos\theta\mathbf{d}
p_x+\sin\theta\mathbf{d} p_y$ respectively. Consequently, the two
spanning sections
$-s^{-1}\sin\theta\xi^1+(\cos^2\theta-s^{-1}y\sin\theta)\xi^2+\sin\theta(\cos\theta+s^{-1}x)\xi^3$
and
$s^{-1}\cos\theta\xi^1+\cos\theta(\sin\theta+s^{-1}y)\xi^2+(\sin^2\theta-s^{-1}x\cos\theta)\xi^3$
of $\mathfrak{g}^{\mathcal{H}}$ lead to the nonholonomic Noether equations $\mathbf{d}
p_x$ and $\mathbf{d} p_y$ respectively. Thus, 
$\mathcal{D}_G=\mathcal{U}+\mathcal{V}=\operatorname{span}\{\partial_x,\partial_y,\partial_\theta\}$
is found easily because $\mathcal{D}_G$ is the kernel of $\{\mathbf{d} p_x,\mathbf{d}
p_y\}$. This is obviously integrable. The induced Dirac structure on a leaf
$f^{-1}(a,b)$ (where $f$ is the projection on $(p_x,p_y)$) of
$\mathcal{D}_G$ is given by 
\[D_{f^{-1}(a,b)}=\operatorname{span}\{(\cos\theta\partial_x+\sin\theta\partial_y,0),
(\partial_\theta,0), (0,\sin\theta\mathbf{d} x-\cos\theta\mathbf{d} y)\}.\]
Here the reduced space $M_{(a,b)}$ is a single point. The reduced Dirac
structure is hence trivial, as can  also be  seen from the formula
$\frac{(D_{f^{-1}(a,b)}\cap\mathcal{K}_{(a,b)}^\perp)+\mathcal{K}_{(a,b)}}{\mathcal{K}_{(a,b)}}/G_{(a,b)}$.

\medskip

We could also consider the action of $\mathbb{S}^1$ on $M$ given by
$\Phi:\mathbb{S}^1\times M\to M$, $(\alpha,\theta,x,y,p_x,p_y)\mapsto
(\theta+\alpha, x\cos\alpha -  y \sin\alpha, x\sin\alpha+ y\cos\alpha, p_x \cos\alpha -p_y \sin\alpha ,p_x \sin\alpha +p_y \cos\alpha )$. Here, we would have
$\mathcal{V}\cap\mathcal{H}=\{0\}$ except for the points satisfying
$x=-\sin\theta$ and $y=-\cos\theta$, so the condition that
$\mathcal{V}\cap\mathcal{H}$ has constant rank is not satisfied (we
have also $\mathcal{V}+\mathcal{H}\neq TM$).

If we consider the action of $\mathbb{R}^2$ on $M$ given by $\Phi:\mathbb{R}^2\times M\to M$, $(r,s,\theta,x,y,p_x,p_y)\mapsto
(\theta, x+r,y+s,p_x,p_y)$, we have
$\mathcal{V}=\operatorname{span}\{\partial_x,\partial_y\}$. Hence,
$\mathcal{V}+\mathcal{H}= TM$ and
$\mathcal{V}\cap\mathcal{H}=\{\cos\theta\partial_x+\sin\theta\partial_y\}$ has
constant rank on $M$. The distribution $\mathcal{D}_G$ is given by
$\mathcal{D}_G=(\ker\{\cos\theta\mathbf{d} p_x+\sin\theta\mathbf{d}
p_y\} \cap \mathcal{H})+\mathcal{V}=\operatorname{span}\{\sin\theta\partial_{p_x}-\cos\theta\partial_{p_y},\partial_\theta,\partial_x,\partial_y\}$.
This vector bundle is not involutive and hence it is not integrable. Since $\mathcal{U}=\ker\{\cos\theta\mathbf{d} p_x+\sin\theta\mathbf{d}
p_y\}\cap\mathcal{H}=\operatorname{span}\{\sin\theta\partial_{p_x}-\cos\theta\partial_{p_y},\partial_\theta,\cos\theta\partial_x+\sin\theta\partial_y\}$,
it is easy to see that $\mathcal{R}=\mathcal{H}^\circ$ and hence
$\mathcal{R}^\circ\cap\mathcal{V}=\mathcal{V}\cap\mathcal{H}=\{\cos\theta\partial_x+\sin\theta\partial_y\}$,
which confirms the fact that the nonholonomic Noether equation yields in this
case no constant of motion.

\medskip 

\noindent \textbf{The Chaplygin skate with a rotor on it.} We have
$\mathcal{V}\cap\mathcal{H}=\operatorname{span}\{\partial_\phi,\partial_\theta,\cos\theta\partial_x+\sin\theta\partial_y\}$.
If we choose the basis $\xi^1:=(1,0,0,0)$, $\xi^2:=(0,1,0,0)$,  $\xi^3:=(0,0,1,0)$,
$\xi^4:=(0,0,0,1)$
of the Lie algebra $\mathbb{R} \times\mathfrak {se}(2)$ we get
$\xi^1_M=\partial_\phi$, $\xi^2_M=\partial_\theta-y\partial_x+x\partial_y$, $\xi^3_M=\partial_x$, and
$\xi^4_M=\partial_y$. 
Hence, the sections $\xi^1$, $\xi^2+y\xi^3-x\xi^4$, and
$\cos\theta\xi^3+\sin\theta\xi^4\in\Gamma(\mathfrak{g}^{\mathcal{H}})$ are spanning
sections of $\mathfrak{g}^{\mathcal{H}}$ and the corresponding Nonholonomic Noether
equations are $\mathbf{d} p_\phi$, $s\cos\theta\mathbf{d} p_y-s\sin\theta\mathbf{d} 
p_x+\mathbf{d} p_\phi$, and $\cos\theta\mathbf{d}
p_x+\sin\theta\mathbf{d} p_y$, respectively.
Thus, the three
spanning sections $\xi^1$,
$s^{-1}\sin\theta\xi^1-s^{-1}\sin\theta\xi^2+(\cos^2\theta-s^{-
1}y\sin\theta)\xi^3+\sin\theta(\cos\theta+s^{-1}x)\xi^4$,
and
$-s^{-1}\cos\theta\xi^1+s^{-1}\cos\theta\xi^2+\cos\theta(\sin\theta+s^{-
1}y)\xi^3+(\sin^2\theta-s^{-1}x\cos\theta)\xi^4$
of $\mathfrak{g}^{\mathcal{H}}$ lead to the nonholonomic Noether equations $\mathbf{d} p_\phi$, 
$\mathbf{d} p_x$, and $\mathbf{d} p_y$, respectively.
 Thus, 
$\mathcal{D}_G=\mathcal{U}+\mathcal{V}=\operatorname{span}\{\partial_x,\partial_y,\partial_\theta,\partial_\phi\}$
is found easily because $\mathcal{D}_G$ is the kernel of $\{\mathbf{d} p_\phi,\mathbf{d} p_x,\mathbf{d}
p_y\}$. This is obviously integrable. The induced Dirac structure on a leaf
$f^{-1}(a,b,c)$ (where $f$ is the projection on $(p_\phi,p_x,p_y)$) of
$\mathcal{D}_G$ is given by 
\[D_{f^{-1}(a,b,c)}=\operatorname{span}\{(\cos\theta\partial_x+\sin\theta\partial_y,0),
(\partial_\theta,0),(\partial_\phi,0), (0,\sin\theta\mathbf{d} x-\cos\theta\mathbf{d} y)\}.\]
Here the reduced space $M_{(a,b,c)}$ is a single point. The reduced Dirac
structure is hence trivial, as can  also be  seen from the formula
$\frac{(D_{f^{-1}(a,b,c)}\cap\mathcal{K}_{(a,b,c)}^\perp)+\mathcal{K}_{(a,b,c)}}{\mathcal{K}_{(a,b,c)}}/G_{(a,b,c)}$.

\medskip

Finally, note that in the last two examples, we have $\mathcal{R}=\{0\}$ and
hence $\mathcal{R}^\circ\cap\mathcal{V} = \mathcal{V}$. 
This is why we get in the first of
the two 
examples the three constants of motion $p_x$, $p_y$, and $p_\theta$ belonging to
the three elements $\xi^2$, $\xi^3$, and $\xi^1$ of $\mathfrak{g}$ and in the second
example
the constants  $p_x$, $p_y$,  $p_\theta$, and $p_\phi$ belonging to
the four elements $\xi^3$, $\xi^4$, $\xi^2$, and $\xi^1$ of $\mathfrak{g}$. Note
that in this case, the constancy of $p_\phi$ follows already from the existence of the constant section
$\xi^1$ of $\mathfrak{g}^{\mathcal{H}}$. 

\medskip

Like in the previous example, the other symmetry groups of the system
(the ``$\theta$-symmetry'' $\mathbb{S}^1$, ``the $\phi$-symmetry''
$\mathbb{S}^1$, $\mathbb{S}^1\times\mathbb{S}^1 $,
$\mathbb{S}^1\times\mathbb{R}^2 $, $\operatorname{SE}(2)$) are not interesting for the method of
reduction presented in this section.

\subsubsection{The Heisenberg particle}
At last, we present an example where the reduced form is not closed. It can be found in \cite{Bloch03}. The configuration space $Q$ is
$\mathbb{R}^3$ with coordinates $(x,y,z)$  subject to the constraint $\dot z=y\dot
x-x\dot y$. The Lagrangian on $TQ$ is given by $L(x,y,z,\dot x,\dot y, \dot
z)=\frac{1}{2}(\dot x^2+\dot y^2+\dot z^2)$ and hence the Legendre transformation yields
\begin{equation*}
p_x=\dot x, \qquad p_y=\dot y,\qquad p_z=\dot z.
\end{equation*}
For $(x,y,z,p_x,p_y,p_z)$, we have $p_z=yp_x-xp_y$. Hence, we have the global
coordinates $(x,y,z,p_x,p_y)$ for $M$ and the $2$-form $\omega_M$ is given by
\begin{align*}
\omega_M&=\mathbf{d} x\wedge \mathbf{d} p_x+\mathbf{d} y\wedge \mathbf{d} p_y+\mathbf{d} z\wedge \mathbf{d}(yp_x-xp_y)\\
&=\mathbf{d} x\wedge (\mathbf{d} p_x +p_y\mathbf{d} z)
+\mathbf{d} y\wedge (\mathbf{d} p_y-p_x\mathbf{d} z)
+\mathbf{d} z\wedge (y\mathbf{d} p_x-x \mathbf{d} p_y).
\end{align*}
The vector bundle $\mathcal{H}$ is given by $\mathcal{H}=\ker\{\mathbf{d} z-y\mathbf{d} x+x\mathbf{d} y\}=\operatorname{span}\{y\partial_z+\partial_x,x\partial_z-\partial_y,\partial_{p_x},\partial_{p_y}\}$
and  we compute 
\begin{align*}
{\mathbf{i}}_{y\partial_z+\partial_x}\omega_M&=y(y\mathbf{d} p_x-x \mathbf{d} p_y)-yp_y\mathbf{d} x+yp_x\mathbf{d} y+(\mathbf{d} p_x +p_y\mathbf{d} z)\\
{\mathbf{i}}_{x\partial_z-\partial_y}\omega_M&=x(y\mathbf{d} p_x-x \mathbf{d} p_y)-xp_y\mathbf{d} x+xp_x\mathbf{d} y-(\mathbf{d} p_y-p_x\mathbf{d} z)\\
{\mathbf{i}}_{\partial_{p_x}}\omega_M&=-\mathbf{d} x-y\mathbf{d} z\\
{\mathbf{i}}_{\partial_{p_y}}\omega_M&=-\mathbf{d} y+x\mathbf{d} z.
\end{align*}
Thus, we get the smooth global spanning sections
\begin{align*}
\big\{&\left(y\partial_z+\partial_x,(y^2+1)\mathbf{d} p_x-xy\mathbf{d} p_y-yp_y\mathbf{d} x+yp_x\mathbf{d}
  y+p_y\mathbf{d} z\right), 
\left(\partial_{p_x},-\mathbf{d} x-y\mathbf{d} z\right),\\
&\left(x\partial_z-\partial_y,-(x^2+1)\mathbf{d} p_y+xy\mathbf{d} p_x-xp_y\mathbf{d} x+xp_x\mathbf{d} y+p_x\mathbf{d} z\right),
\left(\partial_{p_y},-\mathbf{d} y+x\mathbf{d} z\right),\left(0,\mathbf{d} z-y\mathbf{d} x+x\mathbf{d} y\right)\big\}
\end{align*}
for the Dirac structure $D$.

Consider the action $\phi:\mathbb{R}\times Q\to Q$ of the Lie group $G=\mathbb{R}$ on $Q$, given by $\phi(r,x,y,z): =(x,y,z+r)$. This action obviously leaves the
Lagrangian and the constraints invariant. The induced action $\Phi: G\times M\to M$ is given by $\Phi(r,x,y,z,p_x,p_y): = (x,y,z+r,p_x,p_y)$ and
hence the vertical bundle in this example equals $\mathcal{V}=\operatorname{span}\{\partial_z\}$. We get
$\mathcal{H}\cap\mathcal{V}=\{0\}$ and hence $\mathcal{U}=\mathcal{H}$. The
two methods of reduction (in \S\ref{link-nonhol-dirac} and \S\ref{sec:optimal_reduction_nonholonomic}) lead in this case to the same result since the distribution
$\mathcal{D}_G=\mathcal{U}+\mathcal{V}=\mathcal{H}+\mathcal{V}=TM$ is
trivially integrable with $M$ as single leaf.
The reduced Dirac structure $D_{\rm red}$ on $\bar M$ with coordinates
$(x,y,p_x,p_y)$ is thus given by 
\begin{align*}
D_{\rm red}=&\left.\frac{(D\cap\mathcal{K}^\perp)+\mathcal{K}}{\mathcal{K}}\right/G\\
=&\operatorname{span}\big\{\left(\partial_x,(y^2+1)\mathbf{d} p_x-xy\mathbf{d} p_y-yp_y\mathbf{d} x+yp_x\mathbf{d}
  y+p_y(y\mathbf{d} x-x\mathbf{d} y)\right),
\left(\partial_{p_x},-\mathbf{d} x-y(y\mathbf{d} x-x\mathbf{d} y)\right),\\
&\quad\qquad\left(-\partial_y,-(x^2+1)\mathbf{d} p_y+xy\mathbf{d} p_x-xp_y\mathbf{d} x+xp_x\mathbf{d} y+p_x(y\mathbf{d} x-x\mathbf{d} y)\right),
\left(\partial_{p_y},-\mathbf{d} y+x(y\mathbf{d} x-x\mathbf{d} y)\right)\big\}\\
=&\operatorname{span}\big\{\left(\partial_x,(y^2+1)\mathbf{d} p_x-xy\mathbf{d} p_y+(yp_x-xp_y)\mathbf{d} y\right),
\left(\partial_{p_x},-(1+y^2)\mathbf{d} x+xy\mathbf{d} y\right),\\
&\quad\qquad\left(-\partial_y,-(x^2+1)\mathbf{d} p_y+xy\mathbf{d} p_x+(yp_x-xp_y)\mathbf{d} x\right),
\left(\partial_{p_y},-(1+x^2)\mathbf{d} y+xy\mathbf{d} x\right)\big\}.
\end{align*}
Note that this is the graph of the $2$-form 
\[
\omega_{\rm red}=(1+y^2)\mathbf{d}
x\wedge \mathbf{d} p_x+(1+x^2)\mathbf{d} y\wedge\mathbf{d} p_y+ (yp_x-xp_y)\mathbf{d} x\wedge\mathbf{d} y-xy(\mathbf{d}
x\wedge\mathbf{d} p_y+\mathbf{d} y\wedge\mathbf{d} p_x). 
\]
A direct computation shows that the determinant of
$\omega_{\rm red}$ equals $(1+x^2+y^2)^2 \neq 0$ on $\bar{M}$  which
shows that the form $\omega_{\rm red}$ is nondegenerate. The equalities 
\begin{align*}
\mathbf{d}\omega_{\rm red}(\partial_x,\partial_y,\partial_{p_x})=-2y \qquad
\text{and} \qquad \mathbf{d}\omega_{\rm red}(\partial_x,\partial_y,\partial_{p_y})=2x
\end{align*}
show that $\omega_{\rm red}$ is not closed.

Note also that in this example we have $\mathcal{R}=\mathcal{H}^\circ$ and hence $\mathcal{R}^\circ\cap\mathcal{V}=\mathcal{H}\cap\mathcal{V}=\{0\}$.

\def\cprime{$'$} \def\polhk#1{\setbox0=\hbox{#1}{\ooalign{\hidewidth
  \lower1.5ex\hbox{`}\hidewidth\crcr\unhbox0}}} \def\cprime{$'$}
  \def\cprime{$'$} \def\cprime{$'$} \def\cprime{$'$} \def\cprime{$'$}
  \def\cprime{$'$} \def\cprime{$'$}
  \def\polhk#1{\setbox0=\hbox{#1}{\ooalign{\hidewidth
  \lower1.5ex\hbox{`}\hidewidth\crcr\unhbox0}}}
  \def\polhk#1{\setbox0=\hbox{#1}{\ooalign{\hidewidth
  \lower1.5ex\hbox{`}\hidewidth\crcr\unhbox0}}}
  \def\polhk#1{\setbox0=\hbox{#1}{\ooalign{\hidewidth
  \lower1.5ex\hbox{`}\hidewidth\crcr\unhbox0}}}
  \def\polhk#1{\setbox0=\hbox{#1}{\ooalign{\hidewidth
  \lower1.5ex\hbox{`}\hidewidth\crcr\unhbox0}}} \def\cprime{$'$}
  \def\polhk#1{\setbox0=\hbox{#1}{\ooalign{\hidewidth
  \lower1.5ex\hbox{`}\hidewidth\crcr\unhbox0}}}
  \def\polhk#1{\setbox0=\hbox{#1}{\ooalign{\hidewidth
  \lower1.5ex\hbox{`}\hidewidth\crcr\unhbox0}}}
  \def\polhk#1{\setbox0=\hbox{#1}{\ooalign{\hidewidth
  \lower1.5ex\hbox{`}\hidewidth\crcr\unhbox0}}}
  \def\polhk#1{\setbox0=\hbox{#1}{\ooalign{\hidewidth
  \lower1.5ex\hbox{`}\hidewidth\crcr\unhbox0}}}


\begin{thebibliography}{10}

\bibitem{Arnold88}
V.~I. Arnol{\cprime}d, V.~V. Kozlov, and A.~I. Ne{\u\i}shtadt.
\newblock {\em Dynamical {S}ystems. {III}}, volume~3 of {\em Encyclopaedia of
  Mathematical Sciences}.
\newblock Springer-Verlag, Berlin, 1988.
\newblock Translated from the Russian by A. Iacob.

\bibitem{BaSn93}
L.~Bates and J.~\'Sniatycki.
\newblock {Nonholonomic reduction.}
\newblock {\em Rep. Math. Phys.}, {\bf 32}, 99 (1993).

\bibitem{Blankenstein00}
G.~Blankenstein.
\newblock {\em Implicit {H}amiltonian {S}ystems: {S}ymmetry and
  {I}nterconnection}.
\newblock Ph.D.Thesis, University of Twente, 2000.

\bibitem{BlRa04}
G.~Blankenstein and T.S. Ratiu.
\newblock Singular reduction of implicit {H}amiltonian systems.
\newblock {\em Rep. Math. Phys.}, {\bf 53}, 211 (2004).

\bibitem{BlvdS01}
G.~Blankenstein and A.J. van~der Schaft.
\newblock Symmetry and reduction in implicit generalized {H}amiltonian systems.
\newblock {\em Rep. Math. Phys.}, {\bf 47}, 57 (2001).

\bibitem{Bloch03}
A.~M. Bloch.
\newblock {\em Nonholonomic {M}echanics and {C}ontrol}, volume~24 of {\em
  Interdisciplinary Applied Mathematics}.
\newblock Springer-Verlag, New York, 2003.
\newblock With the collaboration of J. Baillieul, P. Crouch and J. Marsden,
  With scientific input from P. S. Krishnaprasad, R. M. Murray and D. Zenkov,
  Systems and Control.

\bibitem{BuCaGu07}
H.~Bursztyn, G.~R. Cavalcanti, and M.~Gualtieri.
\newblock {Reduction of Courant algebroids and generalized complex structures.}
\newblock {\em Adv. Math.}, {\bf 211}, 726 (2007).

\bibitem{BuCr05}
H.~Bursztyn and M.~Crainic.
\newblock Dirac structures, momentum maps, and quasi-{P}oisson manifolds.
\newblock In {\em The breadth of symplectic and {P}oisson geometry}, volume 232
  of {\em Progr. Math.}, pages 1--40. Birkh\"auser Boston, Boston, MA, 2005.

\bibitem{BuCrWeZh04}
H.~Bursztyn, M.~Crainic, A.~Weinstein, and C.~Zhu.
\newblock Integration of twisted {D}irac brackets.
\newblock {\em Duke Math. J.}, {\bf 123}, 549 (2004).

\bibitem{CaDiLeMa98}
F.~Cantrijn, D.~M. de~Diego, M.~de~Le{\'o}n, and J.~C. Marrero.
\newblock Reduction of nonholonomic mechanical systems with symmetries.
\newblock {\em Rep. Math. Phys.}, {\bf 42}, 25 (1998).
\newblock Pacific Institute of Mathematical Sciences Workshop on Nonholonomic
  Constraints in Dynamics (Calgary, AB, 1997).

\bibitem{CeMaRaYo08}
H.~Cendra, J.~E. Marsden, T.S. Ratiu, and H.~Yoshimura.
\newblock {\em In preparation}, 2008.

\bibitem{Courant90a}
T.~J. Courant.
\newblock {Dirac manifolds.}
\newblock {\em Trans. Am. Math. Soc.}, {\bf 319}, 631 (1990).

\bibitem{CoWe88}
T.~J. Courant and A.~Weinstein.
\newblock Beyond {P}oisson structures.
\newblock In {\em Action hamiltoniennes de groupes. {T}roisi\`eme th\'eor\`eme
  de {L}ie ({L}yon, 1986)}, volume~27 of {\em Travaux en Cours}, pages 39--49.
  Hermann, Paris, 1988.

\bibitem{BaCuKeSn95}
R.~Cushman, D.~Kemppainen, J.~{\'S}niatycki, and L.~Bates.
\newblock Geometry of nonholonomic constraints.
\newblock In {\em Proceedings of the XXVII Symposium on Mathematical Physics
  (Toru\'n, 1994)}, volume~36, pages 275--286 (1995).

\bibitem{FaRaSa07}
F.~Fass{\`o}, A.~Ramos, and N.~Sansonetto.
\newblock The reaction-annihilator distribution and the nonholonomic {N}oether
  theorem for lifted actions.
\newblock {\em Regul. Chaotic Dyn.}, {\bf 12}, 579 (2007).

\bibitem{JoRaZa11}
M.~Jotz, T.~Ratiu, and M.~Zambon.
\newblock Invariant frames for vector bundles and applications.
\newblock To appear in {\em Geometriae Dedicata} (2011).

\bibitem{JoRa10b}
M.~Jotz and T.S. Ratiu.
\newblock {Optimal Dirac reduction}.
\newblock {\em Preprint, arXiv:1008.2283}, 2010.

\bibitem{JoRa11}
M.~Jotz and T.S. Ratiu.
\newblock {Induced Dirac structure on isotropy type manifolds}.
\newblock {\em Transform. Groups}, {\bf 16}, 175 (2011).

\bibitem{JoRaSn11}
M.~Jotz, T.S. Ratiu, and J.~{\'S}niatycki.
\newblock Singular {D}irac reduction.
\newblock {\em Trans. Amer. Math. Soc.}, {\bf 363}, 2967 (2011).

\bibitem{KoKo78}
V.~V. Kozlov and N.~N. Kolesnikov.
\newblock On theorems of dynamics.
\newblock {\em J. Appl. Math. Mech.}, {\bf 42}, 28 (1978).

\bibitem{LiMa87}
P.~Libermann and C.-M. Marle.
\newblock {\em {Symplectic {G}eometry and {A}nalytical {M}echanics. Transl.
  from the French by Bertram Eugene Schwarzbach.}}
\newblock {Mathematics and its Applications, 35. Dordrecht etc.: D. Reidel
  Publishing Company, a member of the Kluwer Academic Publishers Group. XVI,
  526 p. }, 1987.

\bibitem{LiWeXu97}
Z.-J. Liu, A.~Weinstein, and P.~Xu.
\newblock Manin triples for {L}ie bialgebroids.
\newblock {\em J. Differential Geom.}, {\bf 45}, 547 (1997).

\bibitem{MaRa99}
J.~E. Marsden and T.S. Ratiu.
\newblock {\em {Introduction to {M}echanics and {S}ymmetry. A {B}asic
  {E}xposition of {C}lassical {M}echanical {S}ystems. 2nd ed.}}
\newblock {Texts in Applied Mathematics. 17. New York, NY: Springer. xviii, 582
  p. }, 1999.

\bibitem{OrRa04}
J.-P. Ortega and T.S. Ratiu.
\newblock {\em {Momentum {M}aps and {H}amiltonian {R}eduction.}}
\newblock {Progress in Mathematics (Boston, Mass.) 222. Boston, MA:
  Birkh\"auser. xxxiv, 497~p. }, 2004.

\bibitem{Pflaum01}
M.J. Pflaum.
\newblock {\em Analytic and {G}eometric {S}tudy of {S}tratified {S}paces},
  volume 1768 of {\em Lecture Notes in Mathematics}.
\newblock Springer-Verlag, Berlin, 2001.

\bibitem{Rosenberg77}
R.~M. Rosenberg.
\newblock {\em Analytical {D}ynamics of {D}iscrete {S}ystems}, volume~4 of {\em
  Mathematical Concepts and Methods in Science and Engineering}.
\newblock Plenum Press, New York, 1977.

\bibitem{Schwarz80}
G.W. Schwarz.
\newblock Lifting smooth homotopies of orbit spaces.
\newblock {\em Inst. Hautes \'Etudes Sci. Publ. Math.}, {\bf 51}, 37 (1980).

\bibitem{SeWe01}
P.~{\v{S}}evera and A.~Weinstein.
\newblock Poisson geometry with a 3-form background.
\newblock {\em Progr. Theoret. Phys. Suppl.}, {\bf 144}, 145 (2001).
\newblock Noncommutative geometry and string theory (Yokohama, 2001).

\bibitem{Stefan74b}
P.~Stefan.
\newblock Accessibility and foliations with singularities.
\newblock {\em Bull. Amer. Math. Soc.}, {\bf 80}, 1142 (1974).

\bibitem{Stefan74a}
P.~Stefan.
\newblock Accessible sets, orbits, and foliations with singularities.
\newblock {\em Proc. London Math. Soc. (3)}, {\bf 29}, 699 (1974).

\bibitem{Stefan80}
P.~Stefan.
\newblock Integrability of systems of vector fields.
\newblock {\em J. London Math. Soc. (2)}, {\bf 21}, 544 (1980).

\bibitem{StXu08}
M.~Sti{\'e}non and P.~Xu.
\newblock Reduction of generalized complex structures.
\newblock {\em J. Geom. Phys.}, {\bf 58}, 105 (2008).

\bibitem{Sussmann73}
H.J. Sussmann.
\newblock Orbits of families of vector fields and integrability of
  distributions.
\newblock {\em Trans. Amer. Math. Soc.}, {\bf 180}, 171 (1973).

\bibitem{Vaisman94}
I.~Vaisman.
\newblock {\em Lectures on the {G}eometry of {P}oisson {M}anifolds}, volume 118
  of {\em Progress in Mathematics}.
\newblock Birkh\"auser Verlag, Basel, 1994.

\bibitem{MaYo06}
H.~Yoshimura and J.~E. Marsden.
\newblock Dirac structures in {L}agrangian mechanics. {I}. {I}mplicit
  {L}agrangian systems.
\newblock {\em J. Geom. Phys.}, {\bf 57}, 133 (2006).

\bibitem{MaYo07}
H.~Yoshimura and J.~E. Marsden.
\newblock Reduction of {D}irac structures and the {H}amilton-{P}ontryagin
  principle.
\newblock {\em Rep. Math. Phys.}, {\bf 60}, 381 (2007).

\bibitem{MaYo09}
H.~Yoshimura and J.E. Marsden.
\newblock Dirac cotangent bundle reduction.
\newblock {\em J. Geom. Mech.}, {\bf 1}, 87 (2009).

\end{thebibliography}
\end{document}